\documentclass[twoside,12pt,a4paper]{amsart}
\usepackage{amsfonts,amsmath,amssymb,amsthm}

\usepackage{color}
\usepackage{graphicx}
\usepackage{latexsym}

\usepackage{caption}
\usepackage{subfigure}
\usepackage{cite}
\usepackage{subfigure}
\usepackage{float}
\usepackage{paralist}
\usepackage{indentfirst}
\usepackage{bm}
\usepackage{tikz}

\usepackage{soul}

\usepackage[hidelinks]{hyperref}

\usepackage[version=4]{mhchem}

\newcommand{\N}{\mathbb{N}}
\newcommand{\R}{\mathbb{R}}

\newcommand{\mE}{\mathcal{E}}
\newcommand{\dive}{\mathop{\mathrm{div}}}
\newtheorem{theorem}{Theorem}[section]
\newtheorem{remark}{Remark}[section]

\newtheorem{lemma}{Lemma}[section]
\newtheorem{prop}{Proposition}[section]

\numberwithin{equation}{section}
\def\al{\alpha}

\def\na{\nabla}
\def\pa{\partial}

\subjclass[2010]{35Q20; 76N10}

\keywords{Compressible Euler system; Compressible Euler-Maxwell system; Relaxation limit; Convergence rate; Stream function}

\begin{document}

\title[Relaxation limits of the compressible Euler and Euler-Maxwell systems]{Global convergence rates in the relaxation limits for the compressible Euler and Euler-Maxwell systems in Sobolev spaces}

\author[T. Crin-Barat]{Timothée Crin-Barat}
\address[T. Crin-Barat]{Université de Toulouse,  Institut de Math\'ematiques de Toulouse, 118 Route de Narbonne, 31062 Toulouse Cedex 9, France}
\email{timothee.crin-barat@math.univ-toulouse.fr}

\author[Y.-J. Peng]{Yue-Jun Peng}
\address[Y.-J. Peng]{Laboratoire de Mathématiques Blaise Pascal,
Université Clermont Auvergne / CNRS
63178 Aubière Cedex, France}
\email{yue-jun.peng@uca.fr}

\author[L.-Y. Shou]{Ling-Yun Shou\textsuperscript{*}}
\address[L.-Y. Shou]{School of Mathematical Sciences {\color{blue}and} Ministry of Education Key Laboratory of NSLSCS, Nanjing Normal University, Nanjing, 210023,
P. R. China}
\email{shoulingyun11@gmail.com}

\footnotetext[1]{\textsuperscript{*}Corresponding author: L.-Y. Shou (shoulingyun11@gmail.com)}

\date{}

\begin{abstract}
We study two relaxation problems in the class of partially dissipative hyperbolic systems: the compressible Euler system and the compressible Euler-Maxwell system. In classical Sobolev spaces, we derive a global convergence rate of $\mathcal{O}(\varepsilon)$ between strong solutions of the relaxed Euler system and the porous medium equation in $\mathbb{R}^d$ ($d\geq1$) for \emph{ill-prepared} initial data. In a well-prepared setting, we derive an enhanced convergence rate of order $\mathcal{O}(\varepsilon^2)$ between the solutions of the relaxed compressible Euler system and their first-order asymptotic approximation. Regarding the relaxed Euler-Maxwell system, we prove the global strong convergence of its solutions to the drift-diffusion model in $\mathbb{R}^3$ in an \emph{ill-prepared} setting. These results are achieved by developing a new asymptotic expansion approach that, combined with stream function techniques, ensures uniform-in-time error estimates.
\end{abstract}

\maketitle

\section{Introduction}
\vspace{2mm}

Relaxation phenomena occur in a wide variety of physical situations, such as modeling blood flow with friction forces, non-equilibrium gas dynamics, kinetic theory, traffic flows, and more (see \cite{cere,vicenti1,Whitham}). These phenomena describe situations in which a physical system is perturbed away from a stable equilibrium state. In this case, systems are described by a set of equations where the source terms contain a large coefficient 
$\varepsilon^{-1}$ with $\varepsilon>0$ representing a time-relaxation parameter. As $\varepsilon$ tends to zero, the solutions are expected to relax towards equilibrium.

In this paper, we consider two classical hyperbolic systems with relaxation in the whole space: the damped compressible Euler system and the damped compressible Euler-Maxwell system. In both cases, global-in-time asymptotic stability as the time-relaxation parameter $\varepsilon$ goes to zero is established in Sobolev spaces.

\subsection{The compressible Euler system with damping}
First, we consider the Euler equations with relaxation in several space variables:
\begin{equation}
\label{Euler}
\left\{
\begin{aligned}
  &\partial_t\rho+\dive(\rho v)=0,\\
  &\partial_t(\rho v)+\dive(\rho v\otimes v)+\na p(\rho)=-\dfrac{\rho v}{\varepsilon},
  \end{aligned}
  \right.
  \end{equation}
with the initial condition
\begin{align}
    (\rho,v)(0,x)=\big(\rho^\varepsilon_0,v^\varepsilon_0\big)(x),\label{Eulerd}
\end{align}
for the time $t>0$ and the position $x=(x_1,x_2,\cdots,x_d)\in\mathbb{R}^d$.
Here, $\varepsilon\in (0,1]$ is the relaxation time, $\rho=\rho(t,x)>0$,
$v=(v^1,v^2,\cdots,v^d)(t,x)$ and $p=p(t,x)$ are
the density, the velocity and the pressure functions, respectively. 


We are interested in strong solutions that are small perturbations of the constant equilibrium
$$
(\bar{\rho},\bar{v})=(\bar{\rho},0), 
$$
where $\bar{\rho}$ is a given positive constant.
Additionally, we assume that the pressure is a function of the density $\rho$,
\begin{equation}
\label{F1.1}
    p\in C^{\infty}_{{\rm{loc}}}(\mathbb{R}^+)\quad \text{and}\quad  p'(\bar{\rho})>0.
\end{equation}
Without loss of generality, we assume that $\bar{\rho}=1$ in the following. For this type of solutions, we will analyze their diffusion limit corresponding to the scaling
\[  \big(\rho^\varepsilon,u^\varepsilon\big)(t,x)
:=\left(\rho,\dfrac{1}{\varepsilon}v\right)\left(\dfrac{t}{\varepsilon},x\right).  \]
The diffusively rescaled Euler system and the initial condition \eqref{Eulerd} become
\begin{equation}
\label{F1.2}
\begin{cases}
  \partial_t\rho^\varepsilon+\dive (\rho^\varepsilon u^\varepsilon)=0,\\[1mm]
\varepsilon^2 \partial_t(\rho^\varepsilon u^\varepsilon)+\varepsilon^2 \dive(\rho^\varepsilon u^\varepsilon\otimes u^\varepsilon)
+\na p(\rho^\varepsilon)=-\rho^\varepsilon u^\varepsilon,
\end{cases}
\end{equation}
and
\begin{equation}
\label{F1.3}
  \big(\rho^\varepsilon,u^\varepsilon\big)(0,x)
=\left(\rho^\varepsilon_0,u_0^\varepsilon\right)(x)\quad \text{with}\quad u_0^\varepsilon:=\dfrac{1}{\varepsilon}v^\varepsilon_0.
\end{equation}
Formally, let us denote
\[  \big(\rho^*,u^*\big):=\lim_{\varepsilon\to 0}\big(\rho^\varepsilon,u^\varepsilon\big).
\]
As $\varepsilon\to 0$, we expect that
\begin{equation}
\label{F1.5}
\begin{cases}
   \partial_t\rho^*+\dive(\rho^*u^*)=0,\\
   \na p(\rho^*)=-\rho^*u^*\quad\text{(Darcy's law)},
\end{cases}
\end{equation}
which yields a filtration equation for the density
\begin{equation}
\label{F1.6}
   \partial_t\rho^*-\Delta p(\rho^*)=0,
\end{equation}
subject to the initial condition
\begin{equation}
\label{F1.4}
    \rho^*(0,x)=\rho^*_0(x).
\end{equation}
In particular, for a $\gamma$-law pressure
$p(\rho)=a^2\rho^\gamma$ with $a>0$ and $\gamma>1$, we recover
the porous medium equation. 


\subsection{The compressible Euler-Maxwell system}
Next, we consider the three-dimensional compressible Euler-Maxwell system:
\begin{equation}
\label{EM}
\begin{cases}
  \partial_t\rho+\dive (\rho v)=0,\\[1mm]
\partial_t(\rho v)+\dive(\rho v\otimes v)
+\na p(\rho)
=-\rho\big(E+ v\times B\big)-\dfrac{\rho v}{\varepsilon},\\
\partial_t E-\na\times B=\rho v,\\
\partial_t B+\na\times E=0,\\
\dive E=1-\rho,\quad\quad \dive B=0,
\end{cases}
\end{equation}
where the unknowns $\rho=\rho(t,x)>0$, $v=(v^1,v^2,v^3)(t,x)$, $E=(E^1,E^2,E^3)(t,x)$ and $B=(B^1,B^2,B^3)(t,x)$ stand for
the density, the velocity, the electric field and the magnetic field, respectively, at the time $t>0$ and the position $x=(x_1,x_2,x_3)$. The term $\rho(E+v\times B)$ corresponds to the Lorentz force, and $p=p(\rho)$ denotes the pressure function that satisfies \eqref{F1.1}. 

The system \eqref{EM} is a fundamental model for plasma dynamics, describing the interactions between compressible ion and electron fluids with a self-consistent electromagnetic field, see \cite{Bitt, CJW,GuoIon}.  We supplement the system \eqref{EM} with the initial condition
\begin{align}
(\rho,v,E,B)(0,x)=(\rho_0^\varepsilon,  v_0^\varepsilon,E_0^\varepsilon,B_0^\varepsilon)(x),\label{EMd}
\end{align}
and we assume that
\begin{align}
    \dive E_0^\varepsilon=1-\rho^\varepsilon_0,\quad\quad \dive B^\varepsilon_0=0.\label{EMd2}
\end{align}
Note that under \eqref{EMd2}, the constraint relaxations in $\eqref{EM}_5$ hold true for all $(t,x)\in\mathbb{R}^+\times\mathbb{R}^3$.

We are interested in strong solutions for \eqref{EM} that are small perturbations of the constant equilibrium
$$   (\bar{\rho},\bar v,\bar E,\bar B)=(1,0,0,B^e)   $$
where $B^e\in\mathbb{R}^3$ is a constant vector. For these solutions, we will
investigate their diffusion limit corresponding to the scaling
\[
\big(\rho^\varepsilon,u^\varepsilon,E^\varepsilon,B^\varepsilon\big)(t,x)
:=\left(\rho,\dfrac{1}{\varepsilon}v,E,B\right)\left(\dfrac{t}{\varepsilon},x\right).  \]
Then, the diffusively rescaled Euler-Maxwell system reads
\begin{equation}
\label{EMeps}
\begin{cases}
  \partial_t\rho^\varepsilon+\dive (\rho^\varepsilon u^\varepsilon)=0,\\[1mm]
\varepsilon^2\partial_t(\rho^\varepsilon u^\varepsilon)+\varepsilon^2\dive(\rho^\varepsilon u^\varepsilon\otimes u^\varepsilon)
+\na p(\rho^\varepsilon)
=-\rho^\varepsilon\big(E^\varepsilon+\varepsilon u^\varepsilon\times B^\varepsilon\big)-\rho^\varepsilon u^\varepsilon,\\
\varepsilon\partial_t E^\varepsilon-\na\times B^\varepsilon=\varepsilon\rho^\varepsilon u^\varepsilon,\\
\varepsilon\partial_t B^\varepsilon+\na\times E^\varepsilon=0,\\
\dive E^\varepsilon=1-\rho^\varepsilon,\qquad \dive B^\varepsilon=0,
\end{cases}
\end{equation}
supplemented with the initial data
\begin{equation}
\label{Emepsd}
   \big(\rho^\varepsilon,u^\varepsilon,E^\varepsilon,B^\varepsilon\big)(0,x)
=\left(\rho^\varepsilon_0,u^\varepsilon_0,E^\varepsilon_0,B^\varepsilon_0\right)(x)\quad \text{with}\quad u_0^\varepsilon:=\dfrac{1}{\varepsilon}v^\varepsilon_0.
\end{equation}
Setting
\[  \big(\rho^*,u^*,E^*,B^*\big):=\lim_{\varepsilon\to 0}\big(\rho^\varepsilon,u^\varepsilon,E^\varepsilon,B^\varepsilon\big),
\]
formally, as $\varepsilon\to 0$, the dynamics of \eqref{EMeps} are expected to be governed by
\begin{equation}
\label{F3.3}
\begin{cases}
\partial_t\rho^*+\dive(\rho^*u^*)=0,\\
  \rho^*u^*=- \na p(\rho^*)-\rho^* E^* \quad\text{(Darcy-type law)},\\
   \na\times B^*=0, \\   \na\times E^*=0,\\
   \dive E^*=1-\rho^*,\qquad \dive B^*=0.
\end{cases}
\end{equation}
The divergence-free and curl-free conditions imply that $B^*$ is a constant vector. Moreover, setting $\Lambda^{-2}=(-\Delta)^{-1}$, there exists a potential function $\phi^*=\Lambda^{-2}(\rho^*-1)$ such that
\[
E^*=\na\phi^*=\Lambda^{-2}\na\rho^*.
\]It follows that
$(\rho^*,\phi^*)$ solves the drift-diffusion system
\begin{equation}
\label{DD}
\begin{cases}
   \partial_t\rho^*-\dive\big(\na p(\rho^*)+\rho^*\na\phi^*\big)=0,\\
   \Delta\phi^*=1-\rho^*,
\end{cases}
\end{equation}
subject to the initial condition
\begin{align}
\rho^*(0,x)=\rho^*_0(x).  \label{DDd}
\end{align}

\subsection{Literature on damped Euler equations and partially dissipative hyperbolic systems}

Before stating the paper's findings, we recall recent efforts devoted to studying partially dissipative hyperbolic systems of the type:
\begin{equation}
  \partial_t V+\sum_{j=1}^dA^j(V)\partial_{x_j}V=\frac{H(V)}{\varepsilon},
  \label{GEQSYM}
\end{equation}
where the unknown $V=V(t,x)$ is a $N$-vector valued function depending on
the time variable $t\in\mathbb{R}^{+}$ and on the space variable
$x\in\mathbb{R}^d\; (d\geq1).$ The $A^j(V)$ ($j=1,..,d$) and $H$ are given
smooth functions on a state space $\mathcal{O}_{V}\in\mathbb{R}^N$.
In the absence of source term $H(V)$, \eqref{GEQSYM} reduces to
a system of conservation laws. In that case, it is well-known that
classical solutions may develop singularities (\textit{e.g.}, shock waves)
in finite time, even if initial data are smooth and small
(see Dafermos \cite{Dafermos1} and Lax \cite{Lax1}). 

A typical example is the isentropic compressible Euler equations with damping \eqref{Euler}. For initial data being small perturbations in Sobolev spaces $H^{s}(\mathbb{R}^{d})$ ($s\geq [\frac{d}{2}]+2$), the global well-posedness and asymptotics of classical solutions for \eqref{Euler} have been studied by Wang and Yang in \cite{Wang}, as well as Sideris, Thomases and Wang in \cite{Sideris}. In the case $\varepsilon=1$, a natural question is to identify the conditions  imposed on $H(V)$ to prevent the finite-time blow-up of classical solutions for \eqref{GEQSYM}. Chen, Levermore and Liu in \cite{chen1994} first formulated a notion of the entropy for \eqref{GEQSYM}. Yong in \cite{Yong} proved the global existence of classical solutions in a neighborhood of constant equilibrium $\bar{V}\in \mathbb{R}^{N}$ satisfying $H(\bar{V})=0$ under the Shizuta–Kawashima condition \cite{SK}. We also mention that Hanouzet and Natalini in \cite{HanouzetNatalini} obtained a similar global existence result for the one-dimensional problem before the work \cite{Yong}. Subsequently, Kawashima and Yong in \cite{KY} gave a new definition of the entropy notion and removed the technical requirement on the dissipative entropy in \cite{Yong,HanouzetNatalini}. Then, Bianchini, Hanouzet and Natalini  in \cite{BHN} derived asymptotic behaviors of smooth solutions to the Cauchy problem for \eqref{GEQSYM}. Recently, Beauchard and Zuazua in \cite{BZ} established the equivalence of the Shizuta-Kawashima condition and the Kalman rank condition from control theory. Thanks to this observation, they provided a new proof for the linear stability using a hypocoercive approach. Then, extensions to critical spaces were obtained in non-homogeneous settings in \cite{XK1,XK2} and in some hybrid homogeneous settings in \cite{CBD2,CBD1,CBD3}.

Regarding the relaxation limit as $\varepsilon\to0$ in systems of the type \eqref{GEQSYM}, the first justification is due to Marcati, Milani and Secchi in \cite{MMS} in a one-dimensional setting. The limiting procedure was carried out by using the theory of compensated compactness. Then, Liu  in \cite{Liu} proved the relaxation to parabolic equations for genuinely nonlinear hyperbolic systems. Marcati and Milani in \cite{mar0} considered the diffusive scaling for the one-dimensional compressible Euler flow
\eqref{Euler} and derived Darcy's law as $\varepsilon\rightarrow0$, which is analogous to the one derived in \cite{MMS}. Later,  Marcati and Rubino in \cite{marcatiparabolicrelaxation}
developed a complete hyperbolic to parabolic relaxation theory for $2\times 2$ genuinely nonlinear hyperbolic balance laws. Junca and Rascle in \cite{Junca} established the relaxation convergence from the isothermal equation \eqref{Euler} to the heat equation for arbitrarily large initial data in $BV(\mathbb{R})$ that are bounded away from vacuum. As for \eqref{GEQSYM} in several dimensions,
Coulombel, Goudon and Lin in  \cite{CoulombelGoudon,CoulombelLin} employed the classical energy approach and constructed uniform-in-$\varepsilon$ smooth solutions to the isentropic
Euler equations \eqref{Euler} and justified the weak relaxation limit in the Sobolev spaces $H^{s}(\mathbb{R}^{d})\:(s\geq [\tfrac{d}{2}]+2)$.
Xu and Wang in \cite{XuWang} improved their result to the setting of critical inhomogeneous Besov space. More precisely, it is shown that the density converges towards the solution of the porous medium equation, as $\varepsilon\rightarrow0$. Peng and Wasiolek in \cite{Peng16AIHP} established the  uniform local existence with respect to $\varepsilon$ and the convergence of general systems \eqref{GEQSYM} to parabolic-type equations as $\varepsilon\rightarrow0$. Subsequently, under the Shizuta–Kawashima stability condition, the authors of  \cite{PengWasiolek} established the uniform global existence and the global-in-time weak convergence from \eqref{GEQSYM} to second-order nonlinear parabolic systems by using Aubin-Lions compactness arguments. 

There are few results deriving explicit convergence rates for the global relaxation limit. In the spirit of the stream function approach of \cite{Junca}, Li, Peng and Zhao \cite{LiPengZhao1d} obtained explicit convergence rates for this relaxation process for $d=1$. Recently, the first author and Danchin \cite{CBD3} justified the strong relaxation limit of diffusively rescaled solutions for \eqref{GEQSYM} globally in time in homogeneous critical Besov spaces with an explicit convergence rate. Finally, we also mention the work of Peng  \cite{Pe24} on this limit for large smooth solutions of the isothermal Euler equations.
\subsection{Literature on the Euler-Maxwell system}

In contrast to the relaxed Euler equations \eqref{Euler}, the Euler-Maxwell system \eqref{EM} exhibits a non-symmetric relaxation mechanism induced by the Lorentz force and Maxwell equations, leading to different behaviors of the solutions and mathematical challenges.

So far there are several results regarding the global existence, large-time behaviour and asymptotic limit for the isentropic Euler-Maxwell system \eqref{EM}.  Chen, Jerome and Wang in \cite{CJW} constructed global
weak solutions to the initial boundary value problem for arbitrarily large initial data in one spatial dimension. In the multidimensional case, the question of global weak solutions is quite open, and mainly smooth solutions have been studied. Jerome in \cite{jerome1}
established the local well-posedness of smooth solutions to the Cauchy problem \eqref{EM}-\eqref{EMd} in the framework of Sobolev spaces $H^{s}(\mathbb{R}^{d})$ with $s\geq 3$ according to the standard theory for symmetrizable
hyperbolic systems. The existence of global smooth solutions near constant equilibrium states has been obtained independently
by Peng, Wang and Gu in \cite{PYJEMSIAM}, and Duan in \cite{Duan2011} in Sobolev spaces and by Xu in \cite{XuEM} in the inhomogeneous critical Besov space. Ueda, Wang and Kawashima in \cite{UedaWangKawa2012} pointed out that the system \eqref{EM} was of regularity-loss type, and time-decay estimates were derived in \cite{Duan2011,UedaKawa2011}. Later, Ueda, Duan and Kawashima in \cite{UDK1} formulated a new structural condition to analyze the weak dissipative mechanism for general hyperbolic systems with non-symmetric relaxation (including the Euler-Maxwell system \eqref{EM}).
Xu, Mori and Kawashima  in \cite{xumorikaw1} developed a time-decay inequality with weaker regularity assumptions of initial data. Concerning the stability of steady-states, we refer to \cite{pengzhu14,peng2015,liuguopeng19}.
In the absence of a damping term in \eqref{EM}, there are results concerning the global or large-time existence using dispersive properties; see \cite{GMEM1,GuoIon,DAI1,IonLie}.

Concerning the relaxation from \eqref{EM} to \eqref{DD}, Hajjej and Peng in \cite{hajjejpeng12} carried out an asymptotic expansion and obtained convergence rates for the relaxation procedure in the case of local-in-time solutions for both well-prepared data and ill-prepared data. Wasiolek in \cite{Wasiolek2016} obtained uniform estimates for global solutions with small perturbations in $H^{s}(\Omega)$ with $s\geq3$ and $\Omega=\mathbb{R}^3$ or $\mathbb{T}^3$, and proved the global weak convergence of the solutions of \eqref{EM} to the solutions of \eqref{DD} via the compactness argument. Recently, Li, Peng and Zhao in \cite{liEM} studied the relaxation limit for global smooth solutions in $\mathbb{T}^3$ and obtained error estimates for smooth periodic solutions between \eqref{EM} and \eqref{DD} using stream function techniques and Poincaré's inequality. Very recently, by developing a new characterization of the dissipation structure and using Fourier analysis tools, the authors of \cite{CBPSX25} provided a rigorous justification of the global-in-time strong convergence of the relaxation process in $\mathbb{R}^3$ with an explicit convergence rate.

\subsection{Aims of the paper}

The relaxation limit problems for both the compressible Euler system \eqref{F1.2} and the compressible Euler–Maxwell system \eqref{EMeps} in the whole space have been addressed in only a few references. On the real line $\mathbb{R}$, the strong relaxation limit was justified in \cite{Junca,LiPengZhao1d} by employing stream function techniques. In higher dimensions, relaxation limits were established in critical Besov spaces in \cite{CBD3,CBPSX25}, relying heavily on elaborate frequency partitions of the solution and the Littlewood–Paley theory. This naturally raises the following question: \textit{Can global-in-time error estimates for such high-dimensional relaxation problems be derived within the classical Sobolev energy framework}?


In the present paper, we answer this question in the affirmative.  In a Sobolev setting, we establish uniform error estimates for ill-prepared initial data. In addition, we prove the convergence in the natural energy space  $L^2(\mathbb{R}^+;L^2)$, which cannot be directly recovered by the results in \cite{CBD3,CBPSX25}. Our proof relies on a high-dimensional version of the stream function techniques and an asymptotic expansion method. The framework we develop avoids using frequency-dependent tools to broaden its  applicability. Specifically, we expect our method to be suitable for cases where the physical domain is a half space or a bounded domain, and for justifying such limits in discrete settings.

\section{Main results}

\subsection{The compressible Euler system}

We first investigate the relaxation problem for the compressible Euler equations \eqref{F1.2}. We define the initial energy
\begin{equation}
\label{F1.7}
  \mE^\varepsilon_0:=\|\rho^\varepsilon_0-1\|^2_{H^m}+\varepsilon^2\|u^\varepsilon_0\|^2_{H^m},
\end{equation}
and the variables for the momentum
\[   q^\varepsilon:=\rho^\varepsilon u^\varepsilon,\qquad q^*:=\rho^*u^*,\qquad
q^\varepsilon_0:=\rho^\varepsilon_0 u^\varepsilon_0. 
\]


We recall the classical uniform global-in-time existence result for \eqref{F1.2}-\eqref{F1.3} from \cite{CoulombelGoudon,CoulombelLin}.
\begin{prop}
{\rm(\!\!\cite{CoulombelGoudon,CoulombelLin})}
\label{P1.1}
Let $d\geq1$, $m\geq [\frac{d}{2}]+2$ and \eqref{F1.1}
hold. There exists a positive constant $\delta$, independent of $\varepsilon$, such
that, if $\mE^\varepsilon_0\le\delta$, then the Cauchy problem \eqref{F1.2}-\eqref{F1.3}
admits a unique global classical solution $(\rho^\varepsilon,u^\varepsilon)$ with $q^\varepsilon=\rho^\varepsilon u^\varepsilon$ which satisfies $(\rho^\varepsilon-1,q^\varepsilon)\in C(\mathbb{R}^{+};H^m)$ and the uniform estimate
\begin{eqnarray}
\label{F1.8}
&\,&\sup_{t\in\mathbb{R}^+}\big(\|\rho^\varepsilon(t)-1\|^2_{H^m}+\varepsilon^2\|u^\varepsilon(t)\|^2_{H^m}+\varepsilon^2\|q^\varepsilon(t)\|^2_{H^m}\big)\\
&\,&+\,\int_0^{+\infty}\big(\|\na\rho^\varepsilon(t)\|^2_{H^{m-1}}
+\|u^\varepsilon(t)\|^2_{H^m}+\|q^\varepsilon(t)\|^2_{H^m}\big)\,dt\nonumber
\le C \mE^\varepsilon_0,
\end{eqnarray}
for a generic constant $C>0$. Moreover, for any given $T>0$, as $\varepsilon\to 0$, up to subsequences,
\begin{equation}
\left\{
\begin{aligned}
    \rho^\varepsilon&\rightarrow \rho^*\;\; \mbox{strongly in}\;\,\,\, 
C([0,T],H^{m-1}_{loc}), \\
q^\varepsilon &\rightharpoonup q^* \;\;\mbox{weakly \,\, in}~~~\;\;
L^2(\mathbb{R}^+,H^m),
\end{aligned}
\right.
\end{equation}
where $(\rho^*,u^*)$ with $u^*=q^*/\rho^*$ is the unique solution of
\eqref{F1.5}-\eqref{F1.4} with  $\rho^*_0$ denoting the weak limit of $(\rho^\varepsilon_0)_{\varepsilon>0}$ in $H^m$.
\end{prop}

For the filtration equation \eqref{F1.6}, one has the following existence result.

\begin{prop}
\label{P1.2}
    Let $d\geq1$, $m\geq [\frac{d}{2}]+2$ and \eqref{F1.1}
hold. There exists a positive constant $\delta^*$ such that,
if $\|\rho^*_0-1\|^2_{H^m} \le\delta^*$, then the Cauchy problem \eqref{F1.6}-\eqref{F1.4}
admits a unique global classical solution $\rho^*$ which
satisfies
\begin{equation}
\label{F1.9}
   \sup_{t\in\mathbb{R}^+}\|\rho^*(t)-1\|^2_{H^m}+\int_0^{+\infty}\|\na\rho^*(t)\|^2_{H^m}\,dt
\le C\|\rho^*_0-1\|^2_{H^m},
\end{equation}
where $C>0$ is a generic constant.
Furthermore, let $u^*$ be given by Darcy's law $\eqref{F1.5}_2$ and set $q^*=\rho^*u^*$. It holds that
\begin{equation}
\label{F1.91}
   \sup_{t\in\mathbb{R}^+}\|q^*(t)\|^2_{H^{m-1}}+\int_0^{+\infty} \|q^*(t)\|^2_{H^m}\,dt
\le C\|\rho^*_0-1\|^2_{H^m}.
\end{equation}
\end{prop}

When studying relaxation limit problems for the rescaled Euler system \eqref{F1.2}, two different heuristics ({\emph{ill-prepared}} and {\emph{well-prepared}} data) have been
used in the literature. We say that an initial datum is {\emph{well-prepared}} if the compatibility condition  \begin{align}\label{CC1}\big(\rho^\varepsilon u^\varepsilon+\nabla p(\rho^\varepsilon)\big)|_{t=0}=q_0^\varepsilon+\nabla p(\rho_0^\varepsilon)\underset{\varepsilon\to0}{\rightarrow0} \quad \text{in } H^{m-1} \end{align}
holds. In particular, this requires the additional condition $\|q_0^\varepsilon\|_{H^{m-1}}=\mathcal{O}(1)$ compared with the weaker one $\|q_0^\varepsilon\|_{H^{m-1}}=\mathcal{O}(\varepsilon^{-1})$ in the existence result (Proposition \ref{P1.1}).

Our first main result concentrates on the case of {\emph{ill-prepared}} data when such a compatibility condition may not be satisfied. We establish the global-in-time strong convergence of the rescaled Euler equations to the filtration equation and Darcy's law as follows. 


\begin{theorem}
\label{T1.1}
    Let $d\geq1$, $m\geq [\frac{d}{2}]+2$ and \eqref{F1.1}
hold. Let $(\rho^\varepsilon,q^\varepsilon)$ and $(\rho^*,q^*)$ be the solutions obtained in Propositions \ref{P1.1} and \ref{P1.2} subject to the initial data $(\rho_0^\varepsilon,q_0^\varepsilon)$ and $\rho_0^*$, respectively. There exists a constant $C>0$ independent of $\varepsilon$ such that
\begin{equation}
\begin{aligned}
\label{F1.11}
  & \sup_{t\in\mathbb{R}^+}\|(\rho^\varepsilon-\rho^*)(t)\|^2_{H^{m-1}}
+\int_0^{+\infty}\|\nabla(\rho^\varepsilon-\rho^*)(t)\|^2_{H^{m-1}}\,dt\le C(\|\rho^\varepsilon_0-\rho_0^*\|^2_{H^{m-1}}+\varepsilon^2),
\end{aligned}
\end{equation}
and
\begin{equation}
\begin{aligned}
\label{F1.11a}
  & \int_0^{+\infty}\|(q^\varepsilon-q^*-q^\varepsilon_I)(t)\|^2_{H^{m-1}}\,dt\le C(\|\rho^\varepsilon_0-\rho_0^*\|^2_{H^{m-1}}+\varepsilon^2),
\end{aligned}
\end{equation}
where $q^\varepsilon_I$ is an initial time-layer correction given by
\begin{equation}
\label{F1.11b}
   q^\varepsilon_I(t,x) := e^{-\frac{t}{\varepsilon^2}} q^\varepsilon_0(x).
\end{equation}
If in addition $\rho^\varepsilon_0-\rho_0^*\in \dot{H}^{-1}$, then the convergence also holds in the following sense:
\begin{equation}
\begin{aligned}
&\int_0^{\infty} \|(\rho^\varepsilon-\rho^*)(t)\|_{L^2}^2\,dt\leq C(\|\rho^\varepsilon_0-\rho_0^*\|^2_{\dot{H}^{-1}}+\varepsilon^2).\label{H-1}
\end{aligned}
\end{equation}
\end{theorem}

\begin{remark}
The error estimate \eqref{F1.11a} in fact leads to the global strong convergence of $(q^\varepsilon)_{\varepsilon>0}$ to $q^*$. Using \eqref{F1.11a} and the fact that $\|q_I^\varepsilon\|_{L^1(\mathbb{R}^+;H^m)}\leq C\varepsilon(\varepsilon \|q_0^\varepsilon\|_{H^m}) \leq C\varepsilon$, it holds, as $\varepsilon\rightarrow 0$, that
\begin{align*}
q^\varepsilon\rightarrow q^*\;\; \mbox{strongly in}\;\;
L^2(\mathbb{R}^+;H^{m-1})+L^1(\mathbb{R}^+;H^m).
\end{align*}
Here $X+Y$ denotes the sum space of $X$ and $Y$. Furthermore, due to the fact that, as $\varepsilon\rightarrow 0$,
\[  \|q_I^\varepsilon\|_{L^2( [1,\infty);H^m)}\leq C\Big(\int_{1}^{\infty} e^{-\frac{t}{\varepsilon^2}}\,d\tau\Big)^{\frac{1}{2}}\|q_0^\varepsilon\|_{H^m} \leq C\varepsilon,  \]
we have
\begin{align*}
     q^\varepsilon\rightarrow q^*\;\; \mbox{strongly in}\;\;
L^2( [1,\infty);H^{m-1}).
\end{align*}
\end{remark}

\begin{remark}
  The $L^2(\mathbb{R}^+;L^2)$ error estimate \eqref{H-1} is a higher-dimensional version of the estimate that can be found in {\rm\cite{Junca}}. If $\rho_0^\varepsilon=\rho_0^*$, then the $\dot{H}^{-1}$ assumption is not required, and we obtain a $\mathcal{O}(\varepsilon)$-bound for $\rho^\varepsilon-\rho^*$ in $L^2(\mathbb{R}^+;H^m)$.
\end{remark}

\vspace{3mm}
Our next aim is to provide a more precise description of the relaxation approximation in terms of one-order asymptotic expansion. More precisely, we show faster convergence rates between the solutions of \eqref{EMeps} and their first-order approximation. Defining the asymptotic expansion
\begin{align}
\rho^\varepsilon_a:=\rho^*+\varepsilon \rho_1,\quad q^\varepsilon_a:=q^*+\varepsilon q_1, \label{asy:Euler}
\end{align}
the pair $(\rho^\varepsilon_a, q^\varepsilon_a)$ will be used to approximate the solution $(\rho^{\varepsilon}, q^{\varepsilon})$ of \eqref{F1.2} in a suitable sense. Observe that
\begin{equation}
\begin{aligned}
 p(\rho^\varepsilon_a)=p(\rho^*+\varepsilon\rho_1)
=p(\rho^*)+\varepsilon p'(\rho^*)\rho_1+O(\varepsilon^2).\label{pressureO3}
\end{aligned}
\end{equation}
Substituting \eqref{asy:Euler} and \eqref{pressureO3} into \eqref{F1.2} and identifying the coefficients in terms of the power of $\mathcal{O}(\varepsilon)$, we can obtain $(\rho_1,q_1)$ by solving the following linear equations
\begin{equation}\label{rho1q1}
\left\{
\begin{aligned}
& \partial_t\rho_1+\dive q_1=0,\\
&q_1=-\na\big(p'(\rho^*)\rho_1\big),
\end{aligned}
\right.
\end{equation}
with
\begin{equation}
\begin{aligned}
&\rho_1|_{t=0}=\rho_{1,0}. \label{rho1q1d} 
\end{aligned}
\end{equation}
From \eqref{rho1q1}, one sees that $\rho_1$ solves the linear filtration
equation
\begin{equation}\label{rho1}
\begin{aligned}
\partial_t\rho_1-\Delta \big(p'(\rho^*)\rho_1\big)=0.
\end{aligned}
\end{equation}
Our next result proves the strong convergence, at a $\mathcal{O}(\varepsilon
^2)$ rate, of $\rho^\varepsilon$ to its asymptotic expansion $\rho^\varepsilon_a$ for {\emph{well-prepared}} data.

\begin{theorem}
\label{T1.2}
    Let $d\geq1$, $m\geq [\frac{d}{2}]+2$, \eqref{F1.1} and the assumptions in Propositions \ref{P1.1} and \ref{P1.2}
hold. Let $(\rho^\varepsilon,q^\varepsilon)$ and $(\rho^*,q^*)$ be the solutions obtained in Propositions \ref{P1.1} and \ref{P1.2} subject to the initial data $(\rho_0^\varepsilon,q_0^\varepsilon)$ and $\rho_0^*$, respectively. In addition, assume $\rho_{1,0}\in H^{m}$, $q_0^*=\lim\limits_{\varepsilon\rightarrow 0}q_0^\varepsilon=-\nabla p(\rho_0^*)$, 
\begin{align}
\|\rho_0^\varepsilon-\rho_0^*-\varepsilon \rho_{1,0}\|_{H^{m-2}}\leq \varepsilon^2 \quad \text{and}\quad \|q_0^\varepsilon-q_0^*\|_{H^{m-1}}\leq \varepsilon.\label{wellpre}
\end{align}
There exists a constant $C>0$ independent of $\varepsilon$ such that
\begin{equation}
\begin{aligned}
  & \sup_{t\in\mathbb{R}^+}\|(\rho^\varepsilon-\rho^*-\varepsilon \rho_1)(t)\|^2_{H^{m-2}}
+\int_0^{+\infty}\|\nabla(\rho^\varepsilon-\rho^*-\varepsilon \rho_1)(t)\|^2_{H^{m-2}}\,dt\le C\varepsilon^4,
\end{aligned}
\end{equation}
and
\begin{equation}
\left\{
\begin{aligned}
&\int_0^{+\infty}\|(q^\varepsilon-q^*)(t)\|^2_{H^{m-2}}\,dt\le C\varepsilon^2,\\
  & \sup_{t\in\mathbb{R}^+}\|(q^\varepsilon-q^*-\varepsilon q_1)(t)\|^2_{H^{m-2}}\leq C \varepsilon^2,\\
  & \int_0^{+\infty}\|(q^\varepsilon-q^*-\varepsilon q_1)(t)\|^2_{H^{m-2}}\,dt\le C\varepsilon^4.
\end{aligned}
\right.
\end{equation}
\end{theorem}


\subsection{Strategy of proof for the Euler system}
We now provide some comments on the proofs of Theorems \ref{T1.1} and \ref{T1.2}. To justify the strong convergence in the energy space $L^2(\mathbb{R}^+;L^2)$, we adapt the stream function techniques (developed in the one-dimensional settings in \cite{Junca,LiPengZhao1d}) to the multi-dimensional framework. More precisely,  we introduce the {\emph{stream function}} associated with the equation of  $\rho^\varepsilon-\rho^*$:
\begin{equation*}
   N^\varepsilon(t,x):=-\int_0^t(q^\varepsilon-q^*)(t',x)\,dt'-\Lambda^{-2}\nabla (\rho^\varepsilon_0-\rho_0^*)(x),
\end{equation*} 
which satisfies
\begin{equation*}
  \partial_t N^\varepsilon=-(q^\varepsilon-q^*)\quad \text{and}\quad  \dive N^\varepsilon=\rho^\varepsilon - \rho^*.
\end{equation*}
This structure enables us to establish $L^2$ energy estimates for $N^\varepsilon$. Under the additional assumption that $\rho^\varepsilon_0 - \rho_0^* \in \dot{H}^{-1}$, this yields an $L^2(\mathbb{R}^+;L^2)$ estimate for $\rho^\varepsilon - \rho^*$ (see Lemma \ref{LemmaH-1}).

Then, we derive higher-order error estimates, without relying on a $\dot{H}^{-1}$ assumption on the initial data. To this end, we perform energy estimates on the error variable $(\tilde{\rho}^\varepsilon,\tilde{q}^\varepsilon) = (\rho^\varepsilon - \rho^*, q^\varepsilon - q^* - q_I^\varepsilon)$, which satisfies
\begin{equation}\label{sfggg}
\begin{cases}
   \partial_t\tilde{\rho}^\varepsilon+\dive\tilde{q}^\varepsilon=-\dive q_I^\varepsilon,\\[1mm]
\varepsilon^2\partial_t\tilde{q}^\varepsilon+\na\big( p(\rho^\varepsilon)-p(\rho^*)\big)+\tilde{q}^\varepsilon
=-\varepsilon^2\partial_tq^*+\mathcal{R}^\varepsilon,
\end{cases}
\end{equation}
with $\mathcal{R}^\varepsilon=\varepsilon^2\dive\big( q^\varepsilon\otimes q^\varepsilon/\rho^\varepsilon)$. The initial layer correction $q_I^\varepsilon$ ensures that $\tilde{q}^\varepsilon|_{t=0} = -\nabla p(\rho_0^*)$, thus avoiding the singularity at time $t=0$. To control the right-hand side of \eqref{sfggg}, we shall make full use of the following bounds:
\begin{align*}
\|\dive q_I^\varepsilon\|_{L^1(\mathbb{R}^+;H^{m-1})}=\mathcal{O}(\varepsilon),\quad \varepsilon^2\|\partial_t q^*\|_{L^2(\mathbb{R}^+;H^{m-2})}=\mathcal{O}(\varepsilon^2),\quad \|\mathcal{R}^\varepsilon\|_{L^2(\mathbb{R}^+;H^{m-1})}=\mathcal{O}(\varepsilon).
\end{align*}
A key challenge arises due to the limited spatial regularity of $\partial_t q^*$. To overcome this, since the term involving $\partial_t q^*$ is of order $\mathcal{O}(\varepsilon^2)$, we perform estimates in $H^{m-1}$ and use the bounds on $q^{\varepsilon}$ and $q^*$ from the existence results. These considerations lead to an $L^{\infty}(\mathbb{R}^+;H^{m-1})$ $\mathcal{O}(\varepsilon)$-bound for $\tilde{\rho}^\varepsilon$ and an $L^{\infty}(\mathbb{R}^+;H^{m-1}) \cap L^2(\mathbb{R}^+;H^{m-1})$ $\mathcal{O}(\varepsilon)$-bound for $\tilde{q}^\varepsilon$. Moreover, exploiting the dissipative structure induced by the pressure term, we also obtain the higher-order  control of $\nabla \tilde{\rho}^\varepsilon$ in $L^2(\mathbb{R}^+;H^{m-1})$ (see Lemmas \ref{L2.1} and \ref{L2.2}).

Concerning the faster convergence rate $\mathcal{O}(\varepsilon^2)$ for the error terms $(\tilde{\rho}_a^\varepsilon, \tilde{q}_a^\varepsilon) := (\rho^\varepsilon - \rho^\varepsilon_a, q^\varepsilon - q^\varepsilon_a)$, we notice that thanks to the well-prepared assumption \eqref{wellpre}, the right-hand side terms of \eqref{sfggg} are of order $\mathcal{O}(\varepsilon^2)$. Then, our strategy is to reformulate the system as a parabolic equation with $\mathcal{O}(\varepsilon^2)$ source terms and perform refined energy estimates on $\tilde{\rho}_a^\varepsilon$. Then, by reverting to the damped formulation for $\tilde{q}_a^\varepsilon$, we obtain the desired convergence rate for $q^\varepsilon$ toward $q^\varepsilon_a$ (cf. Lemmas \ref{Lrho1q1}–\ref{L2.8}). 


\subsection{The compressible Euler-Maxwell system}

We now present our results for the Euler-Maxwell system \eqref{EMeps}. Let the initial energy $\mathcal{X}_0^\varepsilon$ be defined by
\begin{align}
\mathcal{X}^\varepsilon_0:=\|\rho_{0}^\varepsilon-1\|^2_{H^m}+\varepsilon^2\|u^\varepsilon_{0}\|^2_{H^m}+\|E^\varepsilon_0\|^2_{H^m}+\|B^\varepsilon_0-B^e\|^2_{H^m}.
\end{align}
We define $\phi^*:=\Lambda^{-2}(\rho^*-1)$ and $E^*:=\Lambda^{-2}\nabla \rho^*$. We introduce the momentum variables:
\[   q^\varepsilon:=\rho^\varepsilon u^\varepsilon,\quad q^*:=\rho^*u^*,\quad
q^\varepsilon_0:=\rho^\varepsilon_0 u^\varepsilon_0. \quad 
\]

We recall the uniform global existence result for the rescaled Euler-Maxwell system \eqref{EMeps} from \cite{Wasiolek2016}, which leads to the global weak convergence toward the drift-diffusion system \eqref{F3.3}.

\begin{prop}[\!\!\cite{Wasiolek2016}]
\label{existenceEM} Let $d=3$, $m\geq 3$, and \eqref{F1.1}
hold. There exists a constant $\delta_1>0$ independent of $\varepsilon$ such
that, if
\begin{align}
\mathcal{X}^\varepsilon_0\le\delta_{1},\label{X0eps}
\end{align}
the Cauchy problem \eqref{EMeps}-\eqref{Emepsd}
admits a unique global smooth solution $(\rho^\varepsilon,u^\varepsilon,E^\varepsilon,B^\varepsilon)$ which satisfies $(\rho^\varepsilon-1,q^\varepsilon,E^\varepsilon,B^\varepsilon-B^e)\in C(\mathbb{R}^+;H^m)$ and the uniform estimate
\begin{align}
\label{est:EM}
&\sup_{t\in\mathbb{R}^+}\big(\|\rho^\varepsilon(t)-1\|^2_{H^m}+\varepsilon^2\|u^\varepsilon(t)\|^2_{H^m}+\varepsilon^2\|q^\varepsilon(t)\|^2_{H^m}\big)\nonumber\\
&\quad+\sup_{t\in\mathbb{R}^+}\big(\|E^\varepsilon(t)\|^2_{H^m}+\|B^\varepsilon(t)-B^e\|^2_{H^m}\big)\nonumber\\
&\quad+\int_0^{+\infty}\big(\|\rho^\varepsilon(t)\|^2_{H^{m}}
+\|u^\varepsilon(t)\|^2_{H^m}+\|q^\varepsilon(t)\|^2_{H^m}\big)\,dt\nonumber\\
&\quad+\int_0^{+\infty}\big(\|E^\varepsilon(t)\|^2_{H^{m-1}}+\|\nabla B^\varepsilon(t)\|^2_{H^{m-2}}\big)\,dt
\le C\mathcal{X}^\varepsilon_0,
\end{align}
for a constant $C>0$ independent of $\varepsilon$ and the time.
Moreover, for any finite time $T>0$, as $\varepsilon\to 0$, we have, up to subsequences,
\begin{equation}
\left\{
\begin{aligned}
 \rho^\varepsilon\longrightarrow \rho^*\;\; &\mbox{strongly in}\;\;
C([0,T],H^{m-1}_{loc} ),\\
q^\varepsilon \;-\!\!\!\rightharpoonup q^* \;\;&\mbox{weakly \,\, in} \;\;
L^2(\mathbb{R}^+,H^m ), \\
  E^\varepsilon \longrightarrow E^* \;\;&\mbox{strongly in} \;\;
C([0,T],H^{m-1}_{loc} ),  \\
 B^\varepsilon \longrightarrow B^e \;\;& \mbox{strongly in} \;\;
C([0,T];H^{m-1}_{loc} ),
\end{aligned}
\right.
\end{equation}
where $(\rho^*,u^*,E^*,B^e)$ with $u^*=q^*/\rho^*$ is the unique solution of \eqref{F3.3}-\eqref{DDd}
with the initial data $\rho^*_0$ being the weak limit of $(\rho^\varepsilon_0)_{\varepsilon>0}$ in $H^m(\R)$.
\end{prop}

By standard arguments, we have the following result for the limit system \eqref{F3.3}-\eqref{DDd}.

\begin{prop}
\label{existenceDD}
   Let  $d=3$, $m\geq 3$ and \eqref{F1.1}
hold. There is a positive constant $\delta_1^*$ such
that if $\|\rho^*_0-1\|^2_{H^m} \le\delta_1^*$, then the Cauchy problem \eqref{DD}-\eqref{DDd}
admits a unique global classical solution $(\rho^*,\phi^*)$. Define $E^*=\nabla\phi^*$ and $q^*=-\nabla p(\rho^*)-\rho^* E^*$. Then, $(\rho^*,q^*,E^*)$ satisfies
\begin{align}
\label{est:DD}
   &\sup_{t\in\mathbb{R}^+}\big(\|\rho^*(t)-1\|^2_{H^m}+\|\nabla q^*(t)\|_{H^{m-2}}^2+\|\nabla E^*(t)\|_{H^{m-1}}^2\big)\nonumber\\
   &\quad\quad+\int_0^{+\infty}\big(\|\na\rho^*(t)\|^2_{H^m}+\|\nabla q^*(t)\|_{H^{m-1}}^2+\|\nabla E^*(t)\|_{H^{m}}^2\big)\,dt\le C\|\rho^*_0-1\|^2_{H^m},
\end{align}
where $C>0$ is a generic constant.
\end{prop}

For the relaxation limit of the rescaled Euler-Maxwell system \eqref{EMeps}, the initial data are called {\emph{well-prepared}} if the following compatibility condition holds: \begin{align}\label{CC}\big(q^\varepsilon+\nabla p(\rho^\varepsilon)+\rho^\varepsilon E^\varepsilon\big)|_{t=0}=q_0^\varepsilon+\nabla p(\rho_0^\varepsilon)+\rho_0^\varepsilon E_0^\varepsilon\underset{\varepsilon\to0}{\rightarrow0} \quad \text{in } H^{m-1},\end{align} and {\emph{ill-prepared}} otherwise.


We state our result about the global strong convergence of the relaxation limit from the Euler-Maxwell system to the drift-diffusion model and Darcy's law in the {\emph{ill-prepared}} case.

\begin{theorem}
\label{TEM}
Let the assumptions in Propositions \ref{existenceEM} and \ref{existenceDD} hold, and let $(\rho^\varepsilon, q^\varepsilon,E^\varepsilon,B^\varepsilon)$ and $(\rho^*,q^*,E^*)$ be the solutions obtained in these two statements, respectively. Suppose additionally that $\rho^*_0-1\in \dot{H}^{-1}$.
Then, it holds that
\begin{align}
\label{errorEM}
  & \sup_{t\in\mathbb{R}^+}\Big(\|(\rho^\varepsilon-\rho^*)(t)\|^2_{H^{m-1}}+\|(E^\varepsilon-E^*)(t)\|^2_{H^{m-1}}+\|B^\varepsilon(t)-B^e\|^2_{H^{m-1}}\Big)\nonumber\\
&\quad\quad+\int_0^{+\infty}\Big(\|(\rho^\varepsilon-\rho^*)(t)\|_{H^{m}}^2+\|(E^\varepsilon-E^*)(t)\|^2_{H^{m-1}}+\|\nabla B^\varepsilon(t)\|^2_{H^{m-2}}\Big)\,dt\nonumber\\
&\quad\le C\big(\|\rho^\varepsilon_0-\rho^*_0\|^2_{H^{m-1}}+\|E^\varepsilon_0-E^*_0\|^2_{H^{m-1}}+\|B^\varepsilon_0-B^e\|_{H^{m-1}}^2+\varepsilon^2\big),
\end{align}
and
\begin{align}
\label{errorq}
 & \int_0^{+\infty}\|(q^\varepsilon-q^*-q^\varepsilon_I)(t)\|^2_{H^{m-2}}\,dt\nonumber\\
 &\quad\le C\big(\|\rho^\varepsilon_0-\rho^*_0\|^2_{H^{m-1}}+\|E^\varepsilon_0-E^*_0\|^2_{H^{m-1}}+\|B^\varepsilon_0-B^e\|^2_{H^{m-1}}+\varepsilon^2\big),
\end{align}
where $C$ is a constant independent of $\varepsilon$,
$E^*_0=\nabla\Lambda^{-2}\rho^*_0$,
and $q^\varepsilon_I$ is an initial layer correction given by
\begin{equation}
   q^\varepsilon_I(t,x):=q^\varepsilon_0(x) e^{-\frac{t}{\varepsilon^2}}.\label{qIEM}
\end{equation}
\end{theorem}

\begin{remark}
The condition $\rho^*_0-1\in \dot{H}^{-1}$ ensures the $L^2$-regularity for $E^*$ and $q^*$, which is needed to establish error estimates in $L^{\infty}(\mathbb{R}^+;L^2)$. As observed in {\rm\cite{liEM}} for the periodic case, the electric field error $E^\varepsilon-E^*$ can be viewed as a stream-type function for the Euler-Maxwell system, leading to the $L^2(\mathbb{R}^+;L^2)$-error estimate for $\rho^\varepsilon-\rho^*$. 
\end{remark}


\subsection{Strategy of proof for the Euler-Maxwell system}

The proof of Theorem \ref{TEM} builds upon strategies developed for the compressible Euler equations. When analyzing the error $(\tilde{\rho}^\varepsilon,\tilde{q}^\varepsilon,\tilde{E}^\varepsilon,\tilde{B}^\varepsilon)=(\rho^\varepsilon-\rho^*,q^\varepsilon-q^*,E^\varepsilon-E^*,B^\varepsilon-B^e)$, the main difficulty arises from the presence of the zero-order source term  $\partial_t E_*$ in the equation for $\tilde{E}^\varepsilon$ (see \eqref{eqerrorEM}). To address this, we consider the following first-order asymptotic expansion: 
$$
\big(\rho^\varepsilon_a,q^\varepsilon_a,E^\varepsilon_a,B^\varepsilon_a\big):=\big(\rho^*,q^*,E^*,B^e\big)+\varepsilon\big(\rho_1,q_1,E_1,B_1\big).
$$
We observe that the system for $(\tilde{\rho}_a^\varepsilon,\tilde{q}_a^\varepsilon,\tilde{E}_a^\varepsilon,\tilde{B}_a^\varepsilon)=(\rho^\varepsilon-\rho^\varepsilon_a,q^\varepsilon-q^\varepsilon_a,E^\varepsilon-E^\varepsilon_a,B^\varepsilon-B^\varepsilon_a)$ shares a similar dissipative structure to that of the linearized Euler–Maxwell system, modulo source terms that are uniformly bounded by $\mathcal{O}(\varepsilon)$. This enables us to show $\mathcal{O}(\varepsilon)$-bounds for $(\tilde{\rho}_a^\varepsilon,\tilde{q}_a^\varepsilon,\tilde{E}_a^\varepsilon,\tilde{B}_a^\varepsilon)$. Then, since the profile $\big(\rho_1,q_1,E_1,B_1\big)$ is globally defined, one can recover the desired $\mathcal{O}(\varepsilon)$-bounds for the original error unknown $(\tilde{\rho}^\varepsilon,\tilde{q}^\varepsilon,\tilde{E}^\varepsilon,\tilde{B}^\varepsilon)$ by combining refined error estimates with the asymptotic expansion.

\section{Preliminaries and toolbox}\label{section3}

Throughout this paper, $C>0$ denotes a harmless constant independent of $t$ and $\varepsilon$. The notation $\mathcal{F}(f)$ and  $\mathcal{F}^{-1}(f)$ stand for the Fourier transform and the inverse Fourier transform of the function $f$, and we write
$$
\Lambda^{\sigma}:= (-\Delta)^\frac{\sigma}{2}=\mathcal{F}^{-1} \Big( |\xi|^{\sigma} \mathcal{F}(\cdot ) \Big),\quad\sigma\in\mathbb{R}.
$$
In the case $\sigma=2$, we have $-\Delta=\Lambda^{2}$.
Let $s\in\mathbb{R}$, we denote by $H^s$ the Sobolev space
of exponent $s$ with the standard norm
$$
\|f\|_{H^s}:=\Big(\int_{\mathbb{R}^d} (1+|\xi|^2)^s|\mathcal{F}(f) (\xi)|^2 \,d\xi \Big)^{\frac{1}{2}}<\infty.
$$
In particular, when $s$ is a positive integer, we have
$$
\|f\|_{H^s}\sim \sum_{0\leq |\alpha|\leq s}\|\partial^{\alpha} f\|_{L^2}^2.
$$
Furthermore, we denote by $\dot{H}^s(\mathbb{R}^d)$
the homogeneous Sobolev space endowed with the norm
$$
\|f\|_{\dot{H}^s}:=\Big(\int_{\mathbb{R}^d} |\xi|^{2s}|\mathcal{F}(f) (\xi)|^2 \,d\xi \Big)^{\frac{1}{2}}= \|\Lambda^{s}f\|_{L^2}.
$$

The following Sobolev and Moser-type inequalities can be found in
\cite{Majdalocal,Katolocal}. They will be repeatedly used in the proof.
\begin{lemma}[\!\!\cite{Majdalocal,Katolocal}]
The following estimates hold.
\vspace{1mm}
\begin{itemize}
\item Let $s>\frac{d}{2}$. The embedding from $H^{s}$ to
$L^\infty\cap C^0$ is continuous and for $u\in H^{s}$,
\begin{equation}
\label{F1.12}
   \|u\|_{L^\infty}\le C\|u\|_{H^{s}}.
\end{equation}

\item Let $-\frac{d}{2}<s<0$ and $q=\frac{2d}{d-2s}$. The embedding from $L^q$ to $H^{s}$ is continuous and for $u\in L^q$,
\begin{equation}
\label{F1.120}
   \|u\|_{H^s}\le C\|u\|_{L^q}.
\end{equation}

\item Let $s_1,s_2\in \mathbb{R}$. The operator $\Lambda^{s_1}$ is an isomorphism from $H^{s_2}$ to $H^{s_2-s_1}$.

\item For
$u\in H^{m-1}$ and $v\in H^{s}$ with $s>0$,
\begin{equation}
   \|uv\|_{\dot{H}^{s}}\le C \|u\|_{L^{\infty}}\|v\|_{\dot{H}^s}+C\|u\|_{\dot{H}^s}\|v\|_{L^{\infty}}.\label{F1.1311}
\end{equation}
Consequently,  the Sobolev space $H^{s}$ with $s>\frac{d}{2}$ is an algebra and we have
\begin{equation}
\label{F1.13}
  \|uv\|_{H^{s}}\le C \|u\|_{H^{s}}\|v\|_{H^s}.
\end{equation}

\item For $u\in H^{s_1}$ and $v\in \dot{H}^{s}$ with $s_1>\frac{d}{2}$ and $-\frac{d}{2}<s_2<\frac{d}{2}$, we have $uv\in \dot{H}^{s}$ and
\begin{align}
\|uv\|_{\dot{H}^{s}}\leq C\|u\|_{H^{s_1}}\|v\|_{\dot{H}^{s}}.\label{product}
\end{align}


\item Let $g$ be a smooth function on a compact set
$D\subset\mathbb{R}^n$ and let $k\ge 1$ be an integer. For $u\in H^k\cap L^\infty$
satisfying $u(x)\in D$,
\begin{equation}
\label{F1.15}
   \sum_{|\alpha|=k} \|\partial^{\alpha} g(u)\|_{L^2}\le C \|g\|_{C^k(D)}\|u\|^{k-1}_{L^\infty}\|u\|_{\dot{H}^k}.
\end{equation}

\end{itemize}

\end{lemma}
We recall the Gagliardo-Nirenberg inequality.
\begin{lemma}
\label{Lt9}
    Assume that $q,r$ satisfy $1\leq q,r\le \infty$ and $j,m\in \mathbb{Z}^+$
satisfy $0\le j<m$. For any $f\in C_0^{\infty}(\mathbb{R}^d)$, we have
	\begin{equation}
	\|D^jf\|_{L^p}\leq
C\|D^mf\|_{L^r}^{\alpha}\|f\|_{L^q}^{1-\alpha}\label{w11},
	\end{equation}
where
\[  \frac{1}{p}-\frac{j}{d}
=\alpha(\frac{1}{r}-\frac{m}{d})+(1-\alpha)\frac{1}{q},\quad
\frac{j}{m}\le\alpha\le1,  \]
and $C$ depends on $m,d,j,q,r,\alpha$.
\end{lemma}

\section{The Euler system (I): Proof of Theorem \ref{T1.1}}

\subsection{The error equations}
$\,$
\vspace{1mm}

In what follows, we will always denote by $(\rho^\varepsilon, q^\varepsilon)$ and $(\rho^*, q^*)$ the solutions to the Euler equations and the filtration equation obtained in Propositions \ref{P1.1} and  \ref{P1.2}, respectively.

We introduce the error unknowns
\begin{align}
\tilde{\rho}^\varepsilon:=\rho^\varepsilon-\rho^*,\quad\quad \tilde{q}^\varepsilon:=q^\varepsilon-q^*-q^\varepsilon_I,
\end{align}
where the initial layer correction $q^\varepsilon_I$ given by \eqref{F1.11b} allows
us to avoid the time singularity of $\tilde{q}^\varepsilon$ at $t=0$.
Note that $q^\varepsilon_I$ satisfies
\begin{align}
\varepsilon^2\partial_t q^\varepsilon_I+q^\varepsilon_I=0,\quad\quad q^\varepsilon_I|_{t=0}=q^\varepsilon_0=\frac{1}{\varepsilon}\rho_0^\varepsilon u_0^\varepsilon.
\end{align}
With these notations, the error system reads
\begin{equation}\label{eqwiderhoqa}
\begin{cases}
   \partial_t\tilde{\rho}^\varepsilon+\dive\tilde{q}^\varepsilon=\mathcal{R}^\varepsilon_1,\\[1mm]
\varepsilon^2\partial_t\tilde{q}^\varepsilon+p'(1)\na\tilde{\rho}^\varepsilon+\tilde{q}^\varepsilon
=\mathcal{R}^\varepsilon_0
+\mathcal{R}^\varepsilon_2+\mathcal{R}^\varepsilon_3,\\
(\tilde{\rho}^\varepsilon, \tilde{q}^\varepsilon)|_{t=0}=(\rho^\varepsilon-\rho_0^*, \nabla p(\rho_0^*)),
\end{cases}
\end{equation}
where
\begin{equation}
\left\{
\begin{aligned}
\mathcal{R}^\varepsilon_0&=-\big(p'(\rho^\varepsilon)-p'(\rho^*)\big)\na\rho^*-\big( p'(\rho^\varepsilon)-p'(1)\big)\na\tilde{\rho}^\varepsilon,\\
\mathcal{R}^\varepsilon_1&=-\dive q^\varepsilon_I,\\
\mathcal{R}^\varepsilon_2&=-\varepsilon^2 \dive \Big( \frac{q^\varepsilon\otimes q^\varepsilon}{\rho^\varepsilon}\Big),\\
\mathcal{R}^\varepsilon_3&=-\varepsilon^2\partial_t q^*.\label{R1R2R30}
\end{aligned}
\right.
\end{equation}

The term $\mathcal{R}^\varepsilon_0$ will be controlled using the bounds for the error $\tilde{\rho}^\varepsilon$ and the smallness for $\rho^\varepsilon-1$ and $\rho^*-1$. Concerning the remainder terms $\mathcal{R}^\varepsilon_i$ ($i=1,2,3$), we will use the following lemma.
\begin{lemma}
\label{lemmareminderEuler11}
Let $\mathcal{R}_i$ {\rm(}$i=1,2,3${\rm)} be given by \eqref{R1R2R30}.
We have
\begin{align}
&\int_{0}^{\infty}\|\mathcal{R}^\varepsilon_1(t)\|_{H^{m-1}}^2\,dt\leq C,\quad\quad \int_{0}^{\infty}\|\mathcal{R}^\varepsilon_1(t)\|_{H^{m-1} }\,dt\leq C\varepsilon, \label{errorr1}\\
&\int_{0}^{\infty}\|\mathcal{R}^\varepsilon_2(t)\|_{H^{m-1} }^2\,dt\leq C\varepsilon^2, \label{errorr2}\\
&\int_{0}^{\infty}\|\mathcal{R}^\varepsilon_3(t)\|_{H^{m-2} }^2\,dt\leq C\varepsilon^4,\label{errorr3}
\end{align}
where $C>0$ is a constant independent of $\varepsilon$.
\end{lemma}

\noindent {\bf Proof.}
For the term $\mathcal{R}^\varepsilon_1$, using the product law \eqref{F1.13}, we have
\begin{equation}\nonumber
\begin{aligned}
\|q^\varepsilon_0\|_{H^m}&\leq \frac{1}{\varepsilon}\|u^\varepsilon_0\|_{H^m}+\frac{1}{\varepsilon}\|(\rho_0^\varepsilon-1)u^\varepsilon_0\|_{H^m}\leq \frac{C}{\varepsilon}(1+\|\rho_0^\varepsilon-1\|_{H^m})\|u^\varepsilon_0\|_{H^m}.
\end{aligned}
\end{equation}
Together with the definition $\mathcal{R}_1^\varepsilon=-\dive q_{I}^\varepsilon=-e^{-\frac{t}{\varepsilon^2}} \dive  q^\varepsilon_0 $ and the assumptions of $\rho_0^\varepsilon$, $u_0^\varepsilon$ and $\rho_0^*$, this implies that
\begin{equation}\nonumber
\begin{aligned}
\int_0^{\infty}\|\mathcal{R}_1^\varepsilon(t)\|_{H^{m-1}}^2\,dt\leq C\int_0^\infty e^{-\frac{2t}{\varepsilon^2}} \,dt \,\|\dive q_0^\varepsilon\|_{H^{m-1}}^2 \leq C\varepsilon^2\|q_0^\varepsilon\|_{H^{m}}^2\leq C.
\end{aligned}
\end{equation}
Similarly, we have
\begin{equation}\nonumber
\begin{aligned}
\int_0^{\infty}\|\mathcal{R}_1^\varepsilon(t)\|_{H^{m-1}}\,dt=\int_0^\infty e^{-\frac{t}{\varepsilon^2}} \,dt \,\|\dive q_0^\varepsilon\|_{H^{m-1}} \leq C \varepsilon^2 \|q_0^\varepsilon\|_{H^{m}}\leq C\varepsilon.
\end{aligned}
\end{equation}
Thus, \eqref{errorr1} follows.

Concerning $\mathcal{R}_2^\varepsilon$, it holds by \eqref{F1.8}, \eqref{F1.13} and \eqref{F1.15} that
\begin{equation*}
\begin{aligned}
&\quad \int_0^{\infty}\|\mathcal{R}^\varepsilon_2(t)\|_{H^{m-1}}^2\,dt\\
&\leq  \quad C\varepsilon^4 \int_0^{\infty} \Big(1+\Big\|\frac{1}{\rho^\varepsilon}(t)-1\Big\|_{H^m}^2 \Big)\|q^\varepsilon(t)\|_{H^m}^4\, dt\\
&\leq \quad C\varepsilon^2\big(1+\sup_{t\in\mathbb{R}^+}\|\rho^\varepsilon(t)-1\|_{H^m}^2\big)\sup_{t\in\mathbb{R}^+}\varepsilon^2\|q^\varepsilon(t)\|_{H^m}^2\int_0^{\infty}\|q^\varepsilon(t)\|_{H^m}^2\,dt \leq C\varepsilon^2.
\end{aligned}
\end{equation*}

Finally, for $\mathcal{R}_{3}^\varepsilon$, noting that 
\begin{equation}
\begin{aligned}
&\partial_{t}q^*=\nabla (p'(\rho^*)\partial_{t}\rho^*)=-\nabla (p'(\rho^*)\dive q^*),\nonumber
\end{aligned}
\end{equation}
we have
\begin{equation}\label{2151}
\begin{aligned}
\int_0^{\infty}\|\partial_{t}q^*(t)\|_{H^{m-2}}^2\,dt&\leq  \big(1+\sup_{t\in\mathbb{R}^+}\|\rho^*(t)-1\|_{H^{m}}^2\big)\int_0^\infty \|q^*(t)\|_{H^{m}}^2\,dt\leq C,
\end{aligned}
\end{equation}
where we used \eqref{F1.9} and \eqref{F1.13}. This completes
the proof of Lemma \ref{lemmareminderEuler11}. \hfill $\Box$

\subsection{Error estimates in $L^2( \mathbb{R}^+;L^2)$}

\vspace{1mm}
We are now in a position to establish estimates for $(\tilde{\rho}^{\varepsilon},\tilde{q}^{\varepsilon})=(\rho^\varepsilon-\rho^*,q^\varepsilon-q^*-q_I^\varepsilon)$. First, we obtain a rate of convergence for $\tilde{\rho}^\varepsilon$ in $L^2(\mathbb{R}^{+};L^2)$.

\begin{lemma}\label{LemmaH-1}
Assuming $\rho_0^\varepsilon-\rho_0^*\in \dot{H}^{-1}$, we have
\begin{align}\label{H-1:}
&\int_0^\infty \|\tilde{\rho}^\varepsilon(t)\|_{L^2}^2\,dt\le C\|\rho_0^\varepsilon-\rho_0^*\|_{\dot{H}^{-1}}^2+C\varepsilon^2.
\end{align}
\end{lemma}
\noindent {\bf Proof.}
We define the {\emph{stream function}}
\begin{equation}
\label{S}
   N^\varepsilon(t,x):=-\int_0^t(q^\varepsilon-q^*)(t',x)\,dt'+N_0^\varepsilon(x),
\end{equation}
with $N_0^\varepsilon(x):=-\Lambda^{-2}\nabla (\rho^\varepsilon_0-\rho_0^*)(x)$ such that
\begin{equation}
\label{F2.3N}
  \partial_t N^\varepsilon=-(q^\varepsilon-q^*),\qquad \dive N^\varepsilon=\tilde{\rho}^\varepsilon.
\end{equation}
Note that \eqref{eqwiderhoqa} can be rewritten as
\begin{equation}
\label{F2.2N}
\begin{cases}
   \partial_t\tilde{\rho}^\varepsilon+\dive (q^\varepsilon-q^*)=0,\\[1mm]
q^\varepsilon-q^*=-\na\big(p(\rho^\varepsilon)-p(\rho^*)\big)-\varepsilon^2 \partial_t q^\varepsilon-\mathcal{R}_2^\varepsilon.
\end{cases}
\end{equation}
Now we establish the desired estimate \eqref{H-1:} using \eqref{F2.3N} and \eqref{F2.2N}. Taking the $L^2$-inner product of $\eqref{F2.2N}_1$ by $N^\varepsilon$ and using \eqref{F2.3N} and $\eqref{F2.2N}_2$, we have
\begin{align}
   \frac{1}{2}\frac{d}{dt}\|N^\varepsilon\|_{L^2}^2&=-\langle q^\varepsilon-q^*,N^\varepsilon\rangle\nonumber\\
&=\langle\na\big(p(\rho^\varepsilon)-p(\rho^*)\big),N^\varepsilon\rangle
+\varepsilon^2\langle\partial_t q^\varepsilon,N^\varepsilon\rangle
+\langle \mathcal{R}_2^\varepsilon,N^\varepsilon\rangle,  \label{ddtN}
\end{align}
where $\langle\cdot,\cdot\rangle$ stands for the inner product of
$L^2$.

In what follows, we estimate each term on the right-hand side of \eqref{ddtN}. First, since the pressure satisfies   \eqref{F1.1} and $\rho^\varepsilon$ is close to $1$ uniformly, there exists a uniform constant $p_1>0$ such that
\begin{equation}
\label{F2.6}
  \langle\na\big(p(\rho^\varepsilon)-p(\rho^*)\big),N^\varepsilon\rangle
=-\langle p(\rho^\varepsilon)-p(\rho^*),\rho^\varepsilon-\rho^*\rangle
\le -p_1\|\tilde{\rho}^\varepsilon\|_{L^2}^2.
\end{equation}
Next, according to \eqref{F2.3N}, it holds that
\begin{align*}
   \varepsilon^2\langle\partial_t q^\varepsilon,N^\varepsilon\rangle
&=\varepsilon^2\frac{d}{dt}\langle q^\varepsilon,N^\varepsilon\rangle
-\varepsilon^2\langle q^\varepsilon,\partial_t N^\varepsilon\rangle\\
&=\varepsilon^2\frac{d}{dt}\langle q^\varepsilon,N^\varepsilon\rangle
+\varepsilon^2\langle q^\varepsilon,q^\varepsilon-q^*\rangle.
\end{align*}
As $q^*=-\na p(\rho^*)$, the second term on the right-hand side of the above equality can be analyzed by
\[
\varepsilon^2\langle q^\varepsilon,q^\varepsilon-q^*\rangle
\le C\varepsilon^2\big(\|q^\varepsilon\|_{L^2}^2+\|\na\rho^*\|_{L^2}^2\big).
\]
Consequently, we have
\begin{equation}
\label{F2.7}
   \varepsilon^2\langle\partial_t q^\varepsilon,N^\varepsilon\rangle
\le 2\varepsilon^2\frac{d}{dt}\langle q^\varepsilon,N^\varepsilon\rangle
+C\varepsilon^2\big(\|q^\varepsilon\|_{L^2}^2+\|\na\rho^*\|_{L^2}^2\big).
\end{equation}
Moreover, from \eqref{F1.8} and \eqref{R1R2R30}, it is clear that
\begin{align}
  \langle \mathcal{R}_2^\varepsilon,N^\varepsilon\rangle& \leq \|\mathcal{R}_2^\varepsilon\|_{L^2}\|N^\varepsilon\|_{L^2} \nonumber\\
  &\leq C\varepsilon^2 \|q^\varepsilon\|_{H^m}^2\|N^\varepsilon\|_{L^2}\nonumber\\
  &\leq C\varepsilon^4 \|q^\varepsilon\|_{H^m}^2+C \|q^\varepsilon\|_{H^m}^2\|N^\varepsilon\|_{L^2}^2.\label{F2.8N}
\end{align}
Since $\varepsilon\in (0,1]$, combining \eqref{ddtN}-\eqref{F2.8N} together with
\eqref{F1.8}, we obtain
\begin{align}
\label{F2.5}
 &\frac{d}{dt}\big(\|N^\varepsilon(t)\|_{L^2}^2-2\varepsilon^2\langle q^\varepsilon,N^\varepsilon\rangle\big)
+2p_1\|\tilde{\rho}^\varepsilon\|_{L^2}^2\\
&\le C\varepsilon^2\big(\|q^\varepsilon\|_{H^m}^2+\|\na\rho^*\|_{H^m}^2\big)+C\|q^\varepsilon\|_{H^m}^2 \|N^\varepsilon\|_{L^2}^2,\nonumber
\end{align}
which implies that, for all $t\geq0$,
\begin{align}
\label{F2.5N}
&\quad \|N^\varepsilon(t)\|_{L^2}^2+2p_1\int_0^{t}\|\tilde{\rho}^\varepsilon(t')\|_{L^2}^2\,dt'\nonumber\\
& \le \|N^\varepsilon_0\|_{L^2}^2+ 2\varepsilon^2\langle q^\varepsilon,N^\varepsilon\rangle\Big|^t_0+C\varepsilon^2\int_0^t\big(\|q^\varepsilon(t')\|_{H^m}^2+\|\na\rho^*(t')\|_{H^m}^2\big)\,dt'\nonumber\\
&\quad+C\int_0^t\|q^\varepsilon(t')\|_{H^m}^2 \|N^\varepsilon(t')\|_{L^2}^2\,dt'.
\end{align}
Here, one has
$$
2\varepsilon^2\langle q^\varepsilon,N^\varepsilon\rangle\Big|^t_0\leq \frac{1}{2}\|N^\varepsilon\|_{L^2}^2+C\varepsilon^4 \|q^\varepsilon\|_{L^2}^2+C\|N_0^\varepsilon\|_{L^2}^2+C\varepsilon^4 \|q_0^\varepsilon\|_{L^2}^2.
$$
Therefore, applying Gr\"onwall's lemma  to \eqref{F2.5N}, \eqref{F1.8} and the fact that $\varepsilon^2\|q_0^\varepsilon\|_{L^2}^2\leq C$ leads to
\begin{align}
\sup_{t\in\mathbb{R}^+}\|N^\varepsilon(t)\|_{L^2}^2+\int_0^{\infty}\|\tilde{\rho}^\varepsilon(t)\|_{L^2}^2\,dt\leq C\|N_0^\varepsilon\|_{L^2}^2+C\varepsilon^2.\nonumber
\end{align}
Since $-\Lambda^{-2}\nabla$ maps $\dot{H}^{-1}$ to $L^2$, we have $\|N_0^\varepsilon\|_{L^2}\leq C\|\rho_0^\varepsilon-\rho_0^*\|_{\dot{H}^{-1}}$. Therefore, we arrive at \eqref{H-1:}. \qed

\subsection{Error estimates in $L^{\infty}(\mathbb{R}^+;L^2)\cap L^2(\mathbb{R}^+;\dot{H}^1)$}

Without assuming that the initial error belongs to $\dot{H}^{-1}$, we have the following lemma:

\begin{lemma}
\label{L2.1}
    It holds
\begin{align}\label{2.31}
&\sup_{t\in\mathbb{R}^+}\big(\|\tilde{\rho}^\varepsilon(t)\|_{L^2}^2+\varepsilon^2\|\tilde{q}^\varepsilon(t)\|_{L^2}^2\big)+\int_0^\infty \big(\|\nabla\tilde{\rho}^\varepsilon(t)\|_{L^2}^2+\|\tilde{q}^\varepsilon(t)\|_{L^2}^2 \big)\,dt\nonumber\\
&\quad \le C\|\rho_0^\varepsilon-\rho_0^*\|_{L^2}^2+C\varepsilon^2,
\end{align}
where $C>0$ is a constant independent of $\varepsilon$.
\end{lemma}

\noindent {\bf Proof.}
Taking the $L^2$-inner product of $\eqref{eqwiderhoqa}_1$ with $p'(1)\tilde{\rho}^\varepsilon$ and of $\eqref{eqwiderhoqa}_2$ with $\tilde{q}^\varepsilon$, respectively, we arrive at
\begin{align}\label{F3.14a}
&\frac{1}{2}\frac{d}{dt}\Big( p'(1)\|\tilde{\rho}^\varepsilon\|_{L^2}^2+\varepsilon^2\|\tilde{q}^\varepsilon\|
^{2}_{L^2}\Big)+\|\tilde{q}^\varepsilon\|_{L^2}^2\nonumber\\
=&p'(1)\langle \mathcal{R}_1^\varepsilon, \tilde{\rho}^\varepsilon\rangle+\langle \mathcal{R}^\varepsilon_0+\mathcal{R}_2^\varepsilon+\mathcal{R}_3^\varepsilon, \tilde{q}^\varepsilon\rangle\\
\leq& p'(1)\|\mathcal{R}_1^\varepsilon\|_{L^2}\|\tilde{\rho}^\varepsilon\|_{L^2}+C\big(\|\mathcal{R}^\varepsilon_0\|_{L^2}^2+\|\mathcal{R}_2^\varepsilon\|_{L^2}^2+\|\mathcal{R}_3^\varepsilon\|_{L^2}^2\big)+\frac{1}{2}\|\tilde{q}^\varepsilon\|_{L^2}^2. \nonumber
\end{align}
To bound $\mathcal{R}^\varepsilon_0$, we use the fact that
$$
\mathcal{R}^\varepsilon_0=-\pi(\rho^\varepsilon, \rho^*)\tilde{\rho}^\varepsilon \nabla \rho^*- \pi(\rho^\varepsilon, 1)(\rho^\varepsilon-1)\nabla \tilde{\rho}^\varepsilon.
$$
where
$$
   p'(s_1)-p'(s_2)=\pi(s_1, s_2) (s_1-s_2),\quad\quad \pi(s_1, s_2):=\int_0^1 p''(\theta s_1+(1-\theta)s_2) d\theta.
$$
Due to \eqref{F1.8} and \eqref{F1.9}, the densities $\rho^\varepsilon$ and $\rho^*$ have uniform upper bounds, so $\pi(\rho^\varepsilon, \rho^*)$ and $\pi(\rho^\varepsilon, 1)$ are also bounded from above. In the case $d\geq 3$, the Cauchy-Schwarz inequality, together with Sobolev's embeddings $\dot{H}^1\hookrightarrow L^{\frac{2d}{d-2}}$ and $H^{m-1}\hookrightarrow L^{2}\cap L^{\infty}\hookrightarrow L^{d}\cap L^{\infty}$, leads to
\begin{equation*}
\begin{aligned}
 \| \mathcal{R}^\varepsilon_0\|_{L^2}^2
&\le \|\pi(\rho^\varepsilon, \rho^*)\|_{L^{\infty}}^2\|\tilde{\rho}^\varepsilon\|_{L^{\frac{2d}{d-2}}}^2 \|\nabla \rho^*\|_{L^{d}}^2+\|\pi(\rho^\varepsilon,1)\|_{L^{\infty}}^2 \|\rho^\varepsilon-1\|_{L^{\infty}}^2\|\nabla\tilde{\rho}^\varepsilon\|_{L^{2}}^2 \\
&\le C\big(\| \rho^\varepsilon-1\|_{H^{m-1}}^2+\| \rho^*-1\|_{H^{m-1}}^2\big)\|\nabla \tilde{\rho}^\varepsilon\|_{L^2}^2\\
&\le C(\delta+\delta^*)\|\nabla \tilde{\rho}^\varepsilon\|_{L^2}^2 .
\end{aligned}
\end{equation*}
In the case $d=2$, we take advantage of the Gagliardo-Nirenberg-Sobolev inequality in Lemma \ref{Lt9} and \eqref{F1.9} to derive
\begin{equation*}
\begin{aligned}
 \| \mathcal{R}^\varepsilon_0\|_{L^2}^2 &\le \|\pi(\rho^\varepsilon, \rho^*)\|_{L^{\infty}}^2\|\tilde{\rho}^\varepsilon\|_{L^{4}}^2 \|\nabla \rho^*\|_{L^4}^2+\|\pi(\rho^\varepsilon, 1)\|_{L^{\infty}}^2\|\rho^\varepsilon-1\|_{L^{\infty}}^2 \|\nabla \tilde{\rho}^\varepsilon\|_{L^2}^2\\
&\le C\|\tilde{\rho}^\varepsilon\|_{L^{2}}\|\nabla\tilde{\rho}^\varepsilon\|_{L^{2}} \|\nabla \rho^*\|_{L^2}\|\nabla^2 \rho^*\|_{L^2}+C\|\rho^\varepsilon-1\|_{H^{m-1}}^2 \|\nabla \tilde{\rho}^\varepsilon\|_{L^2}^2\\
&\leq C(\delta+\delta^*)\|\nabla\tilde{\rho}^\varepsilon\|_{L^{2}}^2+C\|\nabla^2 \rho^*\|_{L^2}^2\|\tilde{\rho}^\varepsilon\|_{L^2}^2.
\end{aligned}
\end{equation*}
Finally, the case $d=1$ can be addressed similarly:
\begin{equation*}
\begin{aligned}
 \| \mathcal{R}^\varepsilon_0\|_{L^2}^2 &\le \|\pi(\rho^\varepsilon,\rho^*)\|_{L^{\infty}}^2\|\tilde{\rho}^\varepsilon\|_{L^{\infty}}^2 \|\nabla \rho^*\|_{L^2}^2+\|\pi(\rho^\varepsilon, 1)\|_{L^{\infty}}^2\|\rho^\varepsilon-1\|_{L^{\infty}}^2 \|\nabla \tilde{\rho}^\varepsilon\|_{L^2}^2\\
&\le C\|\tilde{\rho}^\varepsilon\|_{L^{2}} \|\nabla\tilde{\rho}^\varepsilon\|_{L^{2}} \|\nabla \rho^*\|_{L^2}^2+C \|\rho^\varepsilon-1\|_{H^{m-1}}^2 \|\nabla \tilde{\rho}^\varepsilon\|_{L^2}^2\\
&\leq C(\delta+\delta^*)\|\nabla\tilde{\rho}^\varepsilon\|_{L^{2}}^2+C\|\nabla \rho^*\|_{L^2}^2\|\tilde{\rho}^\varepsilon\|_{L^2}^2.
\end{aligned}
\end{equation*}
Consequently, for all $d\geq1$, we have
\begin{align}\label{caseall}
\int_0^t \| \mathcal{R}^\varepsilon_0(t')\|_{L^2}^2\,dt'&\leq C(\delta+\delta^*)\int_0^t \|\nabla\tilde{\rho}^\varepsilon(t')\|_{L^{2}}^2\,dt'\nonumber\\
&\quad+C\sup_{t'\in[0,t]} \|\tilde{\rho}^\varepsilon(t')\|_{L^2}^2 \int_0^t \|\nabla \rho^*(t')\|_{H^1}^2\,dt'.
\end{align}
Noting that $\varepsilon^2\| \tilde{q}^\varepsilon|_{t=0}\|_{L^2}^2=\varepsilon^2\|\nabla p(\rho^*_0)\|_{L^2}\leq C\varepsilon^2$, and integrating \eqref{F3.14a} over $[0,t']$ with any $0\leq t'\leq t<\infty$, using \eqref{caseall} and then taking the supremum over $[0,t]$, we have
\begin{equation*}
\begin{aligned}
&\quad \sup_{t'\in[0,t]}\big(\|\tilde{\rho}^\varepsilon(t')\|_{L^2}^2
+\varepsilon^2\| \tilde{q}^\varepsilon(t')\|_{L^2}^2\big)
+\int_0^t\|\tilde{q}^\varepsilon(t')\|_{L^2}^2\,dt'\\
&\le C\|\rho_0^\varepsilon-\rho_0^*\|_{L^2}^2+C\varepsilon^2
+C(\delta+\delta^*)\int_0^t\|\nabla \tilde{\rho}^\varepsilon(t')\|_{L^2}^2\,dt'\\
&\quad+C\sup_{t'\in[0,t]} \|\tilde{\rho}^\varepsilon(t')\|_{L^2}^2\int_0^t\|\nabla \rho^*(t')\|_{H^1}^2\,dt'\\
&\quad+C\sup_{t'\in[0,t]}\|\tilde{\rho}^\varepsilon(t')\|_{L^2}\int_0^t\|\mathcal{R}_1^\varepsilon(t')\|_{L^2}\,dt'+C\int_0^t\big(\|\mathcal{R}_2^\varepsilon(t')\|^2_{L^2}+\|\mathcal{R}_3^\varepsilon(t')\|^2_{L^2}\big)\,dt'\\
&\leq C\|\rho_0^\varepsilon-\rho_0^*\|_{L^2}^2+C\varepsilon^2+\Big(\frac{1}{2}+C\delta^*\Big)\sup_{t'\in[0,t]}\|\tilde{\rho}^\varepsilon(t')\|_{L^2}^2+C(\delta+\delta^*) \int_{0}^{t}\|\nabla \tilde{\rho}^\varepsilon(t')\|_{L^2}^2\,dt',
\end{aligned}
\end{equation*}
where we have used \eqref{errorr1}-\eqref{errorr3} and the fact that
$$
C\sup_{t'\in[0,t]}\|\tilde{\rho}^\varepsilon(t')\|_{L^2}\int_0^t\|  \mathcal{R}_1^\varepsilon(t')\|_{L^2}\,dt'\leq \frac{1}{2} \sup_{t'\in[0,t]}\|\tilde{\rho}^\varepsilon(t')\|_{L^2}^2+ C\Big(\int_0^t\|  \mathcal{R}_1^\varepsilon(t')\|_{L^2}\,dt'\Big)^2.
$$
Since $\delta^*$ can be chosen sufficiently small, we conclude that
\begin{align}\label{tilderhoL2}
&\sup_{t'\in[0,t]}\big(\|\tilde{\rho}^\varepsilon(t')\|_{L^2}^2
+\varepsilon^2\| \tilde{q}^\varepsilon(t')\|_{L^2}^2\big)
+\int_0^t\|\tilde{q}^\varepsilon(t')\|_{L^2}^2\,dt'\nonumber\\
&\quad\leq C\|\rho_0^\varepsilon-\rho_0^*\|_{L^2}^2+C\varepsilon^2 +C(\delta+\delta^*) \int_{0}^{t}\|\nabla \tilde{\rho}^\varepsilon(t')\|_{L^2}^2\,dt'.
\end{align}

In order to address the last term in \eqref{tilderhoL2}, we shall establish the $L^2(0,t;L^2 )$ estimate of $\nabla\tilde{\rho}^\varepsilon$. Taking the inner product of $\eqref{eqwiderhoqa}_2$ with $\nabla\tilde{\rho}^\varepsilon$ leads to
\begin{align}\label{L2cross1}
p'(1) \|\nabla\tilde{\rho}^\varepsilon\|_{L^2}^2&=-\varepsilon^2\langle\partial_{t} \tilde{q}^\varepsilon, \nabla \tilde{\rho}^\varepsilon\rangle -\langle \tilde{q}^\varepsilon, \nabla\tilde{\rho}^\varepsilon\rangle+\langle \mathcal{R}^\varepsilon_0+
\mathcal{R}^\varepsilon_2+\mathcal{R}^\varepsilon_3, \nabla\tilde{\rho}^\varepsilon\rangle \nonumber\\
&\leq -\varepsilon^2\langle\partial_{t} \tilde{q}^\varepsilon, \nabla \tilde{\rho}^\varepsilon\rangle+\frac{1}{2}p'(1)\|\nabla\tilde{\rho}^\varepsilon\|_{L^2}^2+C \|\tilde{q}^\varepsilon\|_{L^2}^2\nonumber\\
&\quad+C\big(\|\mathcal{R}_0^\varepsilon\|^2_{L^2}+\|\mathcal{R}_2^\varepsilon\|^2_{L^2}+\|\mathcal{R}_3^\varepsilon\|^2_{L^2}\big).
\end{align}
Here, using  $\eqref{eqwiderhoqa}_1$ yields
\begin{align}\label{L2cross2}
-\varepsilon^2\langle\partial_{t} \tilde{q}^\varepsilon, \nabla \tilde{\rho}^\varepsilon\rangle&=-\varepsilon^2\frac{d}{dt}\langle  \tilde{q}^\varepsilon, \nabla \tilde{\rho}^\varepsilon\rangle+\varepsilon^2\langle \tilde{q}^\varepsilon, \nabla \partial_{t}\tilde{\rho}^\varepsilon\rangle \nonumber\\
&=-\varepsilon^2\frac{d}{dt}\langle \tilde{q}^\varepsilon, \nabla\tilde{\rho}^\varepsilon\rangle+\varepsilon^2\|\dive \tilde{q}^\varepsilon\|_{L^2}^2-\varepsilon^2\langle \dive\tilde{q}^\varepsilon,\mathcal{R}_1^\varepsilon\rangle\nonumber\\
&\leq- \varepsilon^2\frac{d}{dt}\langle  \tilde{q}^\varepsilon, \nabla\tilde{\rho}^\varepsilon\rangle+C\varepsilon^2 \|\dive \tilde{q}^\varepsilon\|_{L^2}^2+C\varepsilon^2 \|\mathcal{R}_1^\varepsilon\|_{L^2}^2.
\end{align}
Substituting \eqref{L2cross2} into \eqref{L2cross1} and integrating the resulting inequality over $[0,t]$, we have
\begin{align}\label{nablarhof}
&\quad \frac{1}{2}p'(1) \int_0^{t}\|\nabla\tilde{\rho}^\varepsilon(t')\|_{L^2}^2\,dt'\nonumber\\
&\leq  -\varepsilon^2 \langle  \tilde{q}^\varepsilon, \nabla \tilde{\rho}^\varepsilon\rangle\Big|^t_0+C\int_0^t\|\tilde{q}^\varepsilon(t') \|_{L^2}^2\,dt'+C\varepsilon^2\int_0^t \|\dive q^\varepsilon(t')\|_{L^2}^2\,dt'\nonumber\\
&\quad+C \int_0^t \big(\|\mathcal{R}_0^\varepsilon(t')\|^2_{L^2}+\|\mathcal{R}_2^\varepsilon(t')\|^2_{L^2}+\|\mathcal{R}_3^\varepsilon(t')\|^2_{L^2}+\varepsilon^2 \|\mathcal{R}_1^\varepsilon(t')\|^2_{L^2}\big)\,dt'.
\end{align}
To analyze the first term on the right-hand side of \eqref{nablarhof}, we note that, due to \eqref{F1.8}, \eqref{F1.9}, and $\tilde{q}^\varepsilon|_{t=0}=\nabla p(\rho_0^*)$,
\begin{equation}\label{4150}
\begin{aligned}\nonumber
-\varepsilon^2 \langle  \tilde{q}^\varepsilon, \nabla \tilde{\rho}^\varepsilon\rangle\Big|^t_0&\leq \varepsilon^2 \| \tilde{q}^\varepsilon\|_{L^2}^2+C\varepsilon^2\big(\|\rho^\varepsilon-1\|_{H^1}^2+\|\rho^*-1\|_{H^1}^2\big)\\
&\quad+\varepsilon^2\|\nabla p(\rho_0^*)\|_{L^2}\big(\|\nabla\rho_0^\varepsilon\|_{L^2}+\|\nabla\rho_0^*\|_{L^2}\big)\\
&\leq \varepsilon^2 \| \tilde{q}^\varepsilon\|_{L^2}^2+C\varepsilon^2 .
\end{aligned}
\end{equation}
 Moreover, one has
\begin{align}\label{4151}
&\quad\int_0^t \|\dive \tilde{q}^\varepsilon(t')\|_{L^2}^2\,dt'\nonumber\\
&\leq C \int_0^t\big(\|\dive q^\varepsilon(t')\|_{L^2}^2+\|\dive q^*(t')\|_{L^2}^2+\|\dive q^\varepsilon_I(t')\|_{L^2}^2\big)\,dt'\leq C.
\end{align}
In view of \eqref{errorr1}-\eqref{errorr3}, \eqref{caseall}, \eqref{nablarhof} and \eqref{4151}, it holds that
\begin{equation}\label{nablarho}
\begin{aligned}
&\int_0^t\|\nabla\tilde{\rho}^\varepsilon(t')\|_{L^2}^2\,dt'\leq C\varepsilon^2\sup_{t'\in[0,t]}  \| \tilde{q}^\varepsilon(t')\|_{L^2}^2+C\int_0^t\|\tilde{q}^\varepsilon(t') \|_{L^2}^2\,dt'+C\varepsilon^2.
\end{aligned}
\end{equation}
Plugging \eqref{nablarho} into $\eqref{tilderhoL2}$ gives rise to
\begin{equation*}
\begin{aligned}
&\quad \sup_{t'\in[0,t]}\big(\|\tilde{\rho}^\varepsilon(t')\|_{L^2}^2
+\varepsilon^2\| \tilde{q}^\varepsilon(t')\|_{L^2}^2\big)
+\int_0^t\|\tilde{q}^\varepsilon(t')\|_{L^2}^2\,dt'\\
&\leq C \|\rho_0^\varepsilon-\rho_0^*\|_{L^2}^2+C\varepsilon^2\\
&\quad+ C(\delta+\delta^*)\Big(\varepsilon^2\sup_{t'\in[0,t]}  \| \tilde{q}^\varepsilon(t')\|_{L^2}^2+\int_0^t\|\tilde{q}^\varepsilon(t') \|_{L^2}^2\,dt'+\varepsilon^2\Big).
\end{aligned}
\end{equation*}
Since $\delta$ and $\delta^*$ are suitably small, we derive
\begin{equation*}
\begin{aligned}
&\sup_{t'\in[0,t]}\big(\|\tilde{\rho}^\varepsilon(t')\|_{L^2}^2
+\varepsilon^2\| \tilde{q}^\varepsilon(t')\|_{L^2}^2\big)
+\int_0^t\|\tilde{q}^\varepsilon(t')\|_{L^2}^2\,dt'\leq C \|\rho_0^\varepsilon-\rho_0^*\|_{L^2}^2+ C\varepsilon^2,
\end{aligned}
\end{equation*}
which, together with \eqref{nablarho}, yields \eqref{2.31} and finishes the proof of Lemma \ref{L2.1}. \hfill $\Box$






\subsection{Higher-order error estimates}
We have the following lemma.
\begin{lemma}
\label{L2.2}
    It holds
\begin{align}\label{F2.9}
&\quad\sup_{t\in\mathbb{R}^+}\big(\|\nabla \tilde{\rho}^\varepsilon(t)\|_{H^{m-2}}^2+\varepsilon^2\|\nabla\tilde{q}^\varepsilon(t)\|_{H^{m-2}}^2\big)\nonumber\\
&\quad\quad+\int_0^\infty \big(\|\nabla^2 \tilde{\rho}^\varepsilon(t)\|_{H^{m-2}}^2+\|\nabla\tilde{q}^\varepsilon(t)\|_{H^{m-2}}^2 \big)\,dt\nonumber\\
&\le C\|\rho_0^\varepsilon-\rho_0^*\|_{H^{m-1}}^2+C\varepsilon^2,
\end{align}
where $C>0$ is a constant independent of $\varepsilon$.
\end{lemma}

\noindent {\bf Proof.} We perform similar energy estimates as in Lemma \ref{L2.1}. 
Compared with the computations in Lemma \ref{L2.1}, we need to treat $\mathcal{R}_i^\varepsilon$ ($i=0,1,2,3$) in a more careful manner. 
Let $\al\in\N^d$ with $ 1\leq |\al|\leq m-1$. Applying $\partial^{\alpha}$ to \eqref{eqwiderhoqa}, we have
\begin{equation}
\left\{
\label{F2.10}
\begin{aligned}
  &\partial_t \partial^{\alpha}\tilde{\rho}^\varepsilon+\dive \partial^{\alpha} \tilde{q}^\varepsilon=\partial^{\alpha} \mathcal{R}_1^\varepsilon,\\
& \partial_t \partial^{\alpha} \tilde{q}^\varepsilon+p'(1)\na\partial^{\alpha}\tilde{\rho}^\varepsilon+\partial^{\alpha}\tilde{q}^\varepsilon=\partial^{\alpha} \mathcal{R}^\varepsilon_0 +\partial^{\alpha}\mathcal{R}^\varepsilon_2+ \partial^{\alpha}\mathcal{R}^\varepsilon_3.
\end{aligned}
\right.
\end{equation}
This implies the energy equality
\begin{align}\label{F2.111f}
& \frac{1}{2}\frac{d}{dt}\Big( p'(1)\| \partial^{\alpha}\tilde{\rho}^\varepsilon\|_{L^2}^2+\varepsilon^2\|\partial^{\alpha}\tilde{q}^\varepsilon\|_{L^2}^2 \Big)+\| \partial^{\alpha}\tilde{q}^\varepsilon\|_{L^2}^2\nonumber\\
=&p'(1)\langle \partial^{\alpha}\mathcal{R}_1^\varepsilon, \partial^{\alpha} \tilde{\rho}^\varepsilon\rangle+\langle  \partial^{\alpha} \mathcal{R}^\varepsilon_0+\partial^{\alpha} \mathcal{R}^\varepsilon_2+\partial^{\alpha}\mathcal{R}^\varepsilon_3, \partial^{\alpha}\tilde{q}^\varepsilon\rangle.
\end{align}
We now handle the terms on the right-hand side of \eqref{F2.111f}. 
First, the term involving $\mathcal{R}^\varepsilon_0$ is analyzed by
$$
\langle  \partial^{\alpha} \mathcal{R}^\varepsilon_0, \partial^{\alpha}\tilde{q}^\varepsilon\rangle\leq \frac{1}{4}\| \partial^{\alpha}\tilde{q}^\varepsilon\|_{L^2}^2+\|\partial^{\alpha} \mathcal{R}^\varepsilon_0\|_{L^2}^2,
$$
with
\begin{align*}
\|\partial^{\alpha} \mathcal{R}^\varepsilon_0\|_{L^2}^2&\le 2\big\|\partial^{\alpha}\big( (p'(\rho^\varepsilon)-p'(\rho^*))\na \rho^*\big)\big\|_{L^2}^2+2\big\|\partial^{\alpha}\big((p'(\rho^*)-p'(1))\na \tilde{\rho}^\varepsilon\big)\big\|_{L^2}^2.
\end{align*}
It follows from \eqref{F1.9} and {\color{blue}\eqref{F1.1311}} that
\begin{equation*}
\begin{aligned}
&\quad\big\|\partial^{\alpha}\big( (p'(\rho^*)-p'(1)\big)\na \tilde{\rho}^\varepsilon\big)\big\|_{L^2}^2\\
&\leq C\|p'(\rho^*)-p'(1)\|_{L^{\infty}}^2\| \partial^{\alpha}\na \tilde{\rho}^\varepsilon\|_{L^2}^2+C\|\partial^{\alpha} p'(\rho^*)\|_{L^{2}}^2 \|\na \tilde{\rho}^\varepsilon\|_{L^{\infty}}^2\\
&\leq C \|\rho^*-1\|_{H^m}^2  \|\na \tilde{\rho}^\varepsilon\|_{H^{m-1}}^2\\
&\leq C\delta^* \|\na \tilde{\rho}^\varepsilon\|_{L^2}^2+C\delta^* \|\na^2 \tilde{\rho}^\varepsilon\|_{H^{m-2}}^2.
\end{aligned}
\end{equation*}
Similarly, we have
\begin{equation*}
\begin{aligned}
&\big\|\partial^{\alpha}\big((p'(\rho^\varepsilon)-p'(\rho^*))\na  \rho^*\big)\|_{L^2}^2\leq C\|\tilde{\rho}^\varepsilon\|_{L^{\infty}}^2\|\partial^{\alpha} \nabla\rho^*\|_{L^2}^2+C\|\partial^{\alpha}\tilde{\rho}^\varepsilon\|_{L^{2}}^2\|\nabla \rho^*\|_{L^{\infty}}^2.
\end{aligned}
\end{equation*}
To analyze the terms on the right-hand side of the above inequality, we need to consider $d\geq3$ and $d=1,2$,  separately. In the case $d\geq3$, according to the Gagliardo-Nirenberg inequality (Lemma \ref{Lt9}), for any $f\in H^{m-1} $ and $\theta\in(0,1)$ such that $0=\theta(\frac{d}{2}-1)+(1-\theta)(\frac{d}{2}-m+1)$, we discover
$$
\|f \|_{L^{\infty}}\leq C\|f\|_{\dot{H}^1}^{\theta} \|f\|_{\dot{H}^{m-1}}^{1-\theta}\leq C\|\nabla f\|_{H^{m-2}},
$$
from which we infer
\begin{equation*}
\begin{aligned}
\big\|\partial^{\alpha}\big((p'(\rho^\varepsilon)-p'(\rho^*))\na  \rho^*\big)\big\|_{L^2}^2&\leq  C\| \nabla\rho^*\|_{H^{m-1}} ^2\|\nabla\tilde{\rho}^\varepsilon\|_{H^{m-2}}^2.
\end{aligned}
\end{equation*}
In the case $d=1,2$, we also have
$$
\|f \|_{L^{\infty}}\leq C\|f\|_{L^2}^{1-\frac{d}{4}} \|f\|_{\dot{H}^{2}}^{\frac{d}{4}},
$$
so using \eqref{F1.9} and $1\leq |\alpha|\leq m-1$, we have
\begin{equation*}
\begin{aligned}
&\quad \big\|\partial^{\alpha}\big( (p'(\rho^\varepsilon)-p'(\rho^*))\na  \rho^* \big)\big\|_{L^2}^2\\
 &\leq  C\|\tilde{\rho}^\varepsilon\|_{L^{2}}^{2-\frac{d}{2}} \|\nabla^2\tilde{\rho}^\varepsilon\|_{L^{2}}^{\frac{d}{2}}\| \nabla \partial^{\alpha} \rho^*\|_{L^2}^2+ C\|\partial^{\alpha}\tilde{\rho}^\varepsilon\|_{L^{2}}^2\|\nabla \rho^*\|_{H^{m-1}}^2\\
 &\leq C\| \nabla\rho^*\|_{H^{m-1}} ^2\big(\|\tilde{\rho}^\varepsilon\|_{L^{2}}^2+\|\nabla^2\tilde{\rho}^\varepsilon\|_{L^2}^2\big)+C \|\nabla \rho^*\|_{H^{m-1}}^2 \|\nabla\tilde{\rho}^\varepsilon\|_{H^{m-1}}^2 \\
 &\leq C \| \nabla\rho^*\|_{H^{m-1}}^2 \|\tilde{\rho}^\varepsilon\|_{L^{2}}^2+C\delta^* \|\nabla \tilde{\rho}^\varepsilon\|_{L^2}^2 +C \delta^* \|\nabla^2 \tilde{\rho}^\varepsilon\|_{H^{m-2}}^2.
\end{aligned}
\end{equation*}
Hence, combining the above two cases $d\geq3$ and $d=1,2$ and using the error estimate \eqref{2.31}, we arrive at
\begin{align}\label{2.26}
\int_0^t\|\partial^{\alpha} \mathcal{R}^\varepsilon_0(t')\|_{L^2}^2\,dt'& \leq  C \sup_{t'\in[0,t]}\|\tilde{\rho}^\varepsilon(t')\|_{L^{2}}^2\int_0^t\| \nabla\rho^*(t')\|_{H^{m-1}}^2\,dt' \nonumber\\
&\quad+C\delta^*\int_0^t \|\nabla \tilde{\rho}^\varepsilon(t')\|_{L^2}^2\,dt' +C \delta^* \int_0^t\|\nabla^2 \tilde{\rho}^\varepsilon(t')\|_{H^{m-2}}^2\,dt'\nonumber\\
&\leq C\|\rho_0^\varepsilon-\rho_0^*\|_{L^2}^2+C \varepsilon^2 +C\delta^*\int_0^t\|\nabla^2 \tilde{\rho}^\varepsilon(t')\|_{H^{m-2}}^2\,dt'.
\end{align}
Note that the integration associated with $\partial^{\alpha}\mathcal{R}^\varepsilon_i$ ($i=1,2$) can be handled similarly as in Lemma \ref{L2.1}:
\begin{align}
\int_0^t\langle  \partial^{\alpha} \mathcal{R}^\varepsilon_2, \partial^{\alpha}\tilde{q}^\varepsilon\rangle \,dt'&\leq \frac{1}{4}\int_0^t\|\partial^{\alpha}\tilde{q}^\varepsilon(t')\|_{L^2}^2\,dt'+C\int_0^t\|\partial^{\alpha} \mathcal{R}^\varepsilon_2(t')\|_{L^2}^2\,dt'\nonumber\\
&\leq \frac{1}{4}\int_0^t\|\partial^{\alpha}\tilde{q}^\varepsilon(t')\|_{L^2}^2\,dt'+C\varepsilon^2,\label{2.27}
\end{align}
and
\begin{align}\label{4.24}
\int_0^t\langle  \partial^{\alpha} \mathcal{R}^\varepsilon_1, \partial^{\alpha}\tilde{\rho}^\varepsilon\rangle \,dt'&\leq 2\sup_{t'\in [0,t]}\|\partial^{\alpha}\tilde{\rho}^\varepsilon(t')\|_{L^2}\int_0^t\| \mathcal{R}^\varepsilon_1(t')\|_{H^{m-1}} \,dt'\nonumber\\
&\leq \frac{1}{4}\sup_{t'\in [0,t]}\|\partial^{\alpha}\tilde{\rho}^\varepsilon(t')\|_{L^2}^2+C\Big(\int_0^t\| \mathcal{R}^\varepsilon_1(t')\|_{H^{m-1}} \,dt'\Big)^2\nonumber\\
&\leq \frac{1}{4}\sup_{t'\in [0,t]}\|\partial^{\alpha}\tilde{\rho}^\varepsilon(t')\|_{L^2}^2+C\varepsilon^2.
\end{align}
However, the term involving $\partial^{\alpha}\mathcal{R}^\varepsilon_3=-\varepsilon^2\partial^{\alpha}\partial_t q^*$ requires a more elaborate analysis. Indeed, $\partial^{\alpha}\mathcal{R}^\varepsilon_3$ may not be bounded in $L^2(\mathbb{R}^+; L^{2})$ for $|\alpha|=m-1$ but obeys a faster rate (see  \eqref{errorr3}). To address this difficulty, we observe that
$$
\tilde{q}^\varepsilon=q^\varepsilon-q^*-e^{-\frac{t}{\varepsilon^2}} \frac{1}{\varepsilon}\rho_0^\varepsilon u_0^\varepsilon.
$$
Recall the operator $\Lambda$ defined in Section \ref{section3}. By integration by parts, it holds
\begin{equation}\nonumber
\begin{aligned}
&\quad\langle  \partial^{\alpha}\mathcal{R}^\varepsilon_3, \partial^{\alpha}\tilde{q}^\varepsilon\rangle\\
&=\langle \Lambda^{-1} \partial^{\alpha}\mathcal{R}^\varepsilon_3, \Lambda\partial^{\alpha} q^\varepsilon \rangle+\frac{\varepsilon^2}{2}\frac{d}{dt}\|\partial^{\alpha} q^*\|_{L^2}^2+\varepsilon e^{-\frac{t}{\varepsilon^2}} \langle  \Lambda^{-1} \partial^{\alpha}\mathcal{R}^\varepsilon_3, \Lambda\partial^{\alpha} (\rho_0^\varepsilon u_0^\varepsilon)\rangle.
\end{aligned}
\end{equation}
This, together with \eqref{F1.8}, \eqref{F1.9}-\eqref{F1.91}, \eqref{errorr3} and the regularities of initial data, leads to
\begin{align}
&\quad\int_0^t \langle  \partial^{\alpha}\mathcal{R}^\varepsilon_3, \partial^{\alpha}\tilde{q}^\varepsilon\rangle \,dt'\nonumber\\
&\leq \frac{\varepsilon^2}{2}\| q^*\|_{H^{m-1}}^2+\frac{\varepsilon^2}{2}\| \nabla p(\rho_0^*)\|_{H^{m-1}}^2\nonumber\\
&\quad+\Big(\int_0^t \|\mathcal{R}^\varepsilon_3(t')\|_{H^{m-2}}^2\,dt' \Big)^{\frac{1}{2}}\Big( \int_0^t \Big(\|q^\varepsilon(t')\|_{H^m}^2+\varepsilon^2 e^{-\frac{2t'}{\varepsilon^2}}\|\rho_0^\varepsilon u_0^\varepsilon\|_{H^m}^2\Big)\,dt'\Big)^{\frac{1}{2}}\leq C\varepsilon^2 .\label{high2euler0}
\end{align}
Recall that
\[  \varepsilon^2\|\partial^{\alpha}\tilde{q}^\varepsilon|_{t=0}\|_{L^2}^2\leq \varepsilon^2\|\nabla p(\rho^*)|_{t=0}\|_{H^{m-1}}^2\leq C\varepsilon^2.
\]
Integrating \eqref{F2.111f} over $[0,t]$ and substituting \eqref{2.26} and \eqref{2.27}-\eqref{high2euler0} into the resulting inequality, we have
\begin{align}\label{errorHm1}
&\quad\sup_{t'\in [0,t]}\big(\|\nabla \tilde{\rho}^\varepsilon(t)\|_{H^{m-2}}^2+\varepsilon^2\|\nabla\tilde{q}^\varepsilon(t)\|_{H^{m-2}}^2\big)+\int_0^t \|\nabla\tilde{q}^\varepsilon(t')\|_{H^{m-2}}^2\,dt' \nonumber \\
&\leq C\|\rho_0^\varepsilon-\rho_0^*\|_{H^{m-1}}^2+C\varepsilon^2+C\delta^*\int_0^t \|\nabla^2 \tilde{\rho}^\varepsilon(t')\|_{H^{m-2}}^2\,dt'.
\end{align}

Let us now turn to capturing the higher-order dissipation of $\tilde{\rho}^\varepsilon$ required in \eqref{errorHm1}. We rewrite $\eqref{F2.10}_2$ by
\begin{align}
& p'(1)\na\partial^{\alpha}\tilde{\rho}^\varepsilon+\partial^{\alpha}\tilde{q}^\varepsilon=\partial^{\alpha} \mathcal{R}^\varepsilon_0 +\partial^{\alpha}\mathcal{R}^\varepsilon_2-\partial_t \partial^{\alpha} q^\varepsilon.\label{F2.10fff}
\end{align}
Multiplying $\eqref{F2.10fff}$ by $\nabla\partial^{\alpha}\tilde{\rho}^\varepsilon$ and using $\partial_t \partial^{\alpha} \tilde{\rho}^\varepsilon=-\partial^{\alpha}\dive(q^\varepsilon-q^*)$ and integration by parts, we deduce that, after direct computations,
\begin{equation}
\begin{aligned}\nonumber
p'(1)\|\nabla\partial^{\alpha}\tilde{\rho}^\varepsilon\|_{L^2}^2&=-\varepsilon^2\frac{d}{dt}\langle \partial^{\alpha} q^\varepsilon, \nabla \partial^{\alpha}\tilde{\rho}^\varepsilon\rangle-\varepsilon^2 \langle \dive \partial^{\alpha} q^\varepsilon, \partial_t \partial^{\alpha}\tilde{\rho}^\varepsilon\rangle-\langle \partial^{\alpha}\tilde{q}^\varepsilon, \nabla\partial^{\alpha}\tilde{\rho}^\varepsilon\rangle\\
&\quad +\langle \partial^{\alpha}\mathcal{R}^\varepsilon_0
+\partial^{\alpha}\mathcal{R}^\varepsilon_2, \nabla\partial^{\alpha}\tilde{\rho}^\varepsilon\rangle\\
&\leq -\varepsilon^2\frac{d}{dt}\langle \partial^{\alpha} q^\varepsilon, \nabla \partial^{\alpha}\tilde{\rho}^\varepsilon\rangle+C\varepsilon^2\big(\|\partial^{\alpha} \dive q^*\|_{L^2}^2+\| \partial^{\alpha} \dive q^\varepsilon\|_{L^2}^2\big)\\
&\quad+\frac{1}{2}p'(1) \|\nabla\partial^{\alpha}\tilde{\rho}^\varepsilon\|_{L^2}^2+C\|\partial^{\alpha}\tilde{q}^\varepsilon\|_{L^2}^2\\
&\quad+C\big(\|\partial^{\alpha}\mathcal{R}_0^\varepsilon\|_{L^2}^2+\|\partial^{\alpha}\mathcal{R}_2^\varepsilon\|_{L^2}^2\big).
\end{aligned}
\end{equation}
This leads to
\begin{align}
\frac{1}{2}p'(1)\int_0^t \|\nabla\partial^{\alpha}\tilde{\rho}^\varepsilon(t')\|_{L^2}^2\,dt'&\leq  -\varepsilon^2 \langle \partial^{\alpha} q^\varepsilon, \nabla \partial^{\alpha}\tilde{\rho}^\varepsilon\rangle\Big|^t_0+C\int_0^t \|\partial^{\alpha}\tilde{q}^\varepsilon(t')\|_{L^2}^2\,dt'\nonumber\\
&\quad+C\varepsilon^2 \int_0^t \big(\| q^\varepsilon(t')\|_{H^{m}}^2+\| q^*(t')\|_{H^{m}}^2\big)\,dt'\nonumber\\
&\quad+C\int_0^t \big(\|\mathcal{R}_0^\varepsilon(t')\|_{H^{m-1}}^2+\|\mathcal{R}_2^\varepsilon(t')\|_{H^{m-1}}^2\big)\,dt'.\label{sgggb}
\end{align}
Due to $\varepsilon \|q^\varepsilon\|_{H^m}+\varepsilon \|q_0^\varepsilon\|_{H^m}\leq C$, the first term on the right-hand side of \eqref{sgggb} is bounded by
\begin{align}
-\varepsilon^2 \langle \partial^{\alpha} q^\varepsilon, \nabla \partial^{\alpha}\tilde{\rho}^\varepsilon\rangle\Big|^t_0&=\varepsilon^2 \langle \partial^{\alpha} \dive q^\varepsilon, \partial^{\alpha}\tilde{\rho}^\varepsilon\rangle\Big|^t_0\nonumber\\
&\leq \varepsilon \big( \varepsilon \|q^\varepsilon\|_{H^m}\big) \|\nabla \tilde{\rho}^\varepsilon\|_{H^{m-2}}+\varepsilon \big( \varepsilon \|q^\varepsilon_0\|_{H^m}\big) \|\nabla \tilde{\rho}^\varepsilon|_{t=0}\|_{H^{m-2}}\nonumber\\
&\leq C\| \nabla\tilde{\rho}^\varepsilon\|_{H^{m-2}}^2+C\|\nabla(\rho_0^\varepsilon-\rho_0^*)\|_{H^{m-2}}^2+C\varepsilon^2.\label{sggg}
\end{align}
In addition, from the regularities of $ (\rho^\varepsilon,q^\varepsilon)$ and $(\rho^*,q^*)$ one has
\begin{equation}\label{sggg1}
\begin{aligned}
\int_0^t \big(\|q^\varepsilon(t')\|_{H^m}^2+\|q^*(t')\|_{H^m}^2\big)\,dt'\leq C ,
\end{aligned}
\end{equation}
and
\begin{equation}\label{sggg10}
\begin{aligned}
\int_0^t \| \nabla\tilde{\rho}^\varepsilon(t')\|_{H^{m}}^2\,dt'\leq C\int_0^t \big(\|\nabla \rho^\varepsilon(t')\|_{H^m}^2+\|\nabla\rho^*(t')\|_{H^m}^2\big)\,dt'\leq C .
\end{aligned}
\end{equation}
Combining \eqref{errorr1}-\eqref{errorr2}, \eqref{2.26}, \eqref{sgggb}-\eqref{sggg10} and using the smallness of $\delta^*$, we arrive at
\begin{equation}
\begin{aligned}\label{nablarhohigh}
 \int_0^t \|\nabla^2 \tilde{\rho}^\varepsilon(t')\|_{H^{m-2}}^2\,dt'&\leq  C\sup_{t\in[0,t]} \|\nabla \tilde{\rho}^\varepsilon(t')\|_{H^{m-1}}^2+C\|\nabla(\rho_0^\varepsilon-\rho_0^*)\|_{H^{m-2}}^2+C\varepsilon^2.
\end{aligned}
\end{equation}
Plugging \eqref{nablarhohigh} into \eqref{errorHm1}, using the smallness of $\delta^*$ and returning to  \eqref{nablarhohigh},  we get the desired estimate \eqref{F2.9} and finish the proof of Lemma \ref{L2.2}.  \hfill $\Box$

\vspace{2mm}

\section{The Euler system (II): Proof of Theorem \ref{T1.2}}

\subsection{The error equations}

In this section, we prove Theorem \ref{T1.2} and show a faster convergence rate between $(\rho^\varepsilon, q^\varepsilon)$ and its first-order asymptotic expansion
$(\rho^\varepsilon_a, q^\varepsilon_a)$ with
$$
\rho^\varepsilon_a=\rho^*+\varepsilon \rho_1\quad \text{and}\quad q^\varepsilon_a=q^*+\varepsilon q_1,
$$
where the profiles $\rho_1$ and $q_1$ solve the system \eqref{rho1q1}. According to \eqref{F1.5} and \eqref{rho1q1}, $(\rho^\varepsilon_a, q^\varepsilon_a)$ satisfies
\begin{equation}\label{rhoaqa}
\left\{
\begin{aligned}
& \partial_t\rho^\varepsilon_a+\dive q^\varepsilon_a=0,\\
&q^\varepsilon_a =-\nabla \big(p(\rho^*)+\varepsilon p'(\rho^*)\rho_1\big).
\end{aligned}
\right.
\end{equation}
A key ingredient is to establish the error estimates for $(\tilde{\rho}^\varepsilon_a,\tilde{q}^\varepsilon_a)$ given by
\begin{align}
 \tilde{\rho}_a^\varepsilon:=\rho^\varepsilon-\rho^\varepsilon_a,\quad\quad \tilde{q}_a^\varepsilon:=q^\varepsilon- q^\varepsilon_a.\label{erroraEuler}
\end{align}
By a direct computation, the equations for $(\tilde{\rho}^\varepsilon_a,\tilde{q}^\varepsilon_a)$ read
\begin{equation}\label{Eulertildea}
\left\{
\begin{aligned}
& \partial_t\tilde{\rho}_a^\varepsilon+\dive \tilde{q}_a^\varepsilon=0,\\
&\tilde{q}_a^\varepsilon=-\nabla \big(p(\rho^\varepsilon)-p(\rho^\varepsilon_a) \big)+\mathcal{R}_4^\varepsilon-\varepsilon^2\partial_t q^\varepsilon,
\end{aligned}
\right.
\end{equation}
where
\begin{align}
\mathcal{R}_4^\varepsilon=-\varepsilon^2 \dive \Big( \frac{q^\varepsilon\otimes q^\varepsilon}{\rho^\varepsilon}\Big)-\nabla \big(p(\rho^\varepsilon_a)-p(\rho^*)-\varepsilon p'(\rho^*)\rho_1\big).\label{R3R4}
\end{align}
Before analyzing $(\tilde{\rho}^\varepsilon,\tilde{q}^\varepsilon)$, we shall derive regularity estimates for $(\rho_1, q_1)$.
\subsection{Estimates for $(\rho_1,q_1)$}
\begin{lemma}\label{Lrho1q1}
There exists a pair $(\rho_1,q_1)$ solving \eqref{rho1q1}-\eqref{rho1q1d} with the initial data $\rho_1|_{t=0}=\rho_{1,0}\in H^{m}$. In addition, there exists a constant $C$ depending only on $\rho_0^*$ and $\rho_{1,0}$ such that
\begin{align}
&\sup_{t\in\mathbb{R}^{+}}\big( \|\rho_1(t)\|_{H^{m}}^2+\|q_1(t)\|_{H^{m-1}}^2\big)\nonumber\\
&\quad+\int_{0}^{\infty}\big(\|\nabla \rho_1(t)\|_{H^{m}}^2+\|q_1(t)\|_{H^{m}}^2+\|\partial_tq_1(t)\|_{H^{m-2}}^2\big) dt\leq C.\label{rho1q1:es}
\end{align}
\end{lemma}

\noindent {\bf Proof.}
Standard theory for linear parabolic equations ensures that the equation \eqref{rho1} with $\rho_1|_{t=0}=\rho_{1,0}\in H^{m}$ admits a unique global solution $\rho_1\in C(\mathbb{R}^{+}; H^m)$. Direct computations on \eqref{rho1} yield that, for all $\alpha\in \mathbb{R}^d$ such that $0\leq |\alpha|\leq m$,
\begin{align}\label{rho1ES}
&\frac{1}{2}\frac{d}{dt}\|\partial^{\alpha}\rho_1\|_{L^2}^2+p'(1)\|\partial^{\alpha}\nabla\rho_1\|_{L^2}^2\nonumber\\
\leq& \big\|\partial^{\alpha}\nabla \big( (p'(\rho^*)-p'(1))\rho_1\big)\big\|_{L^2} \|\partial^{\alpha}\nabla \rho_1\big\|_{L^2}\nonumber\\
\leq& \frac{p'(1)}{2} \|\partial^{\alpha}\nabla\rho_1\|_{L^2}^2+C\big\|\partial^{\alpha} \nabla \big( (p'(\rho^*)-p'(1))\rho_1\big)\big\|_{L^2}^2.
\end{align}
Using \eqref{F1.12} and \eqref{F1.15}, we have
\begin{equation*}
\begin{aligned}
&\quad\sum_{0\leq |\alpha|\leq m}\big\| \nabla\partial^{\alpha} \big( (p'(\rho^*)-p'(1))\rho_1\big)\big\|_{L^2}\\
&\leq C\|\nabla p'(\rho^*)\rho_1\|_{H^{m}}+C\| (p'(\rho^*)-p'(1))\nabla \rho_1\|_{H^{m}}\\
&\leq C\|\nabla p'(\rho^*)\|_{H^{m}}\|\rho_1\|_{H^m}+\|p'(\rho^*)-p'(1)\|_{H^m}\|\nabla\rho_1\|_{H^m}\\
&\leq C\|\rho_1\|_{H^{m}} \|\nabla \rho^*\|_{H^{m}}+C \|\rho^*-1\|_{H^{m}} \|\nabla \rho_1\|_{H^{m}} .
\end{aligned}
\end{equation*}
Therefore, integrating \eqref{rho1ES} in time and recalling \eqref{F1.9}, we arrive at
\begin{equation*}
\begin{aligned}
&\quad\sup_{t'\in[0,t]} \|\rho_1(t')\|_{H^{m}}^2+\int_{0}^{t}\|\nabla\rho_1(t')\|_{H^{m}}^2\,dt'\\
&\leq C\|\rho_{1,0}\|_{H^m}^2+C\sup_{t'\in[0,t]}\|\rho_1(t')\|_{H^{m}}^{2} \int_0^t\|\nabla \rho^*(t')\|_{H^m}^2\,dt'\\
&\quad+C\sup_{t'\in[0,t]}\|\rho^*(t')-1\|_{H^{m}}^{2} \int_0^t\|\nabla \rho_1(t')\|_{H^m}^2\,dt'\\
&\leq C\|\rho_{1,0}\|_{H^m}^2+C\delta^*  \Big(\sup_{t'\in[0,t]}\|\rho_1(t)\|_{H^{m}}^2+\int_{0}^{t}\|\nabla\rho_1(t')\|_{H^{m}}^2\,dt'\Big).
\end{aligned}
\end{equation*}
Consequently, as the constant $\delta^*>0$  given by Proposition \ref{P1.2} is suitably small, we obtain the desired bounds for  $\rho_1$:
\begin{equation}\label{rho1:es}
\begin{aligned}
&\quad\sup_{t\in\mathbb{R}^{+}} \|\rho_1(t)\|_{H^{m}}^2+\int_{0}^{\infty}\|\nabla\rho_1(t)\|_{H^{m}}^2\,dt\leq C\|\rho_{1,0}\|_{H^m}^2+C\|\rho^*-1\|_{H^m}^2.
\end{aligned}
\end{equation}

Next, we establish the estimates for $q_1$. It follows from  \eqref{rho1q1d} and \eqref{F1.12} that
\begin{equation*}
\begin{aligned}
\|q_1\|_{H^{k}}&\leq C\|\nabla \rho_1\|_{H^{k}}+C\|\nabla (p'(\rho^*)-p'(1))\rho_1)\|_{H^{k}}\\
&\leq C\|\nabla \rho_1\|_{H^{k}}+\|\nabla p'(\rho^*)\rho_1\|_{H^{k}}+\| (p'(\rho^*)-p'(1))\nabla\rho_1\|_{H^{k}}\\
&\leq C\|\nabla \rho_1\|_{H^{k}}+C\|\nabla \rho^*\|_{H^{m-1}} \|\rho_1\|_{H^k}+C \|\nabla \rho^*\|_{H^{k}} \|\rho_1\|_{H^{m-1}}\\
&\quad+ C\|\rho^*-1\|_{H^{m-1}}\|\nabla \rho_1\|_{H^k}+C\|\rho^*-1\|_{H^{k}}\|\nabla \rho_1\|_{H^{m-1}},
\end{aligned}
\end{equation*}
for $k=1,...,m$. Hence, owing to \eqref{F1.9} and \eqref{rho1:es}, $q_1$ satisfies
\begin{equation*}
\begin{aligned}
\sup_{t\in\mathbb{R}^+}\|q_1(t)\|_{H^{m-1}}&\leq C\sup_{t\in\mathbb{R}^+}\big(1+\|\rho^*(t)-1\|_{H^{m}}\big) \|\rho_1(t)\|_{H^m}\leq C,
\end{aligned}
\end{equation*}
and
\begin{equation*}
\begin{aligned}
\int_0^{\infty}\|q^*(t)\|_{H^{m}}^2\,dt \leq & C\big(1+ \sup_{t\in\mathbb{R}^{+}}\|\rho^*(t)-1\|_{H^{m}}^2\big)\int_0^{\infty}\|\nabla \rho_1(t)\|_{H^{m}}^2\,dt\\
&\quad+C\sup_{t\in\mathbb{R}^{+}}\|\rho_1(t)\|_{H^{m}}^2\int_0^{\infty}\|\nabla \rho^*(t)\|_{H^{m}}^2\,dt\leq C.
\end{aligned}
\end{equation*}
Finally, since
$$
\partial_tq_1=-\na\big(p''(\rho^*) \partial_t\rho^* \rho_1 +p'(\rho^*)\partial_t\rho_1\big)=\na\big(p''(\rho^*)\dive q^*\rho_1+p'(\rho^*)\dive q_1\big),
$$
it follows that
\begin{equation*}
\begin{aligned}
\int_0^{\infty}\|\partial_tq_1(t)\|_{H^{m-2}}^2\,dt&\leq C\int_0^{\infty}\big(\|p''(\rho^*)\dive q^*\rho_1(t)\|_{H^{m-1}}^2+\|p'(\rho^*)\dive q_1(t)\|_{H^{m-1}}^2\big)\,dt\\
&\leq C\big(1+\sup_{t\in\mathbb{R}^+}\|\rho^*(t)-1\|_{H^{m-1}}^2\big)\sup_{t\in\mathbb{R}^+}\|\rho_1(t)\|_{H^{m-1}}^2\int_0^{\infty}\|q^*(t)\|_{H^m}^2\,dt\\
&\quad+ C\big(1+\sup_{t\in\mathbb{R}^+}\|\rho^*(t)-1\|_{H^{m-1}}^2\big)\int_0^{\infty}\|q_1(t)\|_{H^m}^2\,dt\leq C.
\end{aligned}
\end{equation*}
This completes the proof of Lemma \ref{Lrho1q1}. \hfill $\Box$

\subsection{Error estimates for $q^\varepsilon-q^*$}
Assuming further that \eqref{wellpre} holds, we can prove an improved estimate for the error $q^\varepsilon-q^*$, which will be useful for deriving the faster convergence rate for the error $(\tilde{\rho}^\varepsilon_a,\tilde{q}^\varepsilon_a)$.

\begin{lemma}
\label{L2.3}
Under the conditions of Theorem \ref{T1.2}, it holds
 \begin{equation}\label{2.610}
 \left\{
\begin{aligned}
& \sup_{t\in\mathbb{R}^+}\|(q^\varepsilon-q^*)(t)\|_{H^{m-1}}^2\leq C,\\
&\int_0^\infty  \|(q^\varepsilon-q^*)(t)\|_{H^{m-1}}^2 \,dt \leq C\varepsilon^2,
\end{aligned}
\right.
\end{equation}
for a uniform constant $C>0$.
\end{lemma}

\noindent {\bf Proof.} Since $\|q_0^\varepsilon-q_0^*\|_{H^{m-1}}\leq \varepsilon$, we deduce that  $q_0^\varepsilon$ is uniformly bounded in $H^{m-1}$ with respect to $\varepsilon$. As a consequence,
\begin{align*}
\int_0^\infty  \|q_I^\varepsilon(t)\|_{H^{m-1}}^2 \,dt\leq \int_0^{\infty} e^{-\frac{2t}{\varepsilon^2}} \,dt \|q_0^\varepsilon\|_{H^{m-1}}^2\leq C\varepsilon^2,
\end{align*}
which, together with \eqref{F1.11a} and $\|\rho_0^\varepsilon-\rho_0^*\|_{H^{m-1}}\leq C(1+\|\rho_{1,0}\|_{H^{m-1}})\varepsilon$, gives rise to
\begin{align*}
\int_0^\infty  \|(q^\varepsilon-q^*)(t)\|_{H^{m-1}}^2 \,dt\leq C\int_0^\infty  \big(\|(q^\varepsilon-q^*-q_I^\varepsilon)(t)\|_{H^{m-1}}^2 + \|q_I^\varepsilon(t)\|_{H^{m-1}}^2 \big)\,dt\leq C\varepsilon^2.
\end{align*}
This yields the second estimate in $\eqref{2.610}$. Meanwhile, the first estimate can be easily deduced from \eqref{2.31}, \eqref{F2.9} and the definition of $q_I^\varepsilon$. \hfill $\Box$

\subsection{Estimates for $\mathcal{R}^\varepsilon_4$}
We have the following lemma pertaining to the $\mathcal{O}(\varepsilon^2)$-control of $\mathcal{R}_4$.
\begin{lemma}\label{lemmareminderEuler}
It holds
\begin{align}
&\int_{0}^{\infty} \|\mathcal{R}_4(t)\|_{H^{m-2}}^2\,dt\leq C\varepsilon^4,\label{errorr4}
\end{align}
where $C>0$ is a constant independent of $\varepsilon$.
\end{lemma}

\noindent {\bf Proof.} We first handle the term $\varepsilon^2 \dive \Big( \frac{q^\varepsilon\otimes q^\varepsilon}{\rho^\varepsilon}\Big)$. In view of the regularities for $\rho^\varepsilon, q^\varepsilon$ and $q^*$ (see \eqref{F1.8} and \eqref{F1.91}), $q^\varepsilon=q^\varepsilon-q^*+q^*$ and the improved estimate \eqref{2.610}, it follows that
\begin{align}
&\quad\int_0^\infty\Big\|\dive \Big( \frac{q^\varepsilon\otimes q^\varepsilon}{\rho^\varepsilon}\Big)(t)\Big\|_{H^{m-2}}^2\,dt\nonumber\\
&\leq C\big(1+\sup_{t\in\mathbb{R}^+}\|\rho^\varepsilon(t)-1\|_{H^m}^2\big)\sup_{t\in\mathbb{R}^+}\big( \|(q^\varepsilon-q^*)(t)\|_{H^{m-1}}^2+ \|q^*(t)\|_{H^{m-1}}^2\big)\int_0^{\infty}\|q^\varepsilon(t)\|_{H^m}^2\,dt\leq C.\nonumber
\end{align}
To analyze the second term in $\mathcal{R}^\varepsilon_4$, we write
\begin{equation*}
\begin{aligned}
&p(\rho^\varepsilon_a)-p(\rho^*)-\varepsilon p'(\rho^*)\rho_1= \varepsilon^2 \Pi^\varepsilon(\rho^*,\rho^\varepsilon_a)\rho_1^2,
\end{aligned}
\end{equation*}
where
\begin{equation*}
\begin{aligned}
  \Pi^\varepsilon(\rho^*,\rho_1):=\int_0^1 p''(\theta\rho^*+(1-\theta)\rho^\varepsilon_a)\theta d\theta,\quad \rho_a^\varepsilon=\rho^*+\varepsilon \rho_1.
\end{aligned}
\end{equation*}
This implies that
\begin{equation*}
\begin{aligned}
&\nabla\big( p(\rho^\varepsilon_a)-p(\rho^*)-\varepsilon p'(\rho^*)\rho_1\big)\\
&=\varepsilon^2\big(\partial_{\rho^*}\Pi^\varepsilon(\rho^*,\rho_1)\nabla\rho^*+ \partial_{\rho_1}\Pi^\varepsilon(\rho^*,\rho_1)\nabla\rho_1\big)\rho^2_1
+2\varepsilon^2\Pi^\varepsilon(\rho^*,\rho_1) \rho_1 \nabla \rho_1.
\end{aligned}
\end{equation*}
Therefore, \eqref{F1.15} as well as \eqref{F1.9} and  \eqref{rho1q1:es} guarantee that  
\begin{equation}\label{R40000}
\begin{aligned}
&\int_0^{\infty}\big\|\nabla \big( p(\rho^\varepsilon_a)-p(\rho^*)-\varepsilon p'(\rho^*)\rho_1\big)(t)\big\|_{H^m}^2\,dt\\
&\quad \leq C\varepsilon^4\int_0^{\infty}\big(\|\nabla \rho^*(t)\|_{H^m}^2+\|\nabla \rho_1(t)\|_{H^m}^2\big)\,dt\leq C\varepsilon^4.
\end{aligned}
\end{equation}
The combination of \eqref{errorr3} and \eqref{R40000} leads to \eqref{errorr4}. \hfill $\Box$

\subsection{Improved error estimates}
$\,$
\vspace{1mm}

We are going to establish faster convergence rates for $\rho^\varepsilon-\rho^*-\varepsilon \rho_1$ and $q^\varepsilon-q^*-\varepsilon q_1-q_I^\varepsilon$, under the condition \eqref{wellpre}. 
Recall that $(\rho^\varepsilon,q^\varepsilon)$ and $(\rho^*,q^*)$ correspond to the solutions for the Euler equations and the porous medium equation obtained in Propositions \ref{P1.1} and  \ref{P1.2}, respectively.

First, we establish $L^2$ estimates for $\rho^\varepsilon-\rho^*-\varepsilon \rho_1$ and $q^\varepsilon-q^*-\varepsilon q_1-q_I^\varepsilon$ with a rate of order  $\varepsilon^2$. 

\begin{lemma}
\label{L2.7}
 It holds
\begin{align}\label{2.3111}
&\sup_{t\in\mathbb{R}^+}\big(\|\tilde{\rho}^\varepsilon_a(t)\|_{L^2}^2+\varepsilon^2\|\tilde{q}^\varepsilon_a(t)\|_{L^2}^2\big)+\int_0^\infty \big(\|\nabla \tilde{\rho}^\varepsilon_a(t)\|_{L^2}^2+\|\tilde{q}^\varepsilon_a(t)\|_{L^2}^2 \big)\,dt\le C\varepsilon^4,
\end{align}
where the error $(\tilde{\rho}^\varepsilon_a,\tilde{q}^\varepsilon_a)$ is defined by \eqref{erroraEuler}, and $C>0$ is a constant independent of $\varepsilon$.
\end{lemma}

\noindent {\bf Proof.}
The system \eqref{Eulertildea} can be reformulated as the inhomogeneous parabolic flow:
\begin{equation}\label{eqwiderhoqa11}
\begin{aligned}
   \partial_t\tilde{\rho}_a^\varepsilon-\dive( p'(\rho^\varepsilon) \nabla \tilde{\rho}_a^\varepsilon) =\dive \big((p'(\rho^\varepsilon)-p'(\rho^\varepsilon_a))\na\rho^\varepsilon_a\big)-\dive \mathcal{R}_4^\varepsilon+\varepsilon^2\partial_t\dive  q^\varepsilon.
\end{aligned}
\end{equation}
Taking the $L^2$-inner product of \eqref{eqwiderhoqa11} with $\tilde{\rho}_a^\varepsilon$ yields
\begin{align}\label{22222}
&\frac{1}{2}\frac{d}{dt}\|\tilde{\rho}_a^\varepsilon\|_{L^2}^2+ \langle p'(\rho^\varepsilon) \nabla \tilde{\rho}_a^\varepsilon, \nabla \tilde{\rho}_a^\varepsilon\rangle\nonumber\\
=&-\langle (p'(\rho^\varepsilon)-p'(\rho^\varepsilon_a)) \na\rho^\varepsilon_a, \nabla \tilde{\rho}_a^\varepsilon\rangle -\langle \mathcal{R}_4^\varepsilon, \nabla \tilde{\rho}_a^\varepsilon\rangle-\varepsilon^2\langle \partial_t  q^\varepsilon, \nabla \tilde{\rho}_a^\varepsilon\rangle.
\end{align}
Since $\rho^\varepsilon$ is in a neighborhood of $1$, the condition \eqref{F1.1} implies that there exists a constant $p_1$ such that
\begin{equation}
\begin{aligned}
 \langle p'(\rho^\varepsilon) \nabla \tilde{\rho}_a^\varepsilon, \nabla \tilde{\rho}_a^\varepsilon\rangle\geq p_1 \|\nabla \tilde{\rho}_a^\varepsilon\|_{L^2}^2.\label{2222221}
\end{aligned}
\end{equation}
Concerning the first term on the right-hand side of \eqref{22222}, note that
$$
p'(\rho^\varepsilon)-p'(\rho^\varepsilon_a)=\pi(\rho^\varepsilon, \rho^\varepsilon_a)\tilde{\rho}^\varepsilon_a,\quad\quad
\pi(\rho^\varepsilon,\rho^\varepsilon_a):=\int_0^1 p''(\theta\rho^\varepsilon+(1-\theta)\rho^\varepsilon_a) d\theta.
$$
In the case $d\geq 3$, one may employ the Cauchy-Schwarz inequality together with Sobolev's embedding and the uniform upper bounds for $\rho^\varepsilon$, $\rho^*$, and $\rho_1$. Consequently, we have
\begin{align}\label{case11}
 -\langle(p'(\rho^\varepsilon)-p'(\rho^\varepsilon_a))\na\rho^\varepsilon_a,\nabla \tilde{\rho}_a^\varepsilon\rangle
&\le C\|\pi(\rho^\varepsilon, \rho^\varepsilon_a)\|_{L^{\infty}}\|\tilde{\rho}_a^\varepsilon\|_{L^{\frac{2d}{d-2}}}
\|\nabla \rho^\varepsilon_a\|_{L^d} \|\nabla \tilde{\rho}_a^\varepsilon\|_{L^2}\nonumber\\
&\le C\|\nabla \tilde{\rho}_a^\varepsilon\|_{L^2} \|\nabla \rho^\varepsilon_a\|_{H^{m-1}} \|\nabla \tilde{\rho}_a^\varepsilon\|_{L^2}\nonumber\\
&\leq C(\sqrt{\delta^*}+\varepsilon) \|\nabla \tilde{\rho}^\varepsilon_a\|_{L^2}^2.
\end{align}
where we used 
$$
\|\nabla \rho^\varepsilon_a\|_{H^{m-1}}\leq C\|\rho^*-1\|_{H^m}+\varepsilon\|\rho_1\|_{H^m}\leq C(\sqrt{\delta^*}+\varepsilon).
$$
In the case $d=2$, we take advantage of the Gagliardo-Nirenberg-Sobolev inequality to derive
\begin{align}\label{case21}
& -\langle(p'(\rho^\varepsilon)-p'(\rho^\varepsilon_a))\na\rho^\varepsilon_a,\nabla \tilde{\rho}_a^\varepsilon\rangle\nonumber\\
&\le C\|\pi(\rho^\varepsilon, \rho^\varepsilon_a)\|_{L^{\infty}}\|\tilde{\rho}^\varepsilon_a\|_{L^4} \|\nabla \rho^\varepsilon_a\|_{L^4} \|\nabla \tilde{\rho}_a^\varepsilon\|_{L^2}\nonumber\\
&\le C\|\tilde{\rho}_a^\varepsilon\|_{L^{2}}^{\frac{1}{2}}\|\nabla\tilde{\rho}_a^\varepsilon\|_{L^{2}}^{\frac{1}{2}} \|\nabla \rho^\varepsilon_a\|_{L^2}^{\frac{1}{2}}\|\nabla^2 \rho^\varepsilon_a\|_{L^2}^{\frac{1}{2}}  \|\nabla \tilde{\rho}_a^\varepsilon\|_{L^2}\nonumber\\
&\leq \frac{1}{8}p_1  \|\nabla \tilde{\rho}_a^\varepsilon\|_{L^2}^2+C\|\nabla \rho^\varepsilon_a\|_{H^1}^4\|\tilde{\rho}_a^\varepsilon\|_{L^2}^2.
\end{align}
 Finally,  the case $d=1$ can be treated analogously:
\begin{align}
& -\langle(p'(\rho^\varepsilon)-p'(\rho^\varepsilon_a))\na\rho^\varepsilon_a,\nabla\tilde{\rho}_a^\varepsilon\rangle\nonumber\\
&\le C\|\pi(\rho^\varepsilon, \rho^\varepsilon_a)\|_{L^{\infty}}\|\tilde{\rho}^\varepsilon_a\|_{L^{\infty}} \|\nabla \rho^\varepsilon_a\|_{L^2} \|\nabla\tilde{\rho}_a^\varepsilon\|_{L^2}\nonumber\\
&\le C\|\tilde{\rho}_a^\varepsilon\|_{L^{2}}^{\frac{1}{2}}\|\nabla\tilde{\rho}_a^\varepsilon\|_{L^{2}}^{\frac{1}{2}} \|\nabla \rho^\varepsilon_a\|_{L^2} \|\nabla\tilde{\rho}_a^\varepsilon\|_{L^2}\nonumber\\
&\leq \frac{1}{8}p_1 \|\nabla \tilde{\rho}_a^\varepsilon\|_{L^2}^2+C\|\nabla \rho^\varepsilon_a\|_{L^2}^4\|\tilde{\rho}_a^\varepsilon\|_{L^2}^2.\label{case31}
\end{align}
Also, one has
\begin{equation}\label{2222222}
\begin{aligned}
&-\langle \mathcal{R}_4^\varepsilon, \nabla \tilde{\rho}_a^\varepsilon\rangle\leq \frac{1}{8}p_1 \|\nabla\tilde{\rho}_a^\varepsilon\|_{L^2}^2+C \|\mathcal{R}_4^\varepsilon\|_{L^2}^2.
\end{aligned}
\end{equation}
For the last term on the right-hand side of \eqref{22222}, we have to handle the singularity of $\partial_t q^\varepsilon $. To this matter, for some constant $\nu>0$ to be chosen later, we use the equation for $\tilde{\rho}_a^\varepsilon$ and obtain
\begin{align}
-\varepsilon^2\langle \partial_t  q^\varepsilon, \nabla\tilde{\rho}_a^\varepsilon\rangle&=\varepsilon^2\langle \partial_t  \dive \tilde{q}_a^\varepsilon,  \tilde{\rho}_a^\varepsilon\rangle+\varepsilon^2\langle \partial_t  \dive (q^*+\varepsilon q_1), \tilde{\rho}_a^\varepsilon\rangle\nonumber\\
&=\varepsilon^2\frac{d}{dt}\langle \dive  \tilde{q}_a^\varepsilon, \tilde{\rho}_a^\varepsilon\rangle+\varepsilon^2\langle   \tilde{q}_a^\varepsilon, \nabla\partial_t \tilde{\rho}_a^\varepsilon\rangle-\varepsilon^2\langle \partial_t (q^*+\varepsilon q_1),\nabla  \tilde{\rho}_a^\varepsilon\rangle\nonumber\\
&\leq \frac{1}{8}p_1\|\nabla \tilde{\rho}_a^\varepsilon\|_{L^2}^2+\varepsilon^2\frac{d}{dt}\langle \dive  \tilde{q}_a^\varepsilon, \tilde{\rho}_a^\varepsilon\rangle+ C\varepsilon^4\|\partial_t(q^*+\varepsilon q_1)\|_{L^2}^2\nonumber\\
&\quad+\nu\|\tilde{q}_a^\varepsilon\|_{L^2}^2+C\nu^{-1}\varepsilon^4 \|\nabla\partial_t\tilde{\rho}^\varepsilon_a\|_{L^2}^2 .\label{2222223}
\end{align}
Combining \eqref{22222} with \eqref{2222221}-\eqref{2222223} and using $\|\nabla \rho^\varepsilon_a\|_{H^{m-1}}\leq C(\sqrt{\delta^*}+\varepsilon)\leq C$  yields
\begin{equation}\nonumber
\begin{aligned}
&\frac{1}{2}\frac{d}{dt}\|\tilde{\rho}_a^\varepsilon\|_{L^2}^2+ (p_- -C\sqrt{\delta^*}-C\varepsilon)\|\nabla \tilde{\rho}_a^\varepsilon\|_{L^2}^2\\
\leq &-\varepsilon^2\frac{d}{dt}\langle \dive  \tilde{q}_a^\varepsilon, \tilde{\rho}_a^\varepsilon\rangle+C\|\nabla \rho^\varepsilon\|_{H^{m}}^2\|\tilde{\rho}_a^\varepsilon\|_{L^2}^2\\
&\quad+C \|\mathcal{R}_4^\varepsilon\|_{L^2}^2+ C\varepsilon^4\|\partial_t(q^*+\varepsilon q_1)\|_{L^2}^2+\nu\|\tilde{q}_a^\varepsilon\|_{L^2}^2+C\nu^{-1}\varepsilon^4 \|\nabla\partial_t\tilde{\rho}^\varepsilon_a\|_{L^2}^2.
\end{aligned}
\end{equation}
In fact, one can let $\varepsilon\leq \varepsilon_0$ with $\varepsilon_0<<1$ and make use of the uniform bounds for $\tilde{\rho}_a^\varepsilon$ to obtain \eqref{2.3111} when $\varepsilon_0\leq \varepsilon\leq 1$. Exploiting the smallness of $\delta^*, \varepsilon$ as well as the fact that $\tilde{\rho}_a^\varepsilon|_{t=0}=\rho_0^\varepsilon-\rho_0^*-\varepsilon \rho_{1,0}=\mathcal{O}(\varepsilon^2)$ in $H^{m-2}$, we get
\begin{align}
&\quad \sup_{t'\in[0,t]} \|\tilde{\rho}_a^\varepsilon(t')\|_{L^2}^2 + C \int_0^t \|\nabla\tilde{\rho}_a^\varepsilon(t')\|_{L^2}^2\,dt'\nonumber\\
&\leq C\varepsilon^4+\varepsilon^2\langle \dive  \tilde{q}_a^\varepsilon, \tilde{\rho}_a^\varepsilon\rangle \Big|^t_0+C\int_0^t\|\mathcal{R}_4^\varepsilon(t')\|_{L^2}^2 \,dt'\nonumber\\
&\quad+C\sup_{t'\in[0,t]}\| \rho^\varepsilon(t')\|_{H^m}^2\int_0^t\|\nabla \rho^\varepsilon(t')\|_{H^{m-1}}^2\|\tilde{\rho}_a^\varepsilon(t')\|_{L^2}^2\,dt'+C\nu\int_0^t \|\tilde{q}_a^\varepsilon(t')\|_{L^2}^2\,dt'\nonumber\\
&\quad+ C\varepsilon^4\int_0^t\big(\|\partial_t(q^*+\varepsilon q_1)(t')\|_{L^2}^2+\nu^{-1}\|\nabla\partial_t\tilde{\rho}^\varepsilon_a(t')\|_{L^2}^2\big)\,dt'.\label{rhoatilde}
\end{align}
With the uniform bounds \eqref{F1.8}, \eqref{F1.9}, \eqref{2151} and \eqref{rho1q1:es} at hand, it holds that
\begin{equation}\label{mml1}
\begin{aligned}
\varepsilon^4\int_0^t\|\partial_t(q^*+\varepsilon q_1)(t')\|_{H^{m-2}}^2\,dt'\leq C\varepsilon^4 .
\end{aligned}
\end{equation}
As
$$\partial_t\tilde{\rho}^\varepsilon_a=\partial_t (\rho^\varepsilon-\rho^*-\varepsilon \rho_1)=\dive(-q^\varepsilon+q^*+\varepsilon q_1),  $$
we also have
\begin{equation}\label{mml111}
\begin{aligned}
&\quad\varepsilon^4\int_0^t \|\nabla\partial_t\tilde{\rho}^\varepsilon_a(t')\|_{H^{m-2}}^2\,dt'\nonumber\\
&\leq C\varepsilon^4\int_0^t \big(\|q^\varepsilon(t')\|_{H^{m}}^2+\|q^*(t')\|_{H^m}^2+\|q_1(t')\|_{H^m}^2\big)\,dt'\leq C\varepsilon^4.
\end{aligned}
\end{equation}
Since $\|\tilde{\rho}_a^\varepsilon\|_{H^{m-2}}\leq \varepsilon^2$ and $\|\tilde{q}_a^\varepsilon|_{t=0}\|_{H^{m-1}}\leq \varepsilon+\varepsilon \|q_{1,0}\|_{H^{m-1}}$ one gets from \eqref{wellpre} and \eqref{2.610} that
\begin{align}
\varepsilon^2\langle \dive  \tilde{q}_a^\varepsilon, \tilde{\rho}_a^\varepsilon\rangle \Big|^t_0&\leq \frac{1}{4}\sup_{t'\in[0,t]}\|\tilde{\rho}_a^\varepsilon(t')\|_{L^2}^2+\varepsilon^4\sup_{t'\in[0,t]}\big(\|(q^\varepsilon-q^*)(t')\|_{H^1}^2+\varepsilon^2\|q_1(t')\|_{H^1}^2\big)\nonumber\\
&\quad+\varepsilon^5 \big(1+\|\dive q_{1,0}\|_{L^2}\big)\leq C\varepsilon^4.\label{mml11}
\end{align}
It thus follows from \eqref{F1.8}, \eqref{errorr4} and \eqref{rhoatilde}-\eqref{mml11} that
\begin{equation}
\label{rhoatilde0}
\sup_{t'\in[0,t]}\|\tilde{\rho}_a^\varepsilon(t')\|_{L^2}^2+ \int_0^t
\|\tilde{\rho}_a^\varepsilon(t')\|_{L^2}^2\,dt'\leq C(1+\nu^{-1})\varepsilon^4+\nu\int_0^t
\|\tilde{q}_a^\varepsilon(t')\|_{L^2}^2\,dt'.
\end{equation}

We then turn to analyze the error $\tilde{q}_a^\varepsilon$ that satisfies a damped equation with sources:
\begin{equation}\label{Q}
\begin{aligned}
&\varepsilon^2\partial_t \tilde{q}_a^\varepsilon+ \tilde{q}_a^\varepsilon=-\nabla (p(\rho^\varepsilon)-p(\rho^\varepsilon_a))+\mathcal{R}^\varepsilon_4-\varepsilon^2(\partial_t q^*+\varepsilon \partial_t q_1).
\end{aligned}
\end{equation}
Taking the $L^2$-inner product of \eqref{Q} with $\tilde{q}^\varepsilon_a$ yields
\begin{equation*}
\begin{aligned}
&\frac{\varepsilon^2}{2}\frac{d}{dt}\|\tilde{q}_a^\varepsilon\|_{L^2}^2+\|\tilde{q}_a^\varepsilon\|_{L^2}^2\\
\leq &\big(\|\nabla (p(\rho^\varepsilon)-p(\rho^\varepsilon_a))\|_{L^2}+\|\mathcal{R}^\varepsilon_4\|_{L^2}+\varepsilon^2\|(\partial_t q^*+\varepsilon \partial_t q_1)\|_{L^2}\big)\|\tilde{q}^\varepsilon_a\|_{L^2}\\
\leq & \frac{1}{2}\|\tilde{q}^\varepsilon_a\|_{L^2}^2+C \|\nabla \tilde{\rho}^\varepsilon_a\|_{L^2}^2+C\|\mathcal{R}^\varepsilon_4\|_{L^2}^2+C\varepsilon^4\|(\partial_t q^*+\varepsilon \partial_t q_1)\|_{L^2}^2.
\end{aligned}
\end{equation*}
Since $\tilde{q}_a^\varepsilon|_{t=0}=\mathcal{O}(\varepsilon)$ in $H^{m-1}$, this leads to
\begin{equation}
\begin{aligned}
&\varepsilon^2\sup_{t'\in[0,t]}\|\tilde{q}^\varepsilon_a(t')\|_{L^2}^2+\int_0^t\|\tilde{q}^\varepsilon_a(t')\|_{L^2}^2\,dt'\\
&\quad\leq C\varepsilon^4+ C\int_0^t \|\nabla \tilde{\rho}^\varepsilon_a(t')\|_{L^2}^2\,dt'+C\int_0^t\|\mathcal{R}^\varepsilon_4(t')\|_{L^2}^2\,dt'\\
&\quad\quad+C\varepsilon^4\int_0^t\|(\partial_t q^*+\varepsilon \partial_t q_1)(t')\|_{L^2}^2\,dt'\\
&\leq C (1+\nu^{-1})\varepsilon^4+C\nu \int_0^t\|\tilde{q}^\varepsilon_a(t')\|_{L^2}^2\,dt',\label{Q1525}
\end{aligned}
\end{equation}
where \eqref{errorr4}, \eqref{mml1} and $\eqref{rhoatilde0}$ have been used. Choosing $\nu$ sufficiently small and combining \eqref{Q1525} with \eqref{rhoatilde0}, we prove \eqref{2.3111}. \hfill $\Box$





\subsection{Higher-order error estimates}
We have the following lemma.
\begin{lemma}
\label{L2.8}
 Let $d\geq2$ be such that $m\geq 3$. Then the following estimate holds:
\begin{align}
&\sup_{t\in\mathbb{R}^+}\big(\|\nabla \tilde{\rho}_a^\varepsilon(t)\|_{H^{m-3}}^2+\varepsilon^2\|\nabla \tilde{q}_a^\varepsilon(t)\|_{H^{m-3}}^2\big)\nonumber\\
&\quad\quad+\int_0^\infty \big(\|\nabla^2 \tilde{\rho}_a^\varepsilon(t)\|_{H^{m-3}}^2+\|\nabla \tilde{q}_a^\varepsilon(t)\|_{H^{m-3}}^2 \big)\,dt\le C\varepsilon^4,\label{F2.900}
\end{align}
where $C>0$ is a constant independent of $\varepsilon$.
\end{lemma}

\noindent {\bf Proof.}
Applying $\partial^{\alpha}$ with $1\leq |\alpha|\leq m-2$ to \eqref{eqwiderhoqa11}, we have
\begin{align}
   &\partial_t\partial^{\alpha}\tilde{\rho}_a^\varepsilon- p'(1) \Delta \partial^{\alpha}\tilde{\rho}_a^\varepsilon\nonumber \\
   &\quad=\dive\partial^{\alpha}\Big( (p'(\rho^\varepsilon)-p'(1)) \nabla\tilde{\rho}_a^\varepsilon+ (p'(\rho^\varepsilon)-p'(\rho^\varepsilon_a))\na\rho^\varepsilon_a-\mathcal{R}_4^\varepsilon+\varepsilon^2\partial_t \partial^{\alpha} q^\varepsilon\Big).\label{eqwiderhoqa111}
\end{align}
This implies the following energy inequality:
\begin{align}
&\frac{1}{2}\frac{d}{dt}\| \partial^{\alpha}\tilde{\rho}^\varepsilon_a\|_{L^2}^2+p'(1)\| \partial^{\alpha}\nabla\tilde{\rho}^\varepsilon_a\|_{L^2}^2\nonumber\\
=&-\Big\langle\Big(\partial^{\alpha}\big( (p'(\rho^\varepsilon)-p'(1)) \nabla\tilde{\rho}_a^\varepsilon+ (p'(\rho^\varepsilon)-p'(\rho^\varepsilon_a))\na\rho^\varepsilon_a-\mathcal{R}_4^\varepsilon+\varepsilon^2\partial_t \partial^{\alpha} q^\varepsilon\Big),\partial^{\alpha}\nabla \tilde{\rho}_a^\varepsilon \Big\rangle \nonumber\\
\leq& \frac{1}{2}p'(1)\| \partial^{\alpha}\nabla\tilde{\rho}^\varepsilon_a\|_{L^2}^2+C\big\|\partial^{\alpha}\big( (p'(\rho^\varepsilon)-p'(1)) \nabla\tilde{\rho}_a^\varepsilon\big)\big\|_{L^2}^2\nonumber\\
&\quad+C\big\| \partial^{\alpha}\big((p'(\rho^\varepsilon)-p'(\rho^\varepsilon_a))\na\rho^\varepsilon_a\big)\big\|_{L^2}^2+C\|\partial^{\alpha}\mathcal{R}_4^\varepsilon\|_{L^2}^2-\varepsilon^2\langle \partial_t \partial^{\alpha} q^\varepsilon, \nabla \partial^{\alpha}\tilde{\rho}^\varepsilon_a\rangle.\label{F2.1110}
\end{align}
We now handle the terms on the right-hand side of \eqref{F2.1110}. 
The term $\big\|\partial^{\alpha}\big( (p'(\rho^\varepsilon)-p'(1)) \nabla\tilde{\rho}_a^\varepsilon\big)\big\|_{L^2}^2$ is analyzed in the following three cases.

\begin{itemize}
\item Case 1: $m\geq [\frac{d}{2}]+3$ or $m$ is even when $m=[\frac{d}{2}]+2$ for $d\geq 2$.
\end{itemize}

By \eqref{F1.8} and the Moser-type inequality  \eqref{F1.1311}, we have
\begin{align*}
&\quad\big\|\partial^{\alpha}\big( (p'(\rho^\varepsilon)-p'(1)) \nabla\tilde{\rho}_a^\varepsilon\big)\big\|_{L^2}^2\\
& \leq C\|p'(\rho^\varepsilon)-p'(1)\|_{L^{\infty}}^2 \|\partial^{\alpha}\nabla\tilde{\rho}_a^\varepsilon\|_{L^2}^2+C\|\partial^{\alpha}(p'(\rho^\varepsilon)-p'(1))\|_{L^2}^2 \|\partial^{\alpha}\nabla\tilde{\rho}_a^\varepsilon\|_{L^{\infty}}^2.
\end{align*}
In this case, note that $m-2>\frac{d}{2}$ implies  the embedding $H^{m-2} \hookrightarrow L^{\infty} $. Consequently, by \eqref{F1.8} and \eqref{F1.15}, for $1\leq |\alpha|\leq m-2$ we obtain
\begin{align*}
&\big\|\partial^{\alpha}\big( (p'(\rho^\varepsilon)-p'(1)) \nabla\tilde{\rho}_a^\varepsilon\big)\big\|_{L^2}^2\nonumber\\
\leq & C\|p'(\rho^\varepsilon)-p'(1)\|_{H^{m-2}}^2 \|\partial^{\alpha}\nabla\tilde{\rho}_a^\varepsilon\|_{L^2}^2+C\|\partial^{\alpha}(p'(\rho^\varepsilon)-p'(1))\|_{L^2}^2 \|\nabla\tilde{\rho}_a^\varepsilon\|_{H^{m-2}}^2\nonumber\\
\leq & C\|\rho^\varepsilon -1\|_{H^{m}}^2\|\nabla^2\tilde{\rho}_a^\varepsilon\|_{H^{m-2}}^2\nonumber\\
\leq & C \delta \|\nabla^2\tilde{\rho}_a^\varepsilon\|_{H^{m-2}}^2.
\end{align*}

\begin{itemize}
\item Case 2: $m=[\frac{d}{2}]+2$ and $m$ is odd for $d\geq 3$.
\end{itemize}

In this case, we have $[\frac{d}{2}]=\frac{d}{2}-\frac{1}{2}$ and $m-2<\frac{d}{2}$. It follows from \eqref{F1.8} and the product law \eqref{product} that
\begin{align*}
\big\|\partial^{\alpha}\big( (p'(\rho^\varepsilon)-p'(1)) \nabla\tilde{\rho}_a^\varepsilon\big)\big\|_{L^2}^2&\leq C\|p'(\rho^\varepsilon)-p'(1)\|_{H^{m-1}}^2 \|\partial^{\alpha}\nabla\tilde{\rho}_a^\varepsilon\|_{L^2}^2\nonumber\\
&\leq  C \delta \|\partial^{\alpha}\nabla\tilde{\rho}_a^\varepsilon\|_{L^2}^2.
\end{align*}

\begin{itemize}
\item Case 3: $m=3$ for $d=2$.
\end{itemize}

Since $m-2=1$, we know $|\alpha|=1$. Thus, it is easy to verify that
\begin{align*}
&\big\|\partial^{\alpha}\big( (p'(\rho^\varepsilon)-p'(1)) \nabla\tilde{\rho}_a^\varepsilon\big)\big\|_{L^2}^2\nonumber\\
\leq & C\| p'(\rho^\varepsilon)-p'(1)\|_{L^{\infty}}^2\|\nabla^2\tilde{\rho}_a^\varepsilon \|_{L^2}^2+C\| \nabla(p'(\rho^\varepsilon)-p'(1))\|_{L^{\infty}}^2\|\nabla\tilde{\rho}_a^\varepsilon \|_{L^2}^2\nonumber\\
\leq & C\delta \|\nabla^2\tilde{\rho}_a^\varepsilon \|_{L^2}^2+C\|\nabla \rho^\varepsilon\|_{H^{m}}^2  \|\nabla\tilde{\rho}_a^\varepsilon \|_{L^2}^2.
\end{align*}
Combining the above three cases, we know
\begin{equation*}
\begin{aligned}
\big\|\partial^{\alpha}\big( (p'(\rho^\varepsilon)-p'(1)) \nabla\tilde{\rho}_a^\varepsilon\big)\big\|_{L^2}^2\leq C\delta \|\nabla^2\tilde{\rho}_a^\varepsilon\|_{H^{m-2}}^2+C\|\nabla \rho^\varepsilon\|_{H^{m}}^2 \| \nabla \tilde{\rho}_a^\varepsilon\|_{H^{m-2}}^2.
\end{aligned}
\end{equation*}
With respect to the last term $-\varepsilon^2\langle \partial^{\alpha}\partial_t \partial^{\alpha} q^\varepsilon, \nabla \partial^{\alpha}\tilde{\rho}^\varepsilon_a\rangle$, we use a similar argument as in \eqref{2222223} to get
\begin{align*}
& -\varepsilon^2\langle \partial^{\alpha}\partial_t   q^\varepsilon, \partial^{\alpha}\nabla \tilde{\rho}_a^\varepsilon\rangle&\nonumber\\
\leq &\frac{1}{8}p'(1)\|\nabla  \partial^{\alpha}\tilde{\rho}_a^\varepsilon\|_{L^2}^2+\varepsilon^2\frac{d}{dt}\langle \dive  \partial^{\alpha} \tilde{q}_a^\varepsilon,  \partial^{\alpha}\tilde{\rho}_a^\varepsilon\rangle\\
&\quad+\nu\|\partial^{\alpha} \tilde{q}_a^\varepsilon\|_{L^2}^2+C\nu^{-1}\varepsilon^4\|\partial_t\nabla\tilde{\rho}_a^\varepsilon\|_{L^2}^2+ C\varepsilon^4\|\partial_t \partial^{\alpha}(q^*+\varepsilon q_1)\|_{L^2}^2 \nonumber.
\end{align*}
Recall that $\delta$ and $\delta^*$ are suitably small, and both $\partial_t (q^*+\varepsilon q_1)$ and $\partial_t \nabla\tilde{\rho}^\varepsilon_a$ are uniformly bounded in $L^2(\mathbb{R}^+;H^{m-2})$ due to \eqref{mml1} and \eqref{mml11}.  Then, after a computation similar to \eqref{Q1525}, we integrate \eqref{F2.1110} in time, make use of \eqref{errorr4} and \eqref{mml1}, and then derive
\begin{equation*}
\begin{aligned}
&\sup_{t\in\mathbb{R}^+}\|\nabla\tilde{\rho}_a^\varepsilon(t)\|_{H^{m-2}}^2+\int_0^\infty \|\nabla^2 \partial^{\alpha}\tilde{\rho}_a^\varepsilon(t)\|_{L^2}^2\,dt \le C(1+\nu^{-1})\varepsilon^4+\nu\int_0^{\infty}\|\nabla\tilde{q}_a^\varepsilon(t)\|_{H^{m-2}}^2\,dt.
\end{aligned}
\end{equation*}
Returning to the damped equation \eqref{Q} and performing similar calculations for \eqref{Q1525}, we also can have the desired estimates for $\tilde{q}_a^\varepsilon$, with a remainder involving the bound for $\tilde{\rho}_a^\varepsilon$. Finally, choosing a suitably small $\nu>0$ leads to the expected convergence rate; the details are omitted.
\hfill $\Box$

\vspace{2mm}

\section{The Euler-Maxwell system: Proof of Theorem \ref{TEM}}


\subsection{The error equations}\label{subsection51}
$\,$
$\,$

Let $(\rho^\varepsilon, u^\varepsilon, E^\varepsilon, B^\varepsilon)$ and $(\rho^*,u^*,E^*,B^*)$ be the global solutions to \eqref{EMeps} and \eqref{F3.3}, respectively, given by
Propositions \ref{existenceEM} and \ref{existenceDD}. Introducing the error unknowns
\[
(\tilde{\rho}^\varepsilon,\tilde{q}^\varepsilon, \tilde{E}^\varepsilon,\tilde{B}^\varepsilon)=(\rho^\varepsilon-\rho^*,q^\varepsilon-q^*,E^{\varepsilon}-E^*,B^\varepsilon-B^e),
\]
we have the error system
\begin{equation}
\left\{
\label{eqerrorEM}
\begin{aligned}
  & \partial_t\tilde{\rho}^\varepsilon+\dive \tilde{q}^\varepsilon=0,\\[1mm]
&\varepsilon^2 \partial_{t}\tilde{q}^\varepsilon+\tilde{q}^\varepsilon+\tilde{\rho}^\varepsilon E^\varepsilon+\rho^*  \tilde{E}^\varepsilon+\na\big(p(\rho^\varepsilon)-p(\rho^*)\big)\\
&\quad=
-\varepsilon q^\varepsilon\times B^\varepsilon-\varepsilon^2 \partial_t q^*-\varepsilon^2\dive\Big(\frac{q^\varepsilon\otimes q^\varepsilon}{\rho^\varepsilon}\Big),\\
&\partial_{t}\tilde{E}^\varepsilon-\frac{1}{\varepsilon}\nabla\times \tilde{B}^\varepsilon=- \partial_{t}E^{*}+\tilde{q}^\varepsilon  ,\\
&\partial_{t}\tilde{B}^\varepsilon+\frac{1}{\varepsilon}\nabla\times \tilde{E}^\varepsilon=0,\\
&\dive \tilde{E}^\varepsilon=-\tilde{\rho}^\varepsilon,\quad\quad \dive  \tilde{B}^\varepsilon=0.
\end{aligned}
\right.
\end{equation}
The presence of the low-order term $- \partial_{t}E^{*}$ in the third equation poses a challenge in establishing the desired convergence rate for $(\tilde{\rho}^\varepsilon,\tilde{q}^\varepsilon, \tilde{E}^\varepsilon, \tilde{B}^\varepsilon)$. In order to overcome this difficulty, we construct the asymptotic expansion $$
\big(\rho^\varepsilon_a,q^\varepsilon_a,E^\varepsilon_a,B^\varepsilon_a\big)=\big(\rho^*,q^*,E^*,B^e\big)+\varepsilon\big(\rho_1,q_1,E_1,B_1\big).
$$
Looking for the equations satisfied by $(\rho_1,q_1,E_1,B_1)$, we formally see that
\begin{align*}
  p(\rho^\varepsilon_a)=p(\rho^*+\varepsilon\rho_1)
=p(\rho^*)+\varepsilon p'(\rho^*)\rho_1
+O(\varepsilon^2),\\
\rho^\varepsilon_aE^\varepsilon_a=\rho^*E^*+\varepsilon\big(\rho^*E_1+\rho_1E^*\big)
+O(\varepsilon^2),\\
 \varepsilon q^\varepsilon_a\times B^\varepsilon_a=\varepsilon q^*\times B^e
+O(\varepsilon^2).
\end{align*}
Inserting these expressions into \eqref{EMeps} and identifying the coefficients in terms of the powers of $\varepsilon$, we obtain
\begin{equation}
\label{F3.11a11}
\begin{cases}
   \partial_t\rho_1+\dive q_1=0,\\
q_1+\na\big(p'(\rho^*)\rho_1\big)=-\big(\rho^*E_1+\rho_1E^*\big)-q^*\times B^e,\\
   \na\times B_1=\partial_t E^*-q^*,\qquad \dive E_1=-\rho_1,\\
   \na\times E_1=0,\qquad\quad\quad\quad\,\,\,\, \dive B_1=0.
\end{cases}
\end{equation}
Note that $E_1$ is given by
\[  -\Delta E_1=\nabla\times\nabla E_1-\nabla\dive E_1=\nabla \rho_1,\] or, equivalently,
\begin{align}
E_1=\Lambda^{-2}\nabla \rho_1.\label{E1EM}
\end{align}
Since $(\rho^*,q^*,E^*,B^e)$ is known by \eqref{F3.3},  we have 
$$
\partial_{t}E^*=\Lambda^{-2}\nabla \partial_{t}\rho^*=-\Lambda^{-2}\nabla\dive q^*.
$$
Similarly, from \eqref{F3.11a11} and \eqref{E1EM}, it is clear that
\begin{align*}
    -\Delta B_1=\nabla\times\nabla B_1-\nabla\dive B_1=\nabla\times (\partial_{t}E^*-q^*)=-\nabla\times q^*.
\end{align*}
Consequently, we have
\begin{align}
B_{1}=-\Lambda^{-2}\nabla\times q^*.\label{B1EM}
\end{align}
Then, $\rho_1$ satisfies the drift-diffusion type equation
\begin{equation}
\label{rho1EM}
   \partial_t\rho_1-\dive\big(\na(p'(\rho^*)\rho_1)+\rho_1E^*+\rho^*\nabla \Lambda^{-2}\rho_1\big)
     =\dive (q^*\times B^e).
\end{equation}
Without loss of generality, $(\rho_1,q_1,E_1,B_1)$ is supplemented by the initial datum $\rho_1|_{t=0}=0$, which implies
\begin{align}
q_1|_{t=0}=-q^*|_{t=0} \times B^e,\quad E_1|_{t=0}=0,\quad B_1|_{t=0}=-\Lambda^{-2}\nabla\times q^*|_{t=0}.\label{5.2}
\end{align}
Then, by \eqref{EMeps} and \eqref{F3.11a11},
$(\rho^\varepsilon_a,q^\varepsilon_a,E^\varepsilon_a,B^\varepsilon_a)$ solves
\begin{equation}
\label{F3.11ca0}
\begin{cases}
   \partial_t\rho^\varepsilon_a+\dive q^\varepsilon_a=0,\\
\varepsilon^2\partial_t q^\varepsilon_a+q^\varepsilon_a+\na\big(p(\rho^*)+\varepsilon p'(\rho^*)\rho_1\big)
=-\rho^* E^\varepsilon_a-\varepsilon \rho_1 E^*-\varepsilon q^*\times  B^e-\varepsilon^2\dive\Big(\frac{q^\varepsilon\otimes q^\varepsilon}{\rho^\varepsilon}\Big),\\[1mm]
   \partial_t E^\varepsilon_a-\dfrac{1}{\varepsilon}\na\times B^\varepsilon_a
=q^\varepsilon_a+\varepsilon\big(\partial_t E_1-q_1\big),\\[3mm]
   \partial_t B^\varepsilon_a+\dfrac{1}{\varepsilon}\na\times E^\varepsilon_a=\varepsilon \partial_t B_1,\\
    \dive E^\varepsilon_a=1-\rho^\varepsilon_a\qquad
\dive B^\varepsilon_a=0.
\end{cases}
\end{equation}
Recall the initial layer correction $q_I^\varepsilon$ defined in \eqref{qIEM}. In order to derive the convergence rate of $\mathcal{O}(\varepsilon)$ for $(\tilde{\rho}^\varepsilon,\tilde{q}^\varepsilon, \tilde{E}^\varepsilon,\tilde{B}^\varepsilon)$, we turn to analyze the error 
\begin{align}
(\tilde{\rho}_a^\varepsilon,\tilde{q}_a^\varepsilon, \tilde{E}_a^\varepsilon,\tilde{B}_a^\varepsilon)=(\rho^\varepsilon-\rho^\varepsilon_a,q^\varepsilon-q^\varepsilon_a-q_I^\varepsilon,E^\varepsilon-E^\varepsilon_a,B^\varepsilon-B^\varepsilon_a),\label{tildeerrora}
\end{align}
such that \eqref{eqerrorEM} becomes
\begin{equation}
\label{F3.11ca}
\begin{cases}
   \partial_t\tilde{\rho}^\varepsilon_a+\dive \tilde{q}^\varepsilon_a=J_1^\varepsilon,\\
\varepsilon^2\partial_t\tilde{q}^\varepsilon_a+\tilde{q}^\varepsilon_a+p'(1)\nabla \tilde{\rho}_a^\varepsilon+\tilde{E}_a^\varepsilon
=J_0^\varepsilon+J_2^\varepsilon+J_3^\varepsilon+J^\varepsilon_4,\\[1mm]
   \partial_t \tilde{E}^\varepsilon_a-\dfrac{1}{\varepsilon}\na\times \tilde{B}^\varepsilon_a
=\tilde{q}^\varepsilon_a+J^\varepsilon_5,\\[3mm]
   \partial_t \tilde{B}^\varepsilon_a+\dfrac{1}{\varepsilon}\na\times \tilde{E}^\varepsilon_a=J_6^\varepsilon,\\
   \dive \tilde{E}_a^\varepsilon=-\tilde{\rho}^\varepsilon_a,\qquad \dive \tilde{B}^\varepsilon_a=0,
\end{cases}
\end{equation}
with
\begin{equation}
J_0^\varepsilon=-\big(p'(\rho^\varepsilon)-p'(\rho^\varepsilon_a)\big)\na\rho^\varepsilon_a-\big( p'(\rho^\varepsilon)-p'(1)\big)\na\tilde{\rho}_a^\varepsilon -\tilde{\rho}^\varepsilon_a E^\varepsilon-(\rho^*-1) \tilde{E}^\varepsilon_a,\label{J0}
\end{equation}
and
\begin{equation}\label{J12345}
\left\{
\begin{aligned}
J_1^\varepsilon&= \dive q_I^\varepsilon,\\
J_2^\varepsilon &=-\varepsilon^2\dive \Big(\dfrac{q^\varepsilon\otimes q^\varepsilon}{\rho^\varepsilon}\Big)-\nabla \big(p(\rho^\varepsilon_a)-p(\rho^*)-\varepsilon p'(\rho^*)\rho_1\big),\\
J_3^\varepsilon&= \varepsilon \big(\rho_1(E^\varepsilon-E^*)-q^\varepsilon\times B^\varepsilon+q^*\times  B^e),\\
J^\varepsilon_4&=-\varepsilon^2\partial_t q^*-\varepsilon^3\partial_t q_1,\\
J^\varepsilon_5&=\varepsilon(\partial_{t}E_1-q_1),\\
J_6^\varepsilon&=\varepsilon\partial_t B_1.
\end{aligned}
\right.
\end{equation}
Since all the source terms on the right-hand side of \eqref{F3.11ca} are of order $\varepsilon$, we can expect to derive $\mathcal{O}(\varepsilon)$-bounds from  \eqref{F3.11ca}.

\subsection{Regularity for the profiles}

Before deriving bounds for $(\rho_1,q_1,E_1,B_1)$, we need the $L^2$ regularity for $\rho^*$ and $E^*$. Due to $E^*=\Lambda^{-2}\nabla \rho^*=\Lambda^{-2}\nabla (\rho^*-1)$, assuming that  $\rho_0^*-1\in\dot{H}^{-1}$ seems necessary. We have the following lemma.

\begin{lemma}\label{lemma61}
Let $(\rho^*,E^*,q^*)$ be given by Proposition \ref{existenceDD} and assume that $\rho_0^*-1\in \dot{H}^{-1}$. We have \eqref{est:DD} and
\begin{align}\label{Hmin1}
&\sup_{t\in\mathbb{R}^+} \big(\|\rho^*(t)-1\|_{\dot{H}^{-1}}^2+\|q^*(t)\|_{L^2}^2+\|E^*(t)\|_{L^2}^2\big)\nonumber\\
&\quad+\int_{0}^{\infty}\big(\|\rho^*(t)-1\|_{\dot{H}^{-1}}^2+\|q^*(t)\|_{L^2}^2+\|E^*(t)\|_{L^2}^2\big)\,dt\leq C\|\rho_0^*-1\|_{ \dot{H}^{-1} \cap H^m}^2.
\end{align}
\end{lemma}

\noindent {\bf Proof.}
Recall that $\phi^*=\Lambda^{-2}(\rho^*-1)$ and $E^*=\Lambda^{-2}\nabla \rho^*$. Applying $\Lambda^{-1}$ to \eqref{DD} gives
 \begin{align}
\label{DD-1}
   &\partial_t\Lambda^{-1}(\rho^*-1)-p'(1)\Delta \Lambda^{-1}(\rho^*-1)+\Lambda^{-1}(\rho^*-1)\nonumber\\
   &\quad=\Lambda^{-1}\dive\big( (p'(\rho^*)-p'(1)) \nabla\rho^*+                           (\rho^*-1)\Lambda^{-2}\nabla\rho^*\big).
   \end{align}
Performing $L^2$ energy estimates on \eqref{DD-1} and using Young's inequality as well as the fact that $\Lambda^{-1}\dive$ is a bounded operator from $L^2$ to $L^2$, we have
\begin{align}\label{6.1}
&\frac{1}{2}\frac{d}{dt}\|\Lambda^{-1}(\rho^*-1)\|_{L^2}^2+\frac{1}{2} p'(1)\|\rho^*-1\|_{L^2}^2+\|\Lambda^{-1}(\rho^*-1)\|_{L^2}^2\nonumber\\
\leq & C\|(p'(\rho^*)-p'(1)) \nabla\rho^*\|_{L^2}^2+C\|(\rho^*-1)\Lambda^{-2}\nabla\rho^*\|_{L^2}^2\nonumber\\
\leq & C \|\rho^*-1\|_{L^{\infty}}^2\big( \|\nabla \rho^*\|_{L^2}^2+\|\Lambda^{-2}\nabla\rho^*\|_{L^2}^2\big)\nonumber\\
\leq & C \delta_1^*  \|\nabla \rho^*\|_{L^2}^2+ C\|\rho^*-1\|_{H^{m-1}}^2\|\rho^*-1\|_{\dot{H}^{-1}}^2,
\end{align}
from which and  \eqref{est:DD} we obtain
        \begin{align}\label{rho-1}
&\quad \sup_{t\in\mathbb{R}^{+}} \|\Lambda^{-1}(\rho^*(t)-1)\|_{L^2}^2+\int_0^\infty \big(\|\rho^*(t)-1\|_{L^2}^2+\|\Lambda^{-1}\rho^*(t)-1\|_{L^2}^2\big)\, dt \nonumber\\
&\leq Ce^{\int_0^\infty\|\rho^*(t)-1\|_{H^{m-1}}^2\,dt}\Big(\|\Lambda^{-1}(\rho^*_0-1)\|_{L^2}^2+\int_0^{\infty}\|\nabla \rho^*(t)\|_{L^2}^2\,dt\Big)\nonumber\\
&\leq C\|\rho_0^*-1\|_{ \dot{H}^{-1}\cap H^m}^2.
\end{align}
Since $E^*=\Lambda^{-2}\nabla \rho^*$, \eqref{rho-1} directly implies
        \begin{equation}\label{E-1}
\begin{aligned}
&\sup_{t\in\mathbb{R}^{+}} \|E^*(t)\|_{L^2}^2+\int_0^\infty\|E^*(t)\|_{L^2}^2\, dt\leq C\|\rho^*-1\|_{ \dot{H}^{-1}\cap H^m}^2.
\end{aligned}
\end{equation}
 Furthermore, employing \eqref{est:DD}, \eqref{E-1} and $q^*=-\nabla p(\rho^*)-\rho^*E^*$, we arrive at
\begin{align}\label{q*L2}
&\sup_{t\in\mathbb{R}^{+}}\|q^*(t)\|_{L^2}^2+\int_{0}^{\infty} \|q^*(t)\|_{L^2}^2\,dt\nonumber\\
&\leq  C(1+\sup_{t\in \mathbb{R}^{+}}\|\rho^*(t)-1\|_{H^{m}}^2) \Big( \sup_{t\in \mathbb{R}^{+}}\big(\|\rho^*(t)-1\|_{H^1}^2+\|E^*(t)\|_{L^2}^2\big)\nonumber\\
&\quad +\int_{0}^{\infty}\big(\|\rho^*(t)-1\|_{H^1}^2+\|E^*(t)\|_{L^2}^2\big)\,dt\Big)\nonumber\\
&\leq C\|\rho_0^*-1\|_{ \dot{H}^{-1}\cap H^m}^2.
\end{align}
Combining the above estimates \eqref{rho-1}, \eqref{E-1} and \eqref{q*L2} yields \eqref{Hmin1}. \hfill $\Box$

\subsection{Uniform estimates for $(\rho_1,q_1, E_1, B_1)$} We have the following lemma.
\begin{lemma}
\label{lemma62}
The Cauchy problem of \eqref{rho1EM} with $\rho_{1}|_{t=0}=0$ admits a unique global solution $\rho_1\in C(\mathbb{R}^{+};H^{m})$. Let $(q_1, E_1, B_1)$ be determined by
$\eqref{F3.11a11}_2$, \eqref{E1EM} and \eqref{B1EM}. Then, it holds that
\begin{align}\label{rhoqEBq}
&\quad\sup_{t\in \mathbb{R}^+} \big(\|\rho_1(t)\|_{\dot{H}^{-1}\cap H^{m}}^2+\|q_1(t)\|_{H^{m-1}}^2+\|E_1(t)\|_{H^{m+1}}^2+\|B_1(t)\|_{H^{m-1}}^2\big)\nonumber\\
&+\int_{0}^{\infty} \big(\|\rho_1(t)\|_{\dot{H}^{-1}\cap H^{m+1}}^2+\|q_1(t)\|_{H^{m}}^2+\|E_1(t)\|_{H^{m+1}}^2+\|B_1(t)\|_{H^{m+1}}^2\big)\,dt\leq C_{*},
\end{align}
where $C_{*}$ is a positive constant depending on $\rho_0^*$.
\end{lemma}


\noindent {\bf Proof.} Note that \eqref{rho1EM} can be rewritten as
\begin{align}\label{rho1linear}
  & \partial_t\rho_1-p'(1)\Delta \rho_1+\rho_1\nonumber\\
  &\quad=\dive\big ( p'(\rho^*)-p'(1))\nabla\rho_1\big) +\dive \big ( q^*\times  B^e+\rho_1 E^*+(\rho^*-1)\Lambda^{-2}\nabla\rho_1),
\end{align}
with the initial datum $\rho_1|_{t=0}=0$. According to the standard theorem for parabolic systems, there is a unique solution $\rho_1\in C(\mathbb{R}^{+};H^m)$ to \eqref{rho1linear} with $\rho_{1}|_{t=0}=0$.  Then, similarly to \eqref{6.1}, we have
\begin{align}\label{rho1L2}
&\quad\sup_{t\in \mathbb{R}^+} \|\Lambda^{-1}\rho_1(t)\|_{L^2}^2+\int_{0}^{\infty}\big(\|\Lambda^{-1}\rho_1(t)\|_{L^2}^2+\|\rho_1(t)\|_{L^2}^2\big)\,dt\nonumber\\
&\leq C\int_{0}^{\infty} \Big( \|( p'(\rho^*)-p'(1))\nabla\rho_1\|_{L^2}^2+\|(\rho^*-1)\nabla\Lambda^{-2}\rho_1(t)\|_{L^2}^2 \nonumber \\
&\quad+\|\rho_1 E^*(t)\|_{L^2}^2+\|q^*(t)\times B^e\|_{L^2}^2\Big)\,dt\nonumber\\
&\leq C\int_{0}^{\infty} \Big(\|\rho^*(t)-1\|_{L^{\infty}}^2\|\nabla\rho_1(t)\|_{L^2}^2+\|\rho^*(t)-1\|_{L^3}^2\|\nabla\Lambda^{-2}\rho_1(t)\|_{L^6}^2\nonumber\\
&\quad+\|\rho_1(t)\|_{L^3}^2\| E^*(t)\|_{L^6}^2+\|q^*(t)\|_{L^2}^2 \Big)\,dt\nonumber\\
&\leq \frac{1}{4}\int_{0}^{\infty} \|\rho_1(t)\|_{H^1}^2\,dt\nonumber\\
&\quad+C\int_{0}^{\infty} \Big(\big(\|\rho^*(t)-1\|_{H^{m-1}}^2+\|\nabla E^*(t)\|^2\big)\|\rho_1(t)\|_{H^1}^2+\|q^*(t)\|_{L^2}^2\Big)\,dt,
\end{align}
where the embeddings
$H^1\hookrightarrow L^3$, $\dot{H}^1\hookrightarrow L^6$ and $H^{m-1}\hookrightarrow L^{\infty}$ have been used. Concerning the higher-order estimates, one deduces from \eqref{F1.13} and \eqref{F1.15} that
\begin{align}\label{rho1Hs}
&\quad\sup_{t\in \mathbb{R}^+} \|\nabla\rho_1(t)\|_{H^{m-1}}^2+\int_{0}^{\infty} \|\nabla \rho_1(t)\|_{H^{m}}^2\,dt\nonumber\\
&\leq C\int_{0}^{\infty} \Big(\|( p'(\rho^*)-p'(1))\nabla\rho_1(t)\|_{H^{m-1}}+\|(\rho^*-1)\nabla\Lambda^{-2}\rho_1(t)\|_{H^{m-1}}\nonumber\\
&\quad+\|\rho_1 E^*(t)\|_{H^{m-1}}+\|q^*(t)\times B^e\|_{H^{m-1}}\Big)\|\nabla\rho_1(t)\|_{H^{m}} \,dt\nonumber\\
&\leq \frac{1}{4}\int_{0}^{\infty} \|\nabla\rho_1(t)\|_{H^{m}}^2\,dt\nonumber\\
&\quad+C\int_{0}^{\infty} \Big(\big(\|\rho^*(t)-1\|_{H^{m-1}}^2+\| E^*(t)\|_{H^{m-1}}^2\big)\|\rho_1(t)\|_{H^{m-1}}^2+\|q^*(t)\|_{H^{m-1}}^2\Big)\,dt.
\end{align}
Adding \eqref{rho1L2} and \eqref{rho1Hs} together and then using Gr\"onwall's lemma, \eqref{est:DD} and \eqref{Hmin1}, we infer that
\begin{align}
&\quad\sup_{t\in \mathbb{R}^+} \|\rho_1(t)\|_{\dot{H}^{-1}\cap H^{m}}^2+\int_{0}^{\infty} \|\rho_1(t)\|_{\dot{H}^{-1}\cap H^{m+1}}^2\, dt\nonumber\\
&\leq Ce^{C\int_0^{\infty}(\|\rho^*(t)-1\|_{H^{m-1}}^2+\|\nabla E^*(t)\|_{H^{m-2}}^2)\,dt} \int_0^{\infty}\|q^*(t)\|_{H^{m-1}}^2\,dt\leq C_*.
\end{align}
This, combined with $E_1=\Lambda^{-2}\nabla \rho_1$, gives rise to
\begin{equation}
\begin{aligned}
&\sup_{t\in \mathbb{R}^+} \|E_1(t)\|_{ H^{m+1}}^2+\int_{0}^{\infty} \|E_1(t)\|_{ H^{m+2}}^2\,dt\leq C_*.
\end{aligned}
\end{equation}
For $B_1$, note that 
$$
q^*=-\nabla P(\rho^*)-\rho^* E^*=-\nabla ( P(\rho^*)+\Lambda^{-2}\rho^*)-(\rho^*-1)\Lambda^{-2}\nabla\rho^*,
$$
and
$$
B_1=-\Lambda^{-2}\nabla\times q^*=-\Lambda^{-2}\nabla\times \big( (\rho^*-1)\Lambda^{-2}\nabla \rho^*\big).
$$
As $\|\Lambda^{-2}\nabla\times\cdot\|_{L^2}\leq C\|\cdot\|_{\dot{H}^{-1}}\leq C\|\cdot\|_{L^{\frac{6}{5}}}$ due to
$L^{\frac{6}{5}}\hookrightarrow \dot{H}^{-1}$, one deduces from \eqref{est:DD} that
\begin{equation*}
\begin{aligned}
&\quad \sup_{t\in \mathbb{R}^+} \|B_1(t)\|_{L^2}^2+\int_0^{\infty}\|B_1(t)\|_{L^2}^2\,dt\\
&\leq C\sup_{t\in \mathbb{R}^+} \|(\rho^*-1)\Lambda^{-2}\nabla \rho^*(t)\|_{L^{\frac{6}{5}}}^2+C\int_0^{\infty}\|(\rho^*-1)\Lambda^{-2}\nabla \rho^*(t)\|_{L^{\frac{6}{5}}}^2\,dt\\
&\leq C\sup_{t\in \mathbb{R}^+}\big(\|\rho^*(t)-1\|_{L^3}^2\|\Lambda^{-2}\nabla \rho^*(t)\|_{L^2}^2\big)\\
&\quad+C\sup_{t\in \mathbb{R}^+}\|\rho^*(t)-1\|_{L^3}^2\int_0^\infty \|\Lambda^{-2}\nabla \rho^*(t)\|_{L^2}^2\,dt\\
&\leq C\sup_{t\in \mathbb{R}^+}\big(\|\rho^*(t)-1\|_{H^m}^2\|\rho^*(t)-1\|_{\dot{H^{-1}}}^2\big)\\
&\quad+C\sup_{t\in \mathbb{R}^+}\|\rho^*(t)-1\|_{H^m}^2\int_0^\infty \|\rho^*(t)-1\|_{\dot{H}^{-1}}^2\,dt\leq C_*.
\end{aligned}
\end{equation*}
Similarly, it follows that
\begin{equation*}
\begin{aligned}
&\quad \sup_{t\in \mathbb{R}^+} \|\nabla B_1(t)\|_{\dot{H}^{m}}^2+\int_0^{\infty}\|\nabla B_1(t)\|_{H^{m-1}}^2\,dt\\
&\leq C\sup_{t\in \mathbb{R}^+}\Big(\|\rho^*(t)-1\|_{H^m}^2\|\rho^*(t)-1\|_{H^m}^2\big)\\
&\quad+C\sup_{t\in \mathbb{R}^+}\|\rho^*(t)-1\|_{H^m}^2\int_0^\infty \|\rho^*(t)-1\|_{H^m}^2\,dt\leq C_*.
\end{aligned}
\end{equation*}
Noting that
\[  q_1=-\na\big(p'(\rho^*)\rho_1\big)-\big(\rho^*E_1+\rho_1E^*\big)-q^*\times  B^e,
\]
we obtain
\begin{equation*}
\begin{aligned}
 &\quad\sup_{t\in \mathbb{R}^+} \|q_1(t)\|_{H^{m-1}}\\
 &\leq C\sup_{t\in \mathbb{R}^+}\big(\|p'(\rho^*)\rho_1(t)\|_{H^{m}}+\|\rho^*E_1(t)\|_{H^{m-1}}+\|\rho_1E^*(t)\|_{H^{m-1}}+\|q^*(t)\|_{H^{m-1}}\big)\\
&\leq C \sup_{t\in \mathbb{R}^+}\Big( \big(1+\|\rho^*(t)-1\|_{H^{m}}+\|\rho_1(t)\|_{H^{m}}\big) \big(\|\rho_1(t)\|_{H^{m-1}}+\|E_1(t)\|_{H^{m-1}}+\|E^*(t)\|_{H^{m-1}}\big)\\
&\quad+\|q^*(t)\|_{H^{m-1}} \Big)\leq C_*,
\end{aligned}
\end{equation*}
and
\begin{equation*}
\begin{aligned}
&\quad\int_{0}^{\infty}\|q_1(t)\|_{H^{m}}^2\,dt\\
&\leq C\int_{0}^{\infty}\Big( (1+\|\rho^*(t)-1\|_{H^{m}}^2+\|\rho_1(t)\|_{H^{m}}^2) \big(\|\rho_1(t)\|_{H^m}^2+\|E_1(t)\|_{H^m}^2+\|E^*(t)\|_{H^{m}}^2\big)\\
&\quad+\|q^*(t)\|_{H^{m}}^2 \Big)\,dt\leq C_*.
\end{aligned}
\end{equation*}
Combining these estimates, we end up with \eqref{rhoqEBq} and finish the proof of Lemma \ref{lemma62}. \hfill $\Box$

\vspace{3mm}
The convergence estimates of the remainder terms $J^\varepsilon_i$ ($i=1,2,3,4,5$) are
as follows.
\begin{lemma}\label{lemma63}
Let $J_i$ {\rm(}$i=1,2,3,4,5${\rm)} be given by \eqref{J12345}. Then, we have
\begin{align}
&\int_{0}^{\infty}\|J^\varepsilon_1(t)\|_{H^{m-1}}^2\,dt\leq C,\quad\quad \int_{0}^{\infty}\|J^\varepsilon_1(t)\|_{H^{m-1}}\,dt\leq C\varepsilon, \label{J1}\\
&\int_{0}^{\infty}\big(\|J^\varepsilon_2(t)\|_{H^{m-1}}^2+\|J^\varepsilon_3(t)\|_{H^{m}}^2+\|J_5^\varepsilon(t)\|_{H^{m}}^2\big)\,dt\leq C\varepsilon^2,\label{J3} \\
&\int_{0}^{\infty}\|J^\varepsilon_4(t)\|_{H^{m-2}}^2\,dt\leq C\varepsilon^4,\label{J4} \\
&\int_{0}^{\infty}\|J^\varepsilon_6(t)\|_{H^{m}}^2\,dt\leq C\varepsilon^2,\quad \quad \int_{0}^{\infty}\|J^\varepsilon_6(t)\|_{H^{m}}\,dt\leq C\varepsilon\label{J5J6},
\end{align}
where $C>0$ is a constant independent of $\varepsilon$.
\end{lemma}

\noindent {\bf Proof.} The estimates concerning $J_1^\varepsilon$, $J_2^\varepsilon$ and $J_3^\varepsilon$ are similar to those of Lemma \ref{lemmareminderEuler11}, so we omit the details. Recall that
\begin{align*}
&q^*=-\nabla p(\rho^*)-\rho^*E^*,\qquad\qquad\qquad \qquad\qquad\quad \,\,\,\, \, E^*=\Lambda^{-2}\nabla \rho^*,\\
&q_1=-\nabla \big(p'(\rho^*)\rho_1)-\rho^*E_1-\rho_1E^*-q^*\times B^e,\qquad E_1=\Lambda^{-2}\nabla\rho_1.
\end{align*}
Thus, in accordance with \eqref{F3.3} and \eqref{F3.11a11}, we have
$$
 \partial_t q^*=\partial_t\big(-\nabla p(\rho^*)-\rho^*E^*\big)=\nabla( p'(\rho^*)\dive q^*)+\dive q_1 E^*+\rho^*\Lambda^{-2}\nabla\dive q_1,
$$
and
 \begin{align*}
 \partial_t q_1&=\partial_t\big(-\nabla(p'(\rho^*)\rho_1)-\rho^*E_1-\rho_1 E^*-q^*\times B^e\big)\\
 &= \nabla\big(p''(\rho^*)\rho^*\dive q_1+ p'(\rho^*)\dive q_1\big)+\dive q^* E_1+\rho^*\Lambda^{-2}\nabla\dive q_1\\
 &\quad+\dive q_1 E^*+\rho_1\Lambda^{-2}\nabla\dive q^*-\partial_t q^*\times B^e.
 \end{align*}
Consequently, using \eqref{est:DD} and \eqref{rhoqEBq}, we have
\begin{align}
\int_0^{\infty}\|J_4^\varepsilon(t)\|_{H^{m-2}}^2\,dt&\leq C\varepsilon^4\sup_{t\in\mathbb{R}^{+}}\|\rho_1(t)\|_{H^m}^2\int_0^{\infty}\big(\|q^*(t)\|_{H^{m}}^2+\|q_1(t)\|_{H^m}^2\big)\,dt\leq C\varepsilon^4,
\end{align}
and
\begin{align}
\int_0^{\infty}\|J_4^\varepsilon(t)\|_{H^{m-2}}^2\,dt&\leq C\varepsilon^4\int_0^{\infty}\big(\|q^*(t)\|_{H^{m}}^2+\|q_1(t)\|_{H^m}^2\big)\,dt\leq C\varepsilon^4.
\end{align}

As
\[   J_5^\varepsilon=\varepsilon(\partial_t E_1-q_1)=-\varepsilon(\Lambda^{-2}\nabla\dive q_1+q_1), \]
one verifies that
\begin{align}
\int_0^{\infty}\|J_5^\varepsilon(t)\|_{H^{m}}^2\,dt&\leq C\varepsilon^2\int_0^{\infty}\|q_1(t)\|_{H^m}^2\,dt\leq C\varepsilon^2.
\end{align}
Finally,
\[  J_6^\varepsilon=\varepsilon \partial_tB_1=-\varepsilon \Lambda^{-2}\nabla\times \partial_t q^*=\varepsilon \Lambda^{-2}\nabla\times\big(\dive q_1 E^*+(\rho^*-1)\Lambda^{-2}\nabla\dive q_1\big).  \]
It holds by \eqref{Hmin1}, \eqref{rhoqEBq} and $L^{\frac{6}{5}}\hookrightarrow \dot{H}^{-1}$ that
\begin{align*}
&\quad\int_0^{\infty}\|J_6^\varepsilon(t)\|_{L^2}^2\,dt\nonumber\\
&\leq C\varepsilon^2 \int_0^{\infty}\big(\|\dive q_1 E^*(t)\|_{\dot{H}^{-1}}^2+\|(\rho^*-1)\Lambda^{-2}\nabla\dive q_1(t)\|_{\dot{H}^{-1}}^2\big)\,dt\\
&\leq C\varepsilon^2 \int_0^{\infty}\big(\|\dive q_1 E^*(t)\|_{L^{\frac{6}{5}}}^2+\|(\rho^*-1)\Lambda^{-2}\nabla\dive q_1(t)\|_{L^{\frac{6}{5}}}^2\big)\,dt\\
&\leq C\varepsilon^2\int_0^{\infty}\big( \|\dive q_1(t)\|_{L^3}^2\|E^*(t)\|_{L^2}^2+\| \rho^*(t)-1\|_{L^3}^2\|\Lambda^{-2}\nabla\dive q_1(t)\|_{L^2}^2\big)\,dt\\
&\leq C\varepsilon^2\int_0^{\infty}\|q_1(t)\|_{H^1}^2\,dt\leq C\varepsilon^2.
\end{align*}
For the higher-order estimates, one also deduces from \eqref{F1.1311} and $H^{m-1}\hookrightarrow L^{\infty}$ that
\begin{align*}
&\quad\int_0^{\infty}\|J_6^\varepsilon(t)\|_{\dot{H}^{m+1}}^2\,dt\nonumber\\
&\leq C\varepsilon^2 \int_0^{\infty}\big(\|\dive q_1 E^*(t)\|_{\dot{H}^{m}}^2+\|(\rho^*-1)\Lambda^{-2}\nabla\dive q_1(t)\|_{\dot{H}^{m}}^2\big)\,dt\\
&\leq C\varepsilon^2 \int_0^{\infty}\big(\|\dive q_1(t)\|_{L^{\infty}}^2 \|E^*(t)\|_{\dot{H}^{m}}^2+\|\dive q_1(t)\|_{\dot{H}^{m}}^2 \|E^*(t)\|_{L^{\infty}}^2\\
&\quad+\|\rho^*(t)-1\|_{L^{\infty}}^2\|\Lambda^{-2}\nabla\dive q_1(t)\|_{\dot{H}^{m}}^2+\|\rho^*(t)-1\|_{\dot{H}^{m}}^2\|\Lambda^{-2}\nabla\dive q_1(t)\|_{L^{\infty}}^2 \big)\,dt\\
&\leq C\varepsilon^2\int_0^{\infty}\big( \|E^*(t)\|_{H^m}^2+\| \rho^*(t)-1\|_{H^m}^2\big)\|q_1(t)\|_{H^m}^2\,dt\\
&\leq C\varepsilon^2\int_0^{\infty}\|q_1(t)\|_{H^m}^2\,dt\leq C\varepsilon^2.
\end{align*}
Consequently, we have the first estimate in \eqref{J5J6}. Similarly, $\|J_6^\varepsilon(t)\|_{H^m}$ has the $L^1$ time integrability:
\begin{align*}
&\quad\int_0^{\infty}\|J_6^\varepsilon(t)\|_{H^m}\,dt\nonumber\\
&\leq C\varepsilon \int_0^{\infty}\big(\|\dive q_1 E^*(t)\|_{\dot{H}^{-1}\cap \dot{H}^{m-1}}+\|(\rho^*-1)\Lambda^{-2}\nabla\dive q_1(t)\|_{\dot{H}^{-1}\cap \dot{H}^{m-1}}\big)\,dt\\
&\leq C\varepsilon \int_0^{\infty}\big( \|E^*(t)\|_{H^m}+\| \rho^*(t)-1\|_{H^m}\big) \| q_1(t)\|_{H^m}\big)\,dt\\
&\leq C\varepsilon \Big(\int_0^{\infty}\big( \|E^*(t)\|_{H^m}^2+\| \rho^*(t)-1\|_{H^m}^2\big)\,dt\Big)^{\frac{1}{2}}\Big(\int_0^{\infty} \| q_1(t)\|_{H^m}^2\,dt\Big)^{\frac{1}{2}}\leq C\varepsilon.
\end{align*}
This finishes the proof of Lemma \ref{lemma63}.
\hfill $\Box$

\subsection{Low-order error estimates}

We prove the error estimates in \texorpdfstring{$L^2$}{L2} for the solutions $(\rho^\varepsilon,q^\varepsilon,E^\varepsilon,B^\varepsilon)$ and $(\rho^*,q^*,E^*,B^*)$ associated with the systems \eqref{EMeps}
and \eqref{F3.3}, respectively.

\begin{lemma}
\label{L3.1}
Let $(\tilde{\rho}^\varepsilon,\tilde{q}^\varepsilon, \tilde{E}^\varepsilon,\tilde{B}^\varepsilon)=(\rho^\varepsilon-\rho^*,q^\varepsilon-q^*,E^{\varepsilon}-E^*,B^\varepsilon- B^e)$.  For all $t\ge 0$, it holds
\begin{align}
\label{F3.14}
\quad &\sup_{t\in\mathbb{R}^{+}}\big(\|\tilde{\rho}^\varepsilon(t)\|_{L^2}^2
+\varepsilon^2\|\tilde{q}^\varepsilon(t)\|_{L^2}^2+\|\tilde{E}^\varepsilon(t)\|_{L^2}^2+\|\tilde{B}^\varepsilon(t)\|_{L^2}^2\big)
\nonumber\\
&\quad\quad+\int_{0}^{\infty}\big(\|\tilde{\rho}^\varepsilon(t)\|_{H^1}^2+\|\tilde{q}^\varepsilon(t)\|_{L^2}^2+\|\tilde{E}^\varepsilon(t)\|_{L^2}^2\big)\,dt\nonumber\\
 &\le C\big(\|\rho_0^{\varepsilon}-\rho_0^*\|_{L^2}^2+\|E_0^{\varepsilon}-E_0^*\|_{L^2}^2+\|B_0^{\varepsilon}-B^e\|_{L^2}^2\big)+C \varepsilon^2,
\end{align}
where $C>0$ is a uniform constant.
\end{lemma}

\noindent {\bf Proof.}
Let $\varepsilon_1\in(0,1)$ be a suitably small constant. In the case $\varepsilon\geq \varepsilon_1$, one easily gets \eqref{F3.14} using the lower bound $\varepsilon_1$ and the  uniform estimates for $(\tilde{\rho}^\varepsilon,\tilde{q}^\varepsilon,\tilde{E}^\varepsilon,\tilde{B}^\varepsilon)$. In what follows, we consider $\varepsilon\leq \varepsilon_1<<1$.

 As emphasized in Subsection \ref{subsection51}, our idea is to derive the estimates for $(\tilde{\rho}_a^\varepsilon,\tilde{q}_a^\varepsilon, \tilde{E}_a^\varepsilon,\tilde{B}_a^\varepsilon)$ defined in \eqref{tildeerrora} and then recover the desired estimates for $(\tilde{\rho}^\varepsilon,\tilde{q}^\varepsilon, \tilde{E}^\varepsilon,\tilde{B}^\varepsilon)$ due to the bounds for the profiles $\rho_1$, $q_1$, $E_1$ and $B_1$ obtained in Lemma \ref{lemma62}.

Recall that the equations $\eqref{F3.11ca}_1$-$\eqref{F3.11ca}_2$ are
\begin{equation}\label{F3.11cb}
\left\{
\begin{aligned}
 &\partial_t\tilde{\rho}^\varepsilon_a+\dive \tilde{q}^\varepsilon_a=J_1^\varepsilon,\\
&\varepsilon^2\partial_t\tilde{q}^\varepsilon_a+\tilde{q}^\varepsilon_a+p'(1)\nabla \tilde{\rho}^\varepsilon_a +\tilde{E}^\varepsilon_a=J_0^\varepsilon+J_2^\varepsilon+J_3^\varepsilon+J^\varepsilon_4.
\end{aligned}
\right.
\end{equation}
By taking the $L^2$-inner products of $\eqref{F3.11cb}_1$ with $p'(1)\tilde{\rho}^\varepsilon_a$ and of $\eqref{F3.11cb}_2$ with $\tilde{\rho}^\varepsilon_a$, respectively, we have
\begin{align}
  &\frac{1}{2}\frac{d}{dt}\big( p'(1)\|\tilde{\rho}_a^\varepsilon\|_{L^2}^2
+\varepsilon^2\| \tilde{q}_a^\varepsilon\|_{L^2}^2\big)
+\|\tilde{q}_a^\varepsilon\|_{L^2}^2+\langle \tilde{E}_a^\varepsilon,\tilde{q}^\varepsilon_a\rangle
\nonumber\\
=& p'(1)\langle J^\varepsilon_1,\tilde{\rho}_a^\varepsilon\rangle+\langle J^\varepsilon_0+J^\varepsilon_2+J^\varepsilon_3+J^\varepsilon_4,\tilde{q}^\varepsilon_a\rangle.\label{Fmmmmm}
\end{align}
To cancel the last term on the left-hand side of \eqref{Fmmmmm}, we shall make use of the structure of $\eqref{F3.11ca}$. In fact, from $\eqref{F3.11ca}_3$-$\eqref{F3.11ca}_4$ one has
\begin{equation}\label{Fmmmmm1}
\begin{aligned}
\frac{1}{2}\frac{d}{dt}\Big(\|\tilde{E}^\varepsilon_a\|_{L^2}^2+\|\tilde{B}^\varepsilon_a\|_{L^2}^2\Big)-\langle \tilde{q}_a^\varepsilon,\tilde{E}^\varepsilon_a\rangle =\langle J^\varepsilon_5, \tilde{E}^\varepsilon_a\rangle+\langle J_6^\varepsilon, \tilde{B}^\varepsilon_a\rangle.
\end{aligned}
\end{equation}
Adding \eqref{Fmmmmm} with \eqref{Fmmmmm1}, we obtain
\begin{align}\label{5.26}
 &\frac{1}{2}\frac{d}{dt}\Big(p'(1)\|\tilde{\rho}_a^\varepsilon\|_{L^2}^2
+\varepsilon^2\|\tilde{q}_a^\varepsilon\|_{L^2}^2+\|\tilde{E}^\varepsilon_a\|_{L^2}^2+\|\tilde{B}^\varepsilon_a\|_{L^2}^2 \Big)
+\|\tilde{q}_a^\varepsilon\|_{L^2}^2
\nonumber\\
=& \,p'(1)\langle J^\varepsilon_1,\tilde{\rho}_a^\varepsilon\rangle
+\langle J^\varepsilon_0+J^\varepsilon_2+J^\varepsilon_3+J^\varepsilon_4,\tilde{q}^\varepsilon_a\rangle+\langle J^\varepsilon_5, \tilde{E}^\varepsilon_a\rangle+\langle J_6^\varepsilon, \tilde{B}^\varepsilon_a\rangle\nonumber\\
\leq & \frac{1}{2}\|\tilde{q}_a^\varepsilon\|_{L^2}^2+p'(1)\|J^\varepsilon_1\|_{L^2}\|\tilde{\rho}_a^\varepsilon\|_{L^2}+C\nu^{-1}\|J^\varepsilon_5\|_{L^2}^2+\nu\|\tilde{E}^\varepsilon_a\|_{L^2}^2+\|J^\varepsilon_6\|_{L^2}\|\tilde{B}^\varepsilon_a\|_{L^2}\nonumber\\
&\quad+C(\|J_0^\varepsilon\|_{L^2}^2+\|J^\varepsilon_2\|_{L^2}^2+\|J^\varepsilon_3\|_{L^2}^2+\|J^\varepsilon_4\|_{L^2}^2).
\end{align}
Here, $\nu>0$ is a small constant to be chosen later. Similarly to that in \eqref{caseall}, using the regularity estimate for $(\rho^\varepsilon,E^\varepsilon,B^\varepsilon)$, $(\rho^*,E^*,B^*)$ and $(\rho_1,E_1,B_1)$, we have
\begin{align}\label{ffd}
&\quad \int_0^t \|J_0^\varepsilon(t')\|_{L^2}^2\,dt'\nonumber \\
&\leq C\sup_{t'\in [0,t]} \big(\|\nabla \rho^\varepsilon_a(t')\|_{L^{\infty}}^2+\|\rho^*(t')-1\|_{L^{\infty}}^2+\|\rho^\varepsilon(t')-1\|_{L^{\infty}}^2+\|E^\varepsilon(t')\|_{L^{\infty}}^2\big)\nonumber\\
&\quad\times \int_0^t\big(\|\tilde{\rho}^\varepsilon_a(t')\|_{H^1}^2+\|\tilde{E}^\varepsilon_a(t')\|_{L^2}^2\big)\,dt'\nonumber\\
&\leq C(\delta_1+\delta_1^*+\varepsilon)\int_0^t\big(\|\tilde{\rho}^\varepsilon_a(t')\|_{H^1}^2+\|\tilde{E}^\varepsilon_a(t')\|_{L^2}^2\big)\,dt'.
\end{align}
Thus, integrating \eqref{5.26} over $[0,t]$ and using \eqref{J1}-\eqref{J5J6} and \eqref{ffd}, we arrive at
\begin{equation*}
\begin{aligned}
&\quad \sup_{t\in\mathbb{R}^+}\big(\|\tilde{\rho}_a^\varepsilon(t)\|_{L^2}^2
+\varepsilon^2\| \tilde{q}_a^\varepsilon(t)\|_{L^2}^2+\|\tilde{E}^\varepsilon_a(t)\|_{L^2}^2+\|\tilde{B}^\varepsilon_a(t)\|_{L^2}^2\big)+\int_0^{\infty}\|\tilde{q}_a^\varepsilon(t)\|_{L^2}^2\,dt\\
&\leq C\big(\|\tilde{\rho}_a^\varepsilon|_{t=0}\|_{L^2}^2
+\varepsilon^2\| \tilde{q}_a^\varepsilon|_{t=0}\|_{L^2}^2+\|\tilde{E}_a^\varepsilon|_{t=0}\|_{L^2}^2+\|\tilde{B}_a^\varepsilon|_{t=0}\|_{L^2}^2\big)\\
&\quad+C(\delta_1+\delta_1^*+\varepsilon+\nu)\int_0^t\big(\|\tilde{\rho}_a(t')\|_{H^1}^2+\|\tilde{E}_a(t')\|_{L^2}^2\big)\,dt'\\
&\quad+C\int_0^{\infty}\big(\|J_2^\varepsilon(t)\|_{L^2}^2+\|J_3^\varepsilon(t)\|_{L^2}^2+\|J^\varepsilon_4(t)\|_{L^2}^2+\nu^{-1}\|J^\varepsilon_5(t)\|_{L^2}^2\big)\,dt\\
&\quad\quad+C\int_0^{\infty}\big(\|J^\varepsilon_1(t)\|_{L^2}+\|J_6^\varepsilon(t)\|_{L^2}\big)\,dt\,\sup_{t\in\mathbb{R}^+}\big(\|\tilde{\rho}_a^\varepsilon(t)\|_{L^2}
+\|\tilde{B}^\varepsilon_a(t)\|_{L^2}\big)\\
&\leq C\big(\|\tilde{\rho}_a^\varepsilon|_{t=0}\|_{L^2}^2
+\varepsilon^2\| \tilde{q}^\varepsilon|_{t=0}\|_{L^2}^2+\|\tilde{E}_a^\varepsilon|_{t=0}\|_{L^2}^2+\|\tilde{B}_a^\varepsilon|_{t=0}\|_{L^2}^2\big)\\
&\quad+C(\delta_1+\delta_1^*+\varepsilon)\int_0^\infty\big(\|\tilde{\rho}_a(t)\|_{H^1}^2+\|\tilde{E}_a(t)\|_{L^2}^2\big)\,dt\\
&\quad+\Big(C\delta_1+C\delta_1^*+C\varepsilon+\frac{1}{2}\Big)\sup_{t\in\mathbb{R}^+}\big(\|\tilde{\rho}_a^\varepsilon(t)\|_{L^2}^2
+\|\tilde{B}^\varepsilon_a(t)\|_{L^2}^2\big)+C(1+\nu^{-1})\varepsilon^2.
\end{aligned}
\end{equation*}
Note that
\begin{align*}
&\big(\tilde{\rho}^\varepsilon_a,\tilde{q}^\varepsilon_a, \tilde{E}^\varepsilon_a,\tilde{B}^\varepsilon_a\big)|_{t=0}\\
&\quad=\big(\rho_0^\varepsilon-\rho_0^*,{\nabla p(\rho_0^*)}, E_0^\varepsilon-E_0^*,B_0^\varepsilon-B^e\big)+\varepsilon\big(0,-q^*|_{t=0}\times  B^e,0, -\Lambda^{-2}\nabla\times q^*|_{t=0}\big).
\end{align*}
Using
\[ q^*=-\nabla p(\rho^*)-\rho^* E^*,\qquad  E^*=\Lambda^{-2}\nabla \rho^*=\nabla\Lambda^{-2}(\rho^*-1),  \]
and the regularity condition $\rho^*_0-1\in \dot{H}^{-1}\cap H^m$ , we know
\[  q^*|_{t=0}=-\nabla p(\rho_0^*)+\rho_0^* \nabla \Lambda^{-2}(\rho_0^*-1)\in L^2.  \]
In addition, we have
\[  \nabla\times q^*|_{t=0}=-\nabla\times \big((\rho_0^*-1) \nabla \Lambda^{-2}(\rho_0^*-1)\big)  \]
and
\begin{align}\label{data1order}
\|\Lambda^{-2}\nabla\times q^*|_{t=0})\|_{L^2}&\leq C\|(\rho_0^*-1) \nabla \Lambda^{-2}(\rho_0^*-1)\|_{\dot{H}^{-1}}\nonumber\\
&\leq C\|(\rho_0^*-1) \nabla \Lambda^{-2}(\rho_0^*-1)\|_{L^{\frac{6}{5}}}\nonumber\\
&\leq C \|\rho_0^*-1\|_{L^3}\|\nabla \Lambda^{-2}(\rho_0^*-1)\|_{L^2}\nonumber\\
&\leq C\|\rho_0^*-1\|_{H^m}\|\rho_0^*-1\|_{\dot{H}^{-1}},
\end{align}
due to $L^{\frac{6}{5}}\hookrightarrow \dot{H}^{-1}$. Consequently, it holds
\begin{align*}
&\|\tilde{\rho}^\varepsilon_a|_{t=0}\|_{L^2}^2+\|\tilde{q}^\varepsilon_a|_{t=0}\|_{L^2}^2+ \|\tilde{E}^\varepsilon_a|_{t=0}\|_{L^2}^2+\|\tilde{B}^\varepsilon_a|_{t=0}\|_{L^2}^2\\
&\quad\leq C\big(\|\rho_0^\varepsilon-\rho_0^*\|_{L^2}^2+\|E_0^\varepsilon-E_0^*\|_{L^2}^2+\|B_0^\varepsilon-B^e\|_{L^2}^2\big)+C\varepsilon^2.
\end{align*}
As $\delta_1$, $\delta_1^*$ and $\varepsilon$ are sufficiently small, the following estimate holds:
\begin{align}\label{5.27}
&\quad \sup_{t\in\mathbb{R}^+}\big(\|\tilde{\rho}_a^\varepsilon(t)\|_{L^2}^2
+\varepsilon^2\| \tilde{q}_a^\varepsilon(t)\|_{L^2}^2+\|\tilde{E}^\varepsilon_a(t)\|_{L^2}^2+\|\tilde{B}^\varepsilon_a(t)\|_{L^2}^2\big)+\int_0^{\infty}\|\tilde{q}^\varepsilon_a(t)\|_{L^2}^2\,dt\nonumber\\
&\leq C\big(\|\rho_0^\varepsilon-\rho_0^*\|_{L^2}^2+\|E_0^\varepsilon-E_0^*\|_{L^2}^2+\|B_0^\varepsilon-B^e\|_{L^2}^2\big)+C(1+\nu^{-1})\varepsilon^2\nonumber\\
&\quad+C(\delta_1+\delta_1^*+\varepsilon+\nu)\int_0^{\infty}\big(\|\tilde{\rho}^\varepsilon_a(t)\|_{H^1}^2+\|\tilde{E}^\varepsilon_a(t)\|_{L^2}^2\big)\,dt.
\end{align}

Next, the dissipation estimates for $\tilde{\rho}_a^\varepsilon$ are established as follows. Taking the inner product of $\eqref{F3.11ca}_2$ by $\nabla \tilde{\rho}^\varepsilon_a$ and making use of $\eqref{F3.11ca}_1$ and $\dive \tilde{E}^\varepsilon_a=-\tilde{\rho}^\varepsilon_a$, we have
\begin{align}
&\quad p'(1)\|\nabla\tilde{\rho}_a^\varepsilon\|_{L^2}^2+\|\tilde{\rho}^\varepsilon_a\|_{L^2}^2\nonumber\\
 &=-\varepsilon^2\langle\partial_{t} \tilde{q}_a^\varepsilon, \nabla \tilde{\rho}_a^\varepsilon\rangle -\langle \tilde{q}_a^\varepsilon, \nabla\tilde{\rho}_a^\varepsilon\rangle+\langle J^\varepsilon_0+J^\varepsilon_2+J^\varepsilon_3+J^\varepsilon_4,\nabla\tilde{\rho}^\varepsilon_a\rangle\nonumber \\
&=-\varepsilon^2\frac{d}{dt}\langle \tilde{q}_a^\varepsilon, \nabla \tilde{\rho}_a^\varepsilon\rangle -\langle \tilde{q}_a^\varepsilon, \nabla\tilde{\rho}_a^\varepsilon\rangle+\varepsilon^2\|\dive\tilde{q}^\varepsilon_a\|_{L^2}^2\nonumber \\
&\quad+\langle J^\varepsilon_0+J^\varepsilon_2+J^\varepsilon_3+J^\varepsilon_4,\nabla\tilde{\rho}^\varepsilon_a\rangle-\varepsilon^2\langle \dive \tilde{q}^\varepsilon_a,J^\varepsilon_1\rangle\nonumber\\
&\leq -\varepsilon^2\frac{d}{dt}\langle \tilde{q}_a^\varepsilon, \nabla \tilde{\rho}_a^\varepsilon\rangle +C\|\tilde{q}^\varepsilon\|_{L^2}^2+C\varepsilon^2 \|\dive\tilde{q}^\varepsilon_a\|_{L^2}^2 +C\|J_0^\varepsilon\|_{L^2}^2+C\|J_2^\varepsilon\|_{L^2}^2+C\|J_3^\varepsilon\|_{L^2}^2\nonumber\\
&\quad+C\|J_4^\varepsilon\|_{L^2}^2+C\varepsilon^2 \|J_1^\varepsilon\|_{L^2}^2+\frac{p'(1)}{2}\|\nabla \tilde{\rho}^\varepsilon_a\|_{L^2}^2.\label{5.280}
\end{align}
Thus, integrating \eqref{5.280} in time and recalling \eqref{ffd}, we discover that
\begin{align}
&\quad \int_0^{t}\Big(\frac{p'(1)}{2} \|\nabla\tilde{\rho}_a^\varepsilon(t')\|_{L^2}^2+\|\tilde{\rho}_a^\varepsilon(t')\|_{L^2}^2\Big)\,dt'\nonumber\\
&\leq -\varepsilon^2\langle \tilde{q}_a^\varepsilon, \nabla \tilde{\rho}_a^\varepsilon\rangle \Big|^{t}_{0}+C(\delta_1+\delta_1^*+\varepsilon)\int_0^t\big(\|\tilde{\rho}_a(t')\|_{H^1}^2+\|\tilde{E}_a(t')\|_{L^2}^2\big)\,dt' \nonumber\\
&\quad+C\int_0^{t}\|\tilde{q}_a^\varepsilon(t')\|_{L^2}^2\,dt'+C \varepsilon^2 \int_0^t\|\dive \tilde{q}^\varepsilon_a(t')\|_{L^2}^2 \,dt'\nonumber\\
&\quad+C\int_0^t\big(\varepsilon^2\|J_1^\varepsilon(t')\|_{L^2}^2+\|J_2^\varepsilon(t')\|_{L^2}^2+\|J_3^\varepsilon(t')\|_{L^2}^2+\|J^\varepsilon_4(t')\|_{L^2}^2\big)\,dt'.\label{5.28}
\end{align}
As the bounds for $\rho^\varepsilon, \rho^*, \rho_1, q^\varepsilon, q^*$ and $q_1$ have been obtained, it holds
\begin{equation}
\begin{aligned}
\sup_{t\in\mathbb{R}^+}\big(\|\tilde{\rho}_a^\varepsilon(t)\|_{H^{m}}^2+\|\tilde{E}_a^\varepsilon(t)\|_{H^{m}}^2\big)+\int_0^{\infty}\big(\|\tilde{q}_a^\varepsilon(t)\|_{H^{m}}^2+\|\nabla\tilde{B}^\varepsilon_a(t)\|_{H^{m-1}}^2\big)\,dt\leq C.\label{5.29r}
\end{aligned}
\end{equation}
By \eqref{5.2}, \eqref{5.29r} and $\rho_1|_{t=0}=0$, we have
\begin{equation}
\begin{aligned}\label{5.30}
-\varepsilon^2\langle \tilde{q}_a^\varepsilon, \nabla \tilde{\rho}_a^\varepsilon\rangle \Big|^{t}_{0}&\leq \varepsilon^2\|\tilde{q}_a^\varepsilon\|_{L^2}\| \nabla\tilde{\rho}_a^\varepsilon\|_{L^2}+\varepsilon^3 \| q^*|_{t=0}\|_{L^2}\|\rho_0^\varepsilon-\rho_0^*\|_{L^2}\leq C\varepsilon^2\|\tilde{q}_a^\varepsilon\|_{L^2}^2+C\varepsilon^2.
\end{aligned}
\end{equation}
Here we used
\begin{align*}
\|q^*|_{t=0}\|_{L^2}&\leq \|\nabla p(\rho^*)\|_{L^2}+\|\rho_0^* E_0^*\|_{L^2} \nonumber\\
&\leq C\|\rho^*-1\|_{H^1}+C(1+\|\rho_0^*-1\|_{L^{\infty}})\|E_0^*\|_{L^2}\leq C. \nonumber
\end{align*}
Putting \eqref{J1}-\eqref{J3} and \eqref{5.30} into \eqref{5.28}, we derive
\begin{align}\label{5.32}
&\quad \int_0^{\infty}\big( \|\nabla\tilde{\rho}_a^\varepsilon(t)\|_{L^2}^2+\|\tilde{\rho}_a^\varepsilon(t)\|_{L^2}^2\big)\,dt\nonumber\\
&\leq C\sup_{t\in\mathbb{R}^+}\varepsilon^2\|\tilde{q}_a^\varepsilon(t)\|_{L^2}^2+C\int_0^{\infty}\|\tilde{q}^\varepsilon_a(t)\|_{L^2}^2\,dt\nonumber\\
&\quad+C(\delta_1+\delta_1^*+\varepsilon)\int_0^\infty\big(\|\tilde{\rho}^\varepsilon_a(t)\|_{H^1}^2+\|\tilde{E}^\varepsilon_a(t)\|_{L^2}^2\big)\,dt+C\varepsilon^2.
\end{align}

Another key point is to establish the dissipation estimate for the error $\tilde{E}^\varepsilon_a$ required in \eqref{5.27}. Taking the inner product of $\eqref{F3.11ca}_2$ by $ \tilde{E}^\varepsilon_a$, making use of $\eqref{F3.11ca}_3$ and keeping in mind that $\dive \tilde{E}^\varepsilon_a=-\tilde{\rho}^\varepsilon_a$, we obtain
\begin{equation*}
\begin{aligned}
 & \quad p'(1)\| \dive \tilde{E}_a^\varepsilon\|_{L^2}^2+\| \tilde{E}^\varepsilon_a\|_{L^2}^2\\
 &=-\varepsilon^2\langle\partial_{t} \tilde{q}_a^\varepsilon, \tilde{E}_a^\varepsilon\rangle -\langle \tilde{q}_a^\varepsilon, \tilde{E}_a^\varepsilon\rangle+\langle J^\varepsilon_0+J^\varepsilon_2+J^\varepsilon_3+J^\varepsilon_4,\tilde{E}_a^\varepsilon\rangle \\
&=-\varepsilon^2\frac{d}{dt}\langle \tilde{q}_a^\varepsilon,  \tilde{E}_a^\varepsilon\rangle-\langle \tilde{q}_a^\varepsilon, \tilde{E}_a^\varepsilon\rangle+\langle J^\varepsilon_0+J^\varepsilon_2+J^\varepsilon_3+J^\varepsilon_4,\tilde{E}_a^\varepsilon\rangle\\
&\quad+\varepsilon\langle\tilde{q}^\varepsilon_a,\nabla\times\tilde{B}^\varepsilon_a\rangle+\varepsilon^2\langle \tilde{q}^\varepsilon_a,J_5^\varepsilon\rangle\\
&\leq -\varepsilon^2\frac{d}{dt}\langle \tilde{q}_a^\varepsilon,  \tilde{E}_a^\varepsilon\rangle+\frac{1}{2}\| \tilde{E}^\varepsilon_a\|_{L^2}^2+C\| \tilde{q}^\varepsilon_a\|_{L^2}^2+C\varepsilon^2\| \tilde{q}^\varepsilon_a\|_{L^2}^2+ C\varepsilon^2 \|\nabla\times \tilde{B}^\varepsilon_a\|_{L^2}^2\\
&\quad+C\|J_0^\varepsilon\|_{L^2}^2+C\|J_2^\varepsilon\|_{L^2}^2+C\|J_3^\varepsilon\|_{L^2}^2+C\|J_4^\varepsilon\|_{L^2}^2+C\varepsilon^2 \|J_5^\varepsilon\|_{L^2}^2.
\end{aligned}
\end{equation*}
Consequently, it follows that
\begin{equation*}
\begin{aligned}
&\quad\int_0^t\Big(p'(1)\|\dive \tilde{E}^\varepsilon_a(t')\|_{L^2}^2+\frac{1}{2}\|\tilde{E}^\varepsilon_a(t')\|_{L^2}^2\Big)\,dt'\\
&\leq -\varepsilon^2\langle \tilde{q}^\varepsilon_a,\tilde{E}^\varepsilon_a\rangle \Big|^t_0+C\int_0^t\|\tilde{q}^\varepsilon_a(t')\|_{L^2}^2\,dt'+C\int_0^t \|J_0^{\varepsilon}(t')\|_{L^2}^2\,dt'\\
&\quad+C\int_0^t\big(\varepsilon^2 \|\nabla\tilde{B}^\varepsilon_a(t')\|_{L^2}^2+\|J_2^\varepsilon(t')\|_{L^2}^2+\|J_3^\varepsilon(t')\|_{L^2}^2+\|J_4^\varepsilon(t')\|_{L^2}^2+\varepsilon^2\|J_5(t')\|_{L^2}^2\big)\,dt',
\end{aligned}
\end{equation*}
which, together with \eqref{J3}-\eqref{J5J6}, \eqref{ffd} and \eqref{5.29r}, leads to
\begin{equation}\label{5.35}
\begin{aligned}
\int_0^\infty \|\tilde{E}^\varepsilon_a(t')\|_{L^2}^2\,dt'&\leq C\sup_{t\in\mathbb{R}^+}\varepsilon^2\|\tilde{q}_a^\varepsilon(t)\|_{L^2}^2+C\int_0^{\infty}\|\tilde{q}^\varepsilon_a(t)\|_{L^2}^2\,dt\\
&\quad+C(\delta_1+\delta_1^*)\int_0^\infty\big(\|\tilde{\rho}^\varepsilon_a(t)\|_{H^1}^2+\|\tilde{E}^\varepsilon_a(t)\|_{L^2}^2\big)\,dt+C\varepsilon^2.
\end{aligned}
\end{equation}

Finally, we let $\eqref{5.27}+\eta \eqref{5.32}+\eta \eqref{5.35}$ with a uniform small constant $0<\eta<<1$ and then obtain
\begin{align*}
&\quad \sup_{t\in\mathbb{R}^{+}}\big(\|\tilde{\rho}_a^\varepsilon(t)\|_{L^2}^2
+\varepsilon^2\|\tilde{q}_a^\varepsilon(t)\|_{L^2}^2+\|\tilde{E}_a^\varepsilon(t)\|_{L^2}^2+\|\tilde{B}_a^\varepsilon(t)\|_{L^2}^2\big)
\nonumber\\
&\quad\quad+\int_{0}^{\infty}\big(\|\tilde{\rho}_a^\varepsilon(t)\|_{H^1}^2+\|\tilde{q}_a^\varepsilon(t)\|_{L^2}^2+\|\tilde{E}_a^\varepsilon(t)\|_{L^2}^2\big)\,dt\nonumber\\
 &\le C\big(\|\rho_0^{\varepsilon}-\rho_0^*\|_{L^2}^2+\|E_0^{\varepsilon}-E_0^*\|_{L^2}^2+\|B_0^{\varepsilon}-B^e\|_{L^2}^2\big)+C(1+\nu^{-1}) \varepsilon^2\\
 &\quad+C(\delta_1+\delta_1^*+\varepsilon+\nu)\int_0^{\infty}\big(\|\tilde{\rho}^\varepsilon_a(t)\|_{H^1}^2+\|\tilde{E}^\varepsilon_a(t)\|_{L^2}^2\big)\,dt.
\end{align*}
Choosing a suitably small constant $\nu>0$ and making use of the smallness for $\delta_1, \delta_1^*$ and $\varepsilon$, we conclude that
\begin{align}\label{widetildeL2EM}
&\quad \sup_{t\in\mathbb{R}^{+}}\big(\|\tilde{\rho}_a^\varepsilon(t)\|_{L^2}^2
+\varepsilon^2\|\tilde{q}_a^\varepsilon(t)\|_{L^2}^2+\|\tilde{E}_a^\varepsilon(t)\|_{L^2}^2+\|\tilde{B}_a^\varepsilon(t)\|_{L^2}^2\big)\nonumber
\\
&\quad\quad+\int_{0}^{\infty}\big(\|\tilde{\rho}_a^\varepsilon(t)\|_{H^1}^2+\|\tilde{q}_a^\varepsilon(t)\|_{L^2}^2+\|\tilde{E}_a^\varepsilon(t)\|_{L^2}^2\big)\,dt\nonumber\\
 &\le C\big(\|\rho_0^{\varepsilon}-\rho_0^*\|_{L^2}^2+\|E_0^{\varepsilon}-E_0^*\|_{L^2}^2+\|B_0^{\varepsilon}-B^e\|_{L^2}^2\big)+C \varepsilon^2.\
\end{align}
In light of \eqref{rhoqEBq}, \eqref{widetildeL2EM} and the fact that $(\tilde{\rho}_a^\varepsilon,\tilde{q}_a^\varepsilon,\tilde{E}_a^\varepsilon,\tilde{B}_a^\varepsilon)=(\tilde{\rho}^\varepsilon,\tilde{q}^\varepsilon,\tilde{E}^\varepsilon,\tilde{B}^\varepsilon)+\varepsilon(\rho_1,q_1,E_1,B_1)$, \eqref{F3.14} holds.  \hfill $\Box$

\subsection{High-order error estimates}
\begin{lemma}
\label{L5.5}
 It holds
\begin{align}
&\sup_{t\in\mathbb{R}^+}\big(\|\nabla \tilde{\rho}^\varepsilon(t)\|_{H^{m-2}}^2+\varepsilon^2\|\nabla\tilde{q}^\varepsilon(t)\|_{H^{m-2}}^2+\|\nabla\tilde{E}^\varepsilon(t)\|_{H^{m-2}}^2+\|\nabla\tilde{B}^\varepsilon(t)\|_{H^{m-2}}^2\big)\nonumber\\
&\quad\quad+\int_0^\infty \big(\|\nabla \tilde{\rho}^\varepsilon(t)\|_{H^{m-1}}^2+\|\nabla\tilde{q}^\varepsilon(t)\|_{H^{m-2}}^2 +\|\nabla \tilde{E}^\varepsilon(t)\|_{H^{m-2}}^2\big)\,dt\nonumber\\
&\quad \le C\big(\|\rho_0^\varepsilon-\rho_0^*\|_{H^{m-1}}^2+\|E_0^\varepsilon-E_0^*\|_{H^{m-1}}^2+\|B_0^\varepsilon-B^e\|_{H^{m-1}}^2\big)+C\varepsilon^2,\label{F5.5}
\end{align}
where $C>0$ is a constant independent of $\varepsilon$.
\end{lemma}

\noindent {\bf Proof.}
Before deriving \eqref{F5.5}, we need to establish estimates for $(\tilde{\rho}^\varepsilon_a,\tilde{q}^\varepsilon_a, \tilde{E}^\varepsilon_a, \tilde{B}^\varepsilon_a)$.  Let $\al\in\N^d$ with $ 1\leq |\al|\leq m-1$. Applying $\partial^{\alpha}$ to \eqref{F3.11ca} leads to
\begin{equation}
\label{551}
\begin{cases}
   \partial_t\partial^{\alpha}\tilde{\rho}^\varepsilon_a+\dive \partial^{\alpha}\tilde{q}^\varepsilon_a=\partial^{\alpha} J_1^\varepsilon,\\
\varepsilon^2\partial_t\partial^{\alpha}\tilde{q}^\varepsilon_a+\partial^{\alpha}\tilde{q}^\varepsilon_a+p'(1)\nabla \partial^{\alpha}\tilde{\rho}^\varepsilon_a
+\partial^{\alpha}\tilde{E}_a^\varepsilon=\partial^{\alpha} J_0^\varepsilon+ \partial^{\alpha} J_2^\varepsilon+\partial^{\alpha} J_3^\varepsilon+\partial^{\alpha} J^\varepsilon_4,\\[1mm]
   \partial_t \partial^{\alpha}\tilde{E}^\varepsilon_a-\dfrac{1}{\varepsilon}\na\times \partial^{\alpha}\tilde{B}^\varepsilon_a
=\partial^{\alpha} \tilde{q}^\varepsilon_a+\partial^{\alpha} J^\varepsilon_5,\\[3mm]
   \partial_t \partial^{\alpha} \tilde{B}^\varepsilon_a+\dfrac{1}{\varepsilon}\na\times \partial^{\alpha}\tilde{E}^\varepsilon_a= \partial^{\alpha} J_6^\varepsilon,\\
   \dive \partial^{\alpha} \tilde{E}_a^\varepsilon=-\partial^{\alpha}\tilde{\rho}^\varepsilon_a,\qquad \dive \partial^{\alpha} \tilde{B}^\varepsilon_a=0.
\end{cases}
\end{equation}
Let $\nu>0$ be chosen later. Carrying out the $L^2$ energy estimate for \eqref{551} as in \eqref{Fmmmmm}-\eqref{5.26}, we obtain
\begin{align}
 &\frac{1}{2}\frac{d}{dt}\Big(p'(1)\|\partial^{\alpha}\tilde{\rho}_a^\varepsilon\|_{L^2}^2
+\varepsilon^2\|\partial^{\alpha}\tilde{q}_a^\varepsilon\|_{L^2}^2+\|\partial^{\alpha}\tilde{E}^\varepsilon_a\|_{L^2}^2+\|\partial^{\alpha}\tilde{B}^\varepsilon_a\|_{L^2}^2 \Big)
+\|\partial^{\alpha}\tilde{q}_a^\varepsilon\|_{L^2}^2\nonumber\\
\leq & \frac{1}{2}\|\partial^{\alpha}\tilde{q}_a^\varepsilon\|_{L^2}^2+C\big(\|\partial^{\alpha} J_0^\varepsilon\|_{L^2}^2+\|\partial^{\alpha} J^\varepsilon_2\|_{L^2}^2+\|\partial^{\alpha} J^\varepsilon_3\|_{L^2}^2\big)\nonumber\\
&\quad+p'(1)\|\partial^{\alpha} J^\varepsilon_1\|_{L^2}\|\partial^{\alpha}\tilde{\rho}_a^\varepsilon\|_{L^2}+\nu\|\partial^{\alpha} \tilde{E}^\varepsilon_a\|_{L^2}^2+C\nu^{-1}\|\partial^{\alpha} J^\varepsilon_5\|_{L^2}^2\nonumber\\
&\quad+\|\partial^{\alpha} J^\varepsilon_6\|_{L^2}\|\partial^{\alpha} \tilde{B}^\varepsilon_a\|_{L^2}+\langle \partial^{\alpha} J_4^\varepsilon, \partial^{\alpha}\tilde{q}^\varepsilon_a\rangle.\label{5.41}
\end{align}
Here, for any $1\leq |\alpha|\leq m-1$, the combination of \eqref{F1.1311}, \eqref{J0} and Sobolev's inequality
\[ \|f\|_{L^{\infty}}\leq C \|\nabla f\|_{L^2}^{\frac{1}{2}}\|f\|_{\dot{H}^2}^{\frac{1}{2}}\leq C \|\nabla f\|_{H^{m-2}}   \]
yields
\begin{align}
\|\partial^{\alpha} J_0^\varepsilon\|_{L^2}^2&\leq C\big(\|\nabla^2 \rho^\varepsilon_a\|_{H^{m-2}}^2+\|\nabla \rho^*\|_{H^{m-2}}^2+\|\nabla E^\varepsilon\|_{H^{m-2}}^2\big) \big(\|\nabla \tilde{\rho}_a^\varepsilon\|_{H^{m-2}}^2+\|\nabla \tilde{E}^\varepsilon_a\|_{H^{m-2}}^2\big)\nonumber\\
&\quad+C\|\nabla \rho^\varepsilon\|_{H^{m-2}}^2\|\nabla^2\tilde{\rho}^\varepsilon_a\|_{H^{m-2}}^2.\label{J0high}
\end{align}
Together with the regularity estimates for $(\rho^\varepsilon,E^\varepsilon,B^\varepsilon)$, $(\rho^*,E^*,B^*)$ and $(\rho_1,E_1,B_1)$, this implies
\begin{align}\label{5.42}
&\quad\int_0^t \|\partial^{\alpha} J_0^\varepsilon(t')\|_{L^2}^2\, dt'\nonumber\\
&\leq C\sup_{t'\in[0,t]}\big(\|\nabla \tilde{\rho}_a^\varepsilon(t')\|_{H^{m-2}}^2+\|\nabla \tilde{E}^\varepsilon_a(t')\|_{H^{m-2}}^2\big)\nonumber\\
&\quad\times \int_0^t \big(\|\nabla^2 \rho^\varepsilon(t')\|_{H^{m-2}}^2+\|\nabla \rho^*(t')\|_{H^{m-2}}^2+\|\nabla E^\varepsilon(t')\|_{H^{m-2}}^2\big)\,dt'\nonumber\\
&\quad+C\sup_{t'\in[0,t]}\|\nabla \rho^\varepsilon(t')\|_{H^{m-2}}^2 \int_0^t\|\nabla^2\tilde{\rho}^\varepsilon_a(t')\|_{H^{m-2}}^2\,dt'\nonumber\\
&\leq C(\delta_1+\delta_1^*+\varepsilon)\sup_{t'\in[0,t]}\big(\|\nabla \tilde{\rho}_a^\varepsilon(t')\|_{H^{m-2}}^2+\|\nabla \tilde{E}^\varepsilon_a(t')\|_{H^{m-2}}^2\big)\nonumber\\
&\quad+C\delta_1 \int_0^t\|\nabla^2\tilde{\rho}^\varepsilon_a(t')\|_{H^{m-2}}^2\,dt'.
\end{align}
We shall use the faster rate for $J_4^\varepsilon$ and overcome its low regularity. The term $\langle \partial^{\alpha} J_4^\varepsilon, \partial^{\alpha}\tilde{q}^\varepsilon_a\rangle$ can be addressed as
\begin{equation}\label{high2euler}
\begin{aligned}
\langle  \partial^{\alpha} J^\varepsilon_4, \partial^{\alpha}\tilde{q}^\varepsilon_a\rangle&=\langle \Lambda^{-1} \partial^{\alpha} J^\varepsilon_4, \Lambda\partial^{\alpha}\tilde{q}^\varepsilon_a\rangle\leq C\|J^\varepsilon_4\|_{H^{m-2}}\|\tilde{ q}_a^\varepsilon\|_{H^{m}}.
\end{aligned}
\end{equation}
Moreover, for $1\leq |\alpha|\leq m-1$ one has
\begin{equation}\label{6455}
\begin{aligned}
&\|\partial^{\alpha} \tilde{\rho}^\varepsilon_a|_{t=0}\|_{L^2}^2+\|\partial^{\alpha}\tilde{q}^\varepsilon_a|_{t=0}\|_{L^2}^2+\|\partial^{\alpha} \tilde{E}^\varepsilon_a|_{t=0}\|_{L^2}^2+\|\partial^{\alpha}\tilde{B}^\varepsilon_a|_{t=0}\|_{L^2}^2\\
&\leq \|\rho_0^\varepsilon-\rho^*\|_{H^{m-1}}^2+\|E_0^\varepsilon-E_0^*\|_{H^{m-1}}^2+\|B_0^\varepsilon-B^e\|_{H^{m-1}}^2+C\varepsilon^2.
\end{aligned}
\end{equation}
By virtue of the uniform regularity estimates for $(\rho^\varepsilon,E^\varepsilon,B^\varepsilon)$, $(\rho^*,E^*,B^*)$ and $(\rho_1,E_1,B_1)$, it holds that
\begin{align}
&\quad\sup_{t\in\mathbb{R}^+}\|\tilde{\rho}_a^\varepsilon(t)\|_{H^{m}}^2+\int_0^{\infty}\big(\|\tilde{q}_a^\varepsilon(t)\|_{H^{m}}^2+\|\nabla\tilde{B}_a^\varepsilon(t)\|_{H^{m-2}}^2+\|\nabla\tilde{\rho}_a^\varepsilon(t)\|_{H^m}^2\big)\,dt\nonumber\\
&\leq C(\delta_1+\delta^*_1+\varepsilon).\label{5.29}
\end{align}
Integrating \eqref{5.41} in time and using the above estimates \eqref{J0high}-\eqref{6455}, we arrive at
\begin{align*}
&\quad \sup_{t\in\mathbb{R}^+}\big(\|\nabla\tilde{\rho}_a^\varepsilon(t)\|_{H^{m-2}}^2
+\varepsilon^2\| \nabla \tilde{q}_a^\varepsilon(t)\|_{H^{m-2}}^2+\|\nabla \tilde{E}^\varepsilon_a(t)\|_{H^{m-2}}^2+\| \nabla\tilde{B}^\varepsilon_a(t)\|_{H^{m-2}}^2\big)\\
&\quad\quad+\int_0^{\infty}\|\nabla\tilde{q}^\varepsilon_a(t)\|_{H^{m-2}}^2\,dt\\
&\leq C\big(\|\rho_0^\varepsilon-\rho_0^*\|_{H^{m-1}}^2+\|E_0^\varepsilon-E_0^*\|_{H^{m-1}}^2+\|B_0^\varepsilon-B^e\|_{H^{m-1}}^2\big)+C\varepsilon^2\\
&\quad+C(\delta_1+\delta_1^*+\varepsilon)\sup_{t\in\mathbb{R}^+}\big(\|\nabla \tilde{\rho}_a^\varepsilon(t)\|_{H^{m-2}}^2+\|\nabla \tilde{E}^\varepsilon_a(t)\|_{H^{m-2}}^2\big)\\
&\quad+C(\delta_1+\nu) \int_0^\infty\big(\|\nabla^2\tilde{\rho}^\varepsilon_a(t)\|_{H^{m-2}}^2+\|\nabla\tilde{E}^\varepsilon_a(t)\|_{H^{m-1}}^2\big)\,dt\\
&\quad+C\int_0^{\infty}\big(\|J_2^\varepsilon(t)\|_{H^{m-1}}^2+\|J_3^\varepsilon(t)\|_{H^{m-1}}^2+\nu^{-1}\|J_5^\varepsilon(t)\|_{H^{m-1}}^2\big)\,dt\\
&\quad+C\Big(\int_0^{\infty}\|J_4^\varepsilon(t)\|_{H^{m-2}}^2\,dt\Big)^{1/2} \Big(\int_0^\infty \|\tilde{q}_a^\varepsilon(t)\|_{H^m}^2\,dt\Big)^{1/2}\\
&\quad+C\int_0^{\infty}\big(\|J^\varepsilon_1(t)\|_{H^{m-1}}+\|J_6^\varepsilon(t)\|_{H^{m-1}}\big)\,dt\\
&\quad\times \sup_{t\in\mathbb{R}^+}\big(\|\nabla\tilde{\rho}_a^\varepsilon(t)\|_{H^{m-2}}^2
+\|\nabla\tilde{E}^\varepsilon_a(t)\|_{H^{m-2}}^2+\|\nabla\tilde{B}^\varepsilon_a(t)\|_{H^{m-2}}^2\big),
\end{align*}
from which as well as \eqref{J1}-\eqref{J5J6} and \eqref{5.29} we infer
\begin{align}\label{6.49}
&\quad \sup_{t\in\mathbb{R}^+}\big(\|\nabla\tilde{\rho}_a^\varepsilon(t)\|_{H^{m-2}}^2
+\varepsilon^2\| \nabla \tilde{q}_a^\varepsilon(t)\|_{H^{m-2}}^2+\|\nabla \tilde{E}^\varepsilon_a(t)\|_{H^{m-2}}^2+\| \nabla\tilde{B}^\varepsilon_a(t)\|_{H^{m-2}}^2\big)\nonumber\\
&\quad\quad+\int_0^{\infty}\|\nabla\tilde{q}^\varepsilon_a(t)\|_{H^{m-2}}^2\,dt\nonumber\\
&\leq C\big(\|\rho_0^\varepsilon-\rho_0^*\|_{H^{m-1}}^2+\|E_0^\varepsilon-E_0^*\|_{H^{m-1}}^2+\|B_0^\varepsilon-B^e\|_{H^{m-1}}^2\big)+C(1+\nu^{-1})\varepsilon^2\nonumber\\
&\quad+C(\delta_1+\nu) \int_0^\infty\big(\|\nabla^2\tilde{\rho}^\varepsilon_a(t)\|_{H^{m-2}}^2+\|\nabla\tilde{E}^\varepsilon_a(t)\|_{H^{m-1}}^2\big)\,dt.
\end{align}

\vspace{3mm}
We now establish the dissipation estimates for $\tilde{\rho}_a^\varepsilon$. To this end, we rewrite $\eqref{551}_2$ as
\begin{align}
\partial^{\alpha}\tilde{q}^\varepsilon_a+p'(1)\nabla \partial^{\alpha}\tilde{\rho}^\varepsilon_a +\partial^{\alpha}\tilde{E}^\varepsilon_a=\partial^{\alpha} J_0^\varepsilon+\partial^{\alpha} J_2^\varepsilon+\partial^{\alpha} J_3^\varepsilon-\varepsilon^2\partial_t \partial^{\alpha} q^\varepsilon. \label{F3.11carrr}
\end{align}
After taking the $L^2$-inner product of $\eqref{F3.11carrr}$ by $\nabla\partial^{\alpha} \tilde{\rho}^\varepsilon_a$ and using
\[  \partial_t\tilde{\rho}_a^\varepsilon=-\dive(q^\varepsilon-q^*-\varepsilon q_1)\quad
\mbox{and} \quad\partial^{\alpha}\tilde{\rho}^\varepsilon_a=-\dive \partial^{\alpha}\tilde{E}^\varepsilon_a, \]
we obtain
\begin{align*}
 &\quad p'(1)\|\nabla\partial^{\alpha} \tilde{\rho}^\varepsilon_a\|_{L^2}^2+\|\partial^{\alpha} \tilde{\rho}^\varepsilon_a\|_{L^2}^2\\
&=-\varepsilon^2\frac{d}{dt}\langle \partial^{\alpha} q^\varepsilon, \nabla \partial^{\alpha}\tilde{\rho}_a^\varepsilon\rangle+\varepsilon^2 \langle \partial^{\alpha} q^\varepsilon,\partial_t \nabla \partial^{\alpha}\tilde{\rho}_a^\varepsilon\rangle \\
&\quad +\langle \partial^{\alpha} J^\varepsilon_2+\partial^{\alpha} J^\varepsilon_3+\partial^{\alpha} J^\varepsilon_4,\nabla\partial^{\alpha} \tilde{\rho}^\varepsilon_a\rangle-\langle \partial^{\alpha} q^\varepsilon_a,\nabla \partial^{\alpha}\tilde{\rho}_a^\varepsilon\rangle\\
&\leq -\varepsilon^2\frac{d}{dt}\langle \partial^{\alpha} q^\varepsilon, \nabla \partial^{\alpha} \tilde{\rho}_a^\varepsilon\rangle +C\varepsilon^2\big( \|\dive \partial^{\alpha} q^\varepsilon\|_{L^2}^2+\|\dive \partial^{\alpha} q^*\|_{L^2}^2+\varepsilon^2\|\dive \partial^{\alpha} q_1\|_{L^2}^2\big)\\
&\quad+C\|\partial^{\alpha} \tilde{q}_a^\varepsilon\|_{L^2}^2+C\|\partial^{\alpha} J_0^\varepsilon\|_{L^2}^2+C\|\partial^{\alpha} J_2^\varepsilon\|_{L^2}^2+\|\partial^{\alpha} J_3^\varepsilon\|_{L^2}^2\\
&\quad+\frac{p'(1)}{2}\|\nabla \partial^{\alpha}\tilde{\rho}^\varepsilon_a\|_{L^2}^2.
\end{align*}
After integrating this in time, we discover that
\begin{align}
&\quad\int_0^{t}\Big(\frac{1}{2}p'(1)\|\nabla\partial^{\alpha}\tilde{\rho}_a^\varepsilon(t')\|_{L^2}^2+\|\partial^{\alpha}\tilde{\rho}_a^\varepsilon(t')\|_{L^2}^2\Big)\,dt'\nonumber\\
&\leq -\varepsilon^2\langle \partial^{\alpha} q^\varepsilon, \nabla \partial^{\alpha}\tilde{\rho}_a^\varepsilon\rangle \Big|^{t}_{0}+C\int_0^t\|\partial^{\alpha} \tilde{q}_a^\varepsilon(t')\|_{L^2}^2\,dt'\nonumber\\
&\quad+C\varepsilon^2\int_0^t \big(\|q^\varepsilon(t')\|_{H^m}^2+\|q^*(t')\|_{H^m}^2+\varepsilon^2\|q_1(t')\|_{H^m}^2\big)\,dt'\nonumber\\
&\quad+C\int_0^t\big(\| J_0^\varepsilon(t')\|_{H^{m-1}}^2+\| J_2^\varepsilon(t')\|_{H^{m-1}}^2+\| J_3^\varepsilon(t')\|_{H^{m-1}}^2\big)\,dt' .\label{6.47}
\end{align}
By $\rho_1|_{t=0}=0$, \eqref{5.2}, \eqref{5.29}, $\varepsilon \|q_0^\varepsilon\|_{H^{m}}+\varepsilon \|q^\varepsilon\|_{H^{m}}\leq C$ and Young's inequality, we have
\begin{equation}
\begin{aligned}\label{5.30r}
-\varepsilon^2\langle \partial^{\alpha} q^\varepsilon, \nabla \partial^{\alpha} \tilde{\rho}_a^\varepsilon\rangle \Big|^{t}_{0}&\leq \varepsilon \big(\varepsilon\|q^\varepsilon\|_{H^{m}}\big)\|\nabla \tilde{\rho}_a^\varepsilon\|_{H^{m-2}}+\varepsilon \big( \varepsilon \| q_0^\varepsilon\|_{H^{m}}\big)\|\nabla(\rho_0^\varepsilon-\rho_0^*)\|_{H^{m-2}}\\
&\leq \nu\|\nabla \tilde{\rho}_a^\varepsilon\|_{H^{m-2}}^2+C\|\nabla(\rho_0^\varepsilon-\rho_0^*)\|_{H^{m-2}}^2+C(1+\nu^{-1})\varepsilon^2.
\end{aligned}
\end{equation}
Putting \eqref{J1}-\eqref{J3}, \eqref{5.42}, \eqref{5.29} and \eqref{5.30r} into \eqref{6.47} and using Young's inequality, we derive
\begin{align}\label{5.32r}
&\quad \int_0^{\infty} \|\nabla\tilde{\rho}_a^\varepsilon(t)\|_{H^{m-1}}^2\,dt\nonumber\\
&\leq C\|\nabla(\rho_0^\varepsilon-\rho_0^*)\|_{H^{m-1}}^2+C(1+\nu^{-1})\varepsilon^2\nonumber\\
&\quad+C(\delta_1+\delta_1^*+\varepsilon+\nu)\sup_{t\in\mathbb{R}^+}\big(\|\nabla \tilde{\rho}_a^\varepsilon(t)\|_{H^{m-2}}^2+\|\nabla \tilde{E}^\varepsilon_a(t)\|_{H^{m-2}}^2\big)\nonumber\\
&\quad+C\int_0^{\infty}\|\tilde{q}^\varepsilon_a(t)\|_{H^m}^2\,dt+C\delta_1 \int_0^\infty\|\nabla^2\tilde{\rho}^\varepsilon_a(t)\|_{H^{m-2}}^2\,dt.
\end{align}

The next step is to derive the dissipation estimate for the error $\tilde{E}^\varepsilon_a$. Taking the $L^2$-inner product for $\eqref{F3.11carrr}$ by $\partial^{\alpha}\tilde{E}^\varepsilon_a$ with $1\leq |\alpha|\leq m-1$ and using $\partial^{\alpha}\tilde{\rho}^\varepsilon_a=-\partial^{\alpha}\dive \tilde{E}^\varepsilon_a$ yields
\begin{equation*}
\begin{aligned}
 &\quad \|\partial^{\alpha} \tilde{E}^\varepsilon_a\|_{L^2}^2+p'(1)\|\partial^{\alpha}\dive \tilde{E}^\varepsilon_a\|_{L^2}^2\\
&=-\varepsilon^2\frac{d}{dt}\langle \partial^{\alpha} q^\varepsilon,  \partial^{\alpha} \tilde{E}_a^\varepsilon\rangle+\varepsilon^2\langle \partial^{\alpha} q^\varepsilon,  \partial_t \partial^{\alpha} \tilde{E}_a^\varepsilon\rangle-\langle \partial^{\alpha} \tilde{q}_a^\varepsilon, \partial^{\alpha} \tilde{E}^\varepsilon_a\rangle\\
&\quad+\langle \partial^{\alpha} J^\varepsilon_0+\partial^{\alpha} J^\varepsilon_2+\partial^{\alpha} J^\varepsilon_3,\partial^{\alpha} \tilde{E}_a^\varepsilon\rangle.
\end{aligned}
\end{equation*}
Note that
\begin{align*}
    \partial_t \partial^{\alpha} \tilde{E}_a^\varepsilon=\frac{1}{\varepsilon}\nabla\times  \partial^{\alpha} \tilde{B}_a^\varepsilon+ \partial^{\alpha}\tilde{q}_a^\varepsilon+J^\varepsilon_5,
\end{align*}
so that
\begin{align*}
\varepsilon^2\langle \partial^{\alpha} q^\varepsilon,  \partial_t \partial^{\alpha} \tilde{E}_a^\varepsilon\rangle&=\varepsilon \langle \nabla\times  \partial^{\alpha}\tilde{q}^\varepsilon, \partial^{\alpha} \tilde{B}_a^\varepsilon\rangle+\varepsilon^2\langle \partial^{\alpha} q^\varepsilon, \partial^{\alpha}\tilde{q}_a^\varepsilon+J^\varepsilon_5\rangle.
\end{align*}
Therefore, for some constant $\nu>0$ to be determined, we obtain
\begin{equation*}
\begin{aligned}
 &\quad \|\partial^{\alpha} \tilde{E}^\varepsilon_a\|_{L^2}^2+p'(1)\|\partial^{\alpha}\dive \tilde{E}^\varepsilon_a\|_{L^2}^2\\
&=-\varepsilon^2\frac{d}{dt}\langle \partial^{\alpha}\tilde{q}^\varepsilon,  \partial^{\alpha} \tilde{E}_a^\varepsilon\rangle+\frac{1}{2 }\|\partial^{\alpha} \tilde{E}^\varepsilon_a\|_{L^2}^2+\nu \|\partial^{\alpha} \tilde{B}_a^\varepsilon\|_{L^2}^2+C\|\partial^{\alpha} \tilde{q}^\varepsilon_a\|_{L^2}^2\\
&\quad+C\varepsilon^2\big(\|\partial^{\alpha} q^\varepsilon\|_{L^2}^2+\nu^{-1}\|\nabla\times\partial^{\alpha} q^\varepsilon\|_{L^2}^2\big)\\
&\quad+C\| \partial^{\alpha} J^\varepsilon_0\|_{L^2}^2+C\|\partial^{\alpha} J^\varepsilon_2\|_{L^2}^2+C\|\partial^{\alpha} J^\varepsilon_3\|_{L^2}^2+C\|\partial^{\alpha} J^\varepsilon_5\|_{L^2}^2.
\end{aligned}
\end{equation*}
Thus, it follows that
\begin{equation*}
\begin{aligned}
&\quad\frac{1}{2}\int_0^t\|\partial^{\alpha} \tilde{E}^\varepsilon_a(t')\|_{L^2}^2\,dt'\\
&\leq -\varepsilon^2\langle \partial^{\alpha}\tilde{q}^\varepsilon,\partial^{\alpha}\tilde{E}^\varepsilon_a\rangle \Big|^t_0+C\int_0^t\|\nabla \tilde{q}^\varepsilon_a(t')\|_{H^{m-2}}^2\,dt'+\nu\int_0^t \|\nabla\tilde{B}^\varepsilon_a(t') \|_{H^{m-2}}^2\,dt'\\
&\quad+C(1+\nu^{-1})\varepsilon^2\int_0^t \|q^\varepsilon(t')\|_{H^{m}}^2\,dt'+C\int_0^t\|J_0^\varepsilon(t')\|_{H^{m-1}}^2\,dt'\\
&\quad+C\int_0^t\big(\|J_2^\varepsilon(t')\|_{H^{m-1}}^2+\|J_3^\varepsilon(t')\|_{H^{m-1}}^2+\|J_4^\varepsilon(t')\|_{L^2}^2+\|J_5(t')\|_{H^{m-1}}^2\big)\,dt',
\end{aligned}
\end{equation*}
which, together with \eqref{J3}-\eqref{J5J6}, \eqref{5.42} and $\varepsilon\sup_{t\in\mathbb{R}^{+}}\|q^\varepsilon(t)\|_{H^{m-1}}\leq C$, leads to
\begin{align}\label{5.35r}
&\quad\int_0^\infty \|\nabla\tilde{E}^\varepsilon_a(t)\|_{H^{m-2}}^2\,dt\nonumber\\
&\leq C\|\nabla(E_0^\varepsilon-E_0)\|_{H^{m-2}}^2+C\varepsilon^2\nonumber\\
&\quad+C(\delta_1+\delta_1^*+\varepsilon)\sup_{t\in\mathbb{R}^{+}}\big(\|\nabla \tilde{\rho}_a^\varepsilon(t)\|_{H^{m-2}}^2+\|\nabla \tilde{E}^\varepsilon_a(t)\|_{H^{m-2}}^2\big)+C\delta_1 \int_0^\infty\|\nabla^2\tilde{\rho}^\varepsilon_a(t)\|_{H^{m-2}}^2\,dt\nonumber\\
&\quad+C\int_0^{\infty}\|\nabla\tilde{q}^\varepsilon_a(t)\|_{H^{m-2}}^2\,dt+\nu\int_0^\infty \|\nabla\tilde{B}^\varepsilon_a(t) \|_{H^{m-2}}^2\,dt.
\end{align}

We are in a position to obtain the dissipation estimate for $\tilde{B}^\varepsilon_a$ required in \eqref{5.35r}.
For $1\leq |{\al'}|\leq m-2$, we deduce from $\eqref{551}_3$-$\eqref{551}_4$ that
\begin{align*}
\|\nabla\times \pa^{\al'} \tilde{B}^\varepsilon_a\|_{L^2}^2&=\varepsilon\langle \partial_t\pa^{\al'} \tilde{E}^\varepsilon_a,  \nabla\times \pa^{\al'} \tilde{B}^\varepsilon_a\rangle-\langle \pa^{\al'}\tilde{q}^\varepsilon_a, \nabla\times \pa^{\al'} \tilde{B}^\varepsilon_a\rangle-\langle \pa^{\al'} J^\varepsilon_5,\nabla\times \pa^{\al'} \tilde{B}^\varepsilon_a\rangle\nonumber\\
&=\varepsilon\frac{d}{dt} \langle \pa^{\al'} \tilde{E}^\varepsilon_a,  \nabla\times \pa^{\al'} \tilde{B}^\varepsilon_a\rangle-\varepsilon\langle \pa^{\al'}\tilde{q}^\varepsilon_a, \nabla\times \pa^{\al'} \tilde{B}^\varepsilon_a\rangle\nonumber\\
&\quad-\langle \pa^{\al'} J^\varepsilon_5,\nabla\times \pa^{\al'} \tilde{B}^\varepsilon_a\rangle+\| \nabla\times \pa^{\al'} \tilde{E}^\varepsilon_a\|_{L^2}^2-\varepsilon\langle \pa^{\al'}\nabla\times \tilde{E}^\varepsilon_a, \pa^{\al'} J^\varepsilon_6\rangle.
\end{align*}
Consequently, one discovers that
\begin{align}
&\quad \frac{1}{2}\int_0^t \|\nabla\times \pa^{\al'} \tilde{B}^\varepsilon_a(t')\|_{L^2}^2\,dt'\nonumber\\
&\leq \varepsilon \langle \pa^{\al'} \tilde{E}^\varepsilon_a,  \nabla\times \pa^{\al'} \tilde{B}^\varepsilon_a\rangle\Big|^t_0+C\int_0^t\big(\varepsilon^2 \|\tilde{q}^\varepsilon_a(t')\|_{H^{m-2}}^2+\|J^\varepsilon_5(t')\|_{H^{m-2}}^2+\|J^\varepsilon_6(t')\|_{H^{m-2}}^2\big)\,dt'\nonumber\\
&\quad+C\int_0^t \|\nabla \tilde{E}^\varepsilon_a(t')\|_{H^{m-1}}^2\,dt'.\label{657}
\end{align}
By \eqref{J5J6}, \eqref{5.29}, \eqref{657} and $\dive\tilde{B}^\varepsilon_a=0$, it yields
\begin{align}\label{Hhigh}
\int_0^{\infty}\|\nabla \tilde{B}^\varepsilon_a(t)\|_{H^{m-2}}^2\,dt&\leq C\varepsilon \sup_{t\in\mathbb{R}^+}\big(\| \nabla  \tilde{E}^\varepsilon_a(t)\|_{H^{m-2}}^2+\|\nabla \tilde{B}^\varepsilon_a(t)\|_{H^{m-2}}^2\big)\nonumber\\
&\quad+C\varepsilon^2+C\int_0^\infty \|\nabla\tilde{E}^\varepsilon_a(t)\|_{H^{m-2}}^2\,dt.
\end{align}
Collecting \eqref{6.49}, $\eta\times \eqref{5.32r}$, $\eta\times\eqref{5.35r}$ and $\eta\times \eqref{Hhigh}$ with some small constant $\eta$, taking $\nu$ sufficiently small and using the smallness for $\delta_1, \delta_1^*$ and $\varepsilon\leq \varepsilon_2$ with some suitably small $\varepsilon_2$, we obtain
\begin{align}
&\quad \sup_{t\in\mathbb{R}^{+}}\big(\|\nabla\tilde{\rho}_a^\varepsilon(t)\|_{H^{m-2}}^2
+\varepsilon^2\|\nabla\tilde{q}_a^\varepsilon(t)\|_{H^{m-2}}^2+\|\nabla\tilde{E}_a^\varepsilon(t)\|_{H^{m-2}}^2+\|\nabla\tilde{B}_a^\varepsilon(t)\|_{H^{m-2}}^2\big)\nonumber
\\
&\quad\quad+\int_{0}^{\infty}\big(\|\nabla\tilde{\rho}_a^\varepsilon(t)\|_{H^{m-1}}^2+\|\nabla\tilde{q}_a^\varepsilon(t)\|_{H^{m-2}}^2+\|\nabla \tilde{E}_a^\varepsilon(t)\|_{H^{m-2}}^2+\|\nabla \tilde{B}_a^\varepsilon(t)\|_{H^{m-2}}^2\big)\,dt\nonumber\\
 &\le C\big(\|\rho_0^{\varepsilon}-\rho_0^*\|_{H^{m-1}}^2+\|E_0^{\varepsilon}-E_0^*\|_{H^{m-1}}^2+\|B_0^{\varepsilon}-B^e\|_{H^{m-1}}^2\big)+C \varepsilon^2.\nonumber
\end{align}
In the case $\varepsilon\geq \varepsilon_2$, it is clear that  \eqref{F5.5} holds thanks to the uniform estimates satisfied by $(\tilde{\rho}^\varepsilon,\tilde{q}^\varepsilon,\tilde{E}^\varepsilon,\tilde{B}^\varepsilon)$. Using that
\[  (\tilde{\rho}_a^\varepsilon,\tilde{q}_a^\varepsilon,\tilde{E}_a^\varepsilon,\tilde{B}_a^\varepsilon)=(\tilde{\rho}^\varepsilon,\tilde{q}^\varepsilon,\tilde{E}^\varepsilon,\tilde{B}^\varepsilon)+\varepsilon(\rho_1,q_1,E_1,B_1),
\]
the desired bounds for $(\tilde{\rho}^\varepsilon,\tilde{q}^\varepsilon,\tilde{E}^\varepsilon,\tilde{B}^\varepsilon)$ follow.
\hfill $\Box$





 \vspace{1cm}
	
	\noindent 

	\noindent \textbf{Conflict of interest.} The authors do not have any possible conflicts of interest.

	\vspace{3mm}

	\noindent \textbf{Data availability statement.}
	Data sharing is not applicable to this article as no data sets were generated or analyzed during the current study.

\vspace{3mm} 


\vspace{3mm} 

\noindent \textbf{Acknowledgments.}
The authors are grateful to the referee for his/her valuable comments, which 
improve the presentation of this manuscript. T. Crin-Barat is supported by the project ANR-24-CE40-3260 -- Hyperbolic Equations, Approximations \& Dynamics (HEAD). L.-Y. Shou is supported by National Natural Science
Foundation of China (Grant No. 12301275).

\bibliographystyle{abbrv}

\end{document}